\documentclass{article}
\usepackage{amsmath}[1996/11/01]
\usepackage{amssymb,amsthm,amsxtra}
\usepackage{epic,eepic}             
\usepackage{epsfig}

\def\[#1\]{\begin{equation}#1\end{equation}}
\makeatletter
\def\beq{%
   \relax\ifmmode
      \@badmath
   \else
      \ifvmode
         \nointerlineskip
         \makebox[.6\linewidth]%
      \fi
      $$
   \fi
}
\def\eeq{%
   \relax\ifmmode
      \ifinner
         \@badmath
      \else
         $$
      \fi
   \else
      \@badmath
   \fi
   \ignorespaces
}

\def\enddisplaymath{\eeq\global\@ignoretrue}
\makeatother

\newtheorem{thm}{Theorem}
\newtheorem{cor}[thm]{Corollary}
\newtheorem{lem}[thm]{Lemma}
\newtheorem{prop}[thm]{Proposition}
\theoremstyle{remark}
\newtheorem*{rem}{Remark}

\theoremstyle{definition}
\newtheorem{defn}{Definition}

\numberwithin{equation}{section}
\numberwithin{thm}{section}

\DeclareMathOperator{\Exp}{E}

\renewcommand{\Re}{\operatorname{Re}}

\DeclareMathOperator{\Prob}{Pr}
\DeclareMathOperator{\dist}{dist}
\DeclareMathOperator{\argum}{arg}

\DeclareMathOperator{\erf}{erf}

\newcommand{\wt}{\widetilde}

\newcommand{\C}{\mathbb C}
\newcommand{\R}{\mathbb R}

\newcommand\psymmU{%
\begin{picture}(1,1)(0,0)%
\allinethickness{0.5pt}%
\path(0,0)(0,1)(1,1)(1,0)(0,0)%
\end{picture}}
\newcommand\psymmUU{%
\begin{picture}(1,1)(0,0)%
\allinethickness{0.5pt}%
\path(0,0)(0,1)(1,1)(1,0)(0,0)%
\put(0.5,0.5){\makebox(0,0){$\cdot$}}%
\end{picture}}
\newcommand\psymmO{%
\begin{picture}(1,1)(0,0)%
\allinethickness{0.5pt}%
\path(0,0)(0,1)(1,1)(1,0)(0,0)%
\path(0,0)(1,1)%
\end{picture}}
\newcommand\psymmS{%
\begin{picture}(1,1)(0,0)%
\allinethickness{0.5pt}%
\path(0,0)(0,1)(1,1)(1,0)(0,0)%
\path(1,0)(0,1)%
\end{picture}}
\newcommand\psymmu{%
\begin{picture}(1,1)(0,0)%
\allinethickness{0.5pt}%
\path(0,0)(0,1)(1,1)(1,0)(0,0)%
\path(0,0)(1,1)%
\path(0,1)(1,0)%
\end{picture}}

\newbox\tsymmUbox
\newbox\tsymmUUbox
\newbox\tsymmObox
\newbox\tsymmSbox
\newbox\tsymmubox
\setbox\tsymmUbox =\hbox{\kern0.75pt\setlength{\unitlength}{6pt}\psymmU \kern0.75pt}
\setbox\tsymmUUbox=\hbox{\kern0.75pt\setlength{\unitlength}{6pt}\psymmUU\kern0.75pt}
\setbox\tsymmObox =\hbox{\kern0.75pt\setlength{\unitlength}{6pt}\psymmO \kern0.75pt}
\setbox\tsymmSbox =\hbox{\kern0.75pt\setlength{\unitlength}{6pt}\psymmS \kern0.75pt}
\setbox\tsymmubox =\hbox{\kern0.75pt\setlength{\unitlength}{6pt}\psymmu \kern0.75pt}

\def\tsymmO{{\copy\tsymmObox}}
\def\tsymmS{{\copy\tsymmSbox}}
\def\tsymmu{{\copy\tsymmubox}}
\newbox\symmUbox
\newbox\symmUUbox
\newbox\symmObox
\newbox\symmSbox
\newbox\symmubox
\setbox\symmUbox =\hbox{\kern0.75pt\setlength{\unitlength}{4.5pt}\psymmU \kern0.75pt}
\setbox\symmUUbox=\hbox{\kern0.75pt\setlength{\unitlength}{4.5pt}\psymmUU\kern0.75pt}
\setbox\symmObox =\hbox{\kern0.75pt\setlength{\unitlength}{4.5pt}\psymmO \kern0.75pt}
\setbox\symmSbox =\hbox{\kern0.75pt\setlength{\unitlength}{4.5pt}\psymmS \kern0.75pt}
\setbox\symmubox =\hbox{\kern0.75pt\setlength{\unitlength}{4.5pt}\psymmu \kern0.75pt}

\def\symmO{{\copy\symmObox}}
\def\symmS{{\copy\symmSbox}}
\def\symmu{{\copy\symmubox}}

\def\tsymmg{\circledast}
\def\symmg{\circledast}

\begin{document}

\title{{\bf The asymptotics of monotone subsequences of involutions}}
\author{ {\bf Jinho Baik}\footnote{
Princeton University and Institite for Advanced Study, New Jersey,
jbaik@math.princeton.edu}
\ \ and \ \ {\bf Eric M. Rains}\footnote{AT\&T Research, New Jersey,
rains@research.att.com}
}
\date{January 31, 2001}
\maketitle

\begin{abstract}
%
%

We compute the limiting distributions of the lengths of the longest
monotone subsequences of random (signed) involutions with or without
conditions on the number of fixed points (and negated points) as the sizes
of the involutions tend to infinity. The resulting distributions are,
depending on the number of fixed points, (1) the Tracy-Widom distributions
for the largest eigenvalues of random GOE, GUE, GSE matrices, (2) the
normal distribution, or (3) new classes of distributions which interpolate
between pairs of the Tracy-Widom distributions. We also consider the second
rows of the corresponding Young diagrams. In each case the convergence of
moments is also shown. The proof is based on the algebraic work of the
authors in \cite{PartI} which establishes a connection between the
statistics of random involutions and a family of orthogonal polynomials,
and an asymptotic analysis of the orthogonal polynomials which is obtained
by extending the Riemann-Hilbert analysis for the orthogonal polynomials by
Deift, Johansson and the first author in \cite{BDJ}.

\end{abstract}


%
%

\section{Introduction}

\subsection*{$\beta$-Plancherel measure}

\begin{figure}[ht]
 \centerline{\epsfig{file=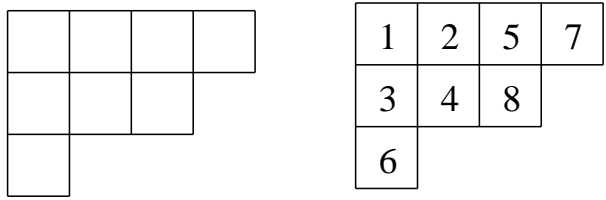, width=8cm}}
 \caption{Young diagram and standard Young tableau}
\label{fig-partition}
\end{figure}

 In the last few years, it has been observed by many authors that there are
certain connections between random permutations and/or Young tableaux, and
random matrices. One of the earliest clues to this relationship appeared in the
work of Regev \cite{Re} in 1981. A Young diagram, or equivalently a partition
$\lambda=(\lambda_1,\lambda_2,\cdots)\vdash n$ ($\lambda_1\ge \lambda_2\ge
\cdots$, $\sum \lambda_j=n$) is an array of $n$ boxes with top and left
adjusted as in the first picture of Figure \ref{fig-partition}, which
represents the example $\lambda=(4,3,1)\vdash 8$.
A standard Young tableau $Q$ is a filling of the diagram $\lambda$ by numbers
$1,2,\cdots, n$ such that numbers are increasing along each row and along each
column. In this case, we say that the tableau $Q$ has the shape $\lambda$. The
second picture in Figure \ref{fig-partition} is an example of a standard Young
tableau with shape $\lambda=(4,3,1)$. Let $d_\lambda$ denote the number of
standard Young tableaux of shape $\lambda$. A result of \cite{Re} is that for
fixed $\beta >0$ and fixed $l$, as $n\to\infty$,
\begin{equation}\label{int-as-1}
   \sum_{\substack{\lambda\vdash n \\ \lambda_1\le l}}
\bigl(d_\lambda\bigr)^\beta
\sim \biggl[ \frac{l^{l^2/2}l^n}{(\sqrt{2\pi})^{(l-1)/2}n^{(l-1)(l+2)/4}}
\biggr]^\beta \frac{n^{(l-1)/2}}{l!}
\int_{\R^l} e^{-\frac12\beta l\sum_{j}x_j^2}
\prod_{j<k} |x_j-x_k|^\beta d^l x.
\end{equation}
The multiple integral on the right hand side is called the Selberg integral,
which can be computed exactly for each $\beta$ in terms of the Gamma function.
In particular, when $\beta=1,2,4$, this integral is the normalization constant
of the eigenvalue density of a random matrix taken from the Gaussian orthogonal
ensemble (GOE), Gaussian unitary ensemble (GUE), Gaussian symplectic ensemble
(GSE), respectively (see e.g. \cite{Mehta}). Motivated by this result, we
define the \emph{$\beta$-Plancherel measure} on the set $Y_n$ of Young diagrams
(or partitions) of size $n$ by
\begin{equation}\label{int-as-2}
   M^\beta_n(\lambda) :=
\frac{d_\lambda^\beta}{\sum_{\mu\vdash n}d_\mu^\beta}, \qquad
\lambda\in Y_n.
\end{equation}
A natural question is the limiting statistics of a random $\lambda\in Y_n$
under $M^\beta_n$ as $n\to\infty$.

The case when $\beta=2$ is quite well studied. In this case, $M^2_n$ is
called the \emph{Plancherel measure} which arises in the representation
theory of the symmetric group $S_n$. Denote by $L_n^{(k)}$ the random
variable $\lambda_k$ under the Plancherel measure $M^2_n$, and set
$L_n=L_n^{(1)}$.
In 1977, the limiting expected shape of $\lambda$ under $M^2_n$ is
obtained in \cite{VK1}, also independently in \cite{LS} for the so-called
Poissonized Plancherel measure. In particular, it is shown that
\begin{equation}\label{qz1.3}
  \lim_{n\to\infty} \frac{\Exp(L_n)}{\sqrt{n}}=2.
\end{equation}
A central limit theorem for $L_n$ is then obtained by \cite{BDJ}:
\begin{equation}\label{qz1.4}
  \lim_{n\to\infty} \Prob\biggl(\frac{L_n-2\sqrt{n}}{n^{1/6}}\biggr)=
  F_2(x),
\end{equation}
where $F_2$ is the so-called Tracy-Widom distribution function, which is
expressed in terms of a solution to the Painlev\'e II equation (see
Definition \ref{def1} in Section \ref{sec:dist} for the definition). The
connection to random matrix theory comes from this function $F_2$: in
1994, Tracy and Widom \cite{TW1} proved that under proper centering and
scaling (which is different from the scaling for $L_n$ in \eqref{qz1.4}),
the largest eigenvalue of a random matrix taken from the Gaussian unitary
ensemble (GUE) has the same limiting distribution given by $F_2$. In other
words, after proper centering and scaling, the first row of a random Young
diagram under the Plancherel measure behaves statistically for large $n$
like the largest eigenvalue of a random GUE matrix. Then in the same paper
\cite{BDJ}, it was conjectured that $L_n^{(k)}$ of a random $\lambda\in
Y_n$ under $M^2_n$ have the same limiting distribution as the $k^{th}$
largest eigenvalue of a random GUE matrix for each $k$. This conjecture was
supported by numerical simulations of Odlyzko and the second author, and
was proved to be true for the second row $L_n^{(2)}$ in
\cite{BDJ2}. The full conjecture for the general row $L_n^{(k)}$ was
subsequently proved by \cite{Ok}, \cite{BOO} and \cite{kurtj:disc},
independently. The authors in \cite{Ok, BOO, kurtj:disc} proved the
convergence in joint distribution for general rows, and also in \cite{BOO,
kurtj:disc}, the authors obtained discrete sine kernel representations for
the so-called bulk scaling limit of correlation functions, an analogue of
the sine kernel that appears in the GUE matrix case. In \cite{BDJ} and
\cite{BDJ2}, in addition to convergence in distribution, the authors
also proved convergence of moments for $L_n^{(1)}$ and $L_n^{(2)}$,
respectively.  Convergence of (joint) moments for the general rows is
obtained recently by
\cite{BDR}. We also mention that there are many works on similar relationships
between tableaux/combinatorics and GUE random matrices. See for example,
\cite{TW3, Bo, kurtj:shape, TW:randomwords, ITW, kuperberg, stanley}. We refer
the readers to \cite{AD, OR, deiftnotice} for a survey and history of $L_n$,
and to \cite{Mehta} for general reference on random matrix theory (see also the
recent book \cite{deift}).

One of the main topics in this paper is the limiting statistics of $\lambda\in
Y_n$ under $M^1_n$. From \eqref{int-as-1} and the results for the case when
$\beta=2$, one might guess that for $\beta=1$ the limiting statistics of a
random $\lambda$ is same as the limiting statistics of the eigenvalues of GOE
(Gaussian orthogonal ensemble) matrices. We establish this fact for the first
two rows. More precisely, we prove (see Theorems \ref{thm22} and
\ref{thm-second2} together with the remark that follows), denoting by
$\tilde{L}_n^{(k)}$ the random variable $\lambda_k$ of a random $\lambda\in
Y_n$ under $M^1_n$,
\begin{equation}\label{qz1.5}
  \lim_{n\to\infty} \Prob\biggl( \frac{\tilde{L}_n^{(k)}-2\sqrt{n}}{n^{1/6}} \le x
  \biggr) = F_1^{(k)}(x), \qquad k=1,2,
\end{equation}
where $F_1^{(k)}(x)$ is the limiting distribution function \cite{TW2} for the
(scaled) $k^{th}$ largest eigenvalue of a random matrix taken from GOE. We also
prove convergence of moments. As in the case of $\beta=2$, we expect that
the above result should extend to the general rows $k\ge 3$, and also to the
joint distributions. For general values of $\beta>0$, again from
\eqref{int-as-1}, we expect that in the large $n$ limit, the rows of a random
Young diagram under $M^\beta_n$ correspond to the Coulomb charges on the real
line with the quadratic potential at the inverse temperature $\beta$, which
specializes to GOE, GUE, GSE eigenvalue distributions for the cases
$\beta=1,2,4$, respectively. This conjecture seems natural from the perspective
of the discrete Coulomb gas interpretation for the Plancherel measure case
$M^2_n$ given by Johansson \cite{kurtj:shape, kurtj:disc}.


\subsection*{Random involutions}

The Plancherel measure $M^2_n$ has a nice combinatorial interpretation. The
well-known Robinson-Schensted correspondence \cite{Sc} establishes a bijection
between the permutations $\pi$ of size $n$ and the pairs of standard Young
tableaux $(P,Q)$ where the shape of $P$ and the shape of $Q$ are the same and
the shape of $P$ (or $Q$), denoted by $\lambda(\pi)$, is a partition of $n$.
Thus the Plancherel measure $M^2_n$ on $Y_n$ is the push forward of the uniform
probability measure on $S_n$. Moreover, under this correspondence,
$\lambda_1(\pi)$ is equal to the length  of the longest increasing subsequence
of $\pi$. More generally, a theorem of Greene \cite{Greene} says that
$\lambda_1(\pi)+\cdots+\lambda_k(\pi)$ is equal to the length of the longest
so-called $k$-increasing subsequence of $\pi$. Thus the difference of the
lengths of the longest $k$-increasing subsequence and the longest
$(k-1)$-increasing subsequence of $\pi\in S_n$ under the uniform probability
measure is equal to $\lambda_k$ of $\lambda\in Y_n$ under the Plancherel
measure $M^2_n$ in the sense of joint distributions. Thus for example,
\eqref{qz1.3} and \eqref{qz1.4} can be restated for the results on the longest
increasing subsequence of a random permutation. On the other hand, the sum of
the lengths of the first $k$ columns of $\lambda$ is equal to the length of the
longest $k$-\emph{decreasing} subsequence of corresponding $\pi$. But since the
transpose $\lambda^t$ have the same statistics as $\lambda$ under $M^2_n$, the
results \eqref{qz1.3} and \eqref{qz1.4} also hold for the longest decreasing
subsequence of a random permutation.

The measure $M^1_n$ also has a combinatorial interpretation. If $\pi$ is
mapped to $(P,Q)$ under the Robinson-Schensted correspondence, then
$\pi^{-1}$ is mapped to $(Q,P)$ (see e.g. Section 5.1.4 of
\cite{Kn}). Therefore, the set of involutions $\pi=\pi^{-1}\in S_n$ is in
bijection with the set of standard Young tableaux whose shapes are
partitions of $n$. Consequently, the uniform probability measure on the set
of involutions
\begin{equation}
  \tilde{S}_n=\{\pi\in S_n : \pi=\pi^{-1}\}
\end{equation}
is pushed forward to the $1$-Plancherel measure $M^1_n$ on $Y_n$. Thus the
result \eqref{qz1.5} for $k=1$ implies that in the large $n$ limit, the length
of the longest increasing (also decreasing) subsequence of a random involution
behaves statistically like the largest eigenvalue of a random GOE matrix.

An involution $\pi\in \tilde{S}_n$ consists only of 1 cycles and 2 cycles. It
turns out that if we put a condition on the number of 1 cycles (or fixed
points) of $\pi$, the limiting distribution is different. Introduce a new
ensemble
\begin{equation}\label{qz1.6}
  S_{n,m} = \{ \pi\in \tilde{S}_{2n+m} : |\{x : \pi(x)=x\}|=m \}.
\end{equation}
For an involution $\pi$, the number of fixed points is equal to the number of
odd parts of $\lambda^t$ (see \cite{Kn}). Equivalently, the number of fixed
points of $\pi$ is equal to $\lambda_1-\lambda_2+\lambda_3-\cdots$. Thus the
uniform probability measure on the set $S_{n,m}$ is push-forwarded to the
measure
\begin{equation}\label{qz1.8}
  \frac{d_\lambda}{\sum_{\mu\in Y_{n,m}} d_\mu}, \qquad
  \lambda\in Y_{n,m},
\end{equation}
where
\begin{equation}\label{qz1.7}
  Y_{n,m}= \{ \lambda=(\lambda_1,\lambda_2,\cdots)\in Y_{2n+m} :
  \sum_j (-1)^{j-1} \lambda_j = m \}.
\end{equation}
Note that the rows and columns of $\lambda\in Y_{n,m}$ now have different
distributions. We denote by $L_{n,m}^{\symmO,(k)}$ and $L_{n,m}^{\symmS, (k)}$
the random variables given by the lengths of the $k^{th}$ \emph{row} and the
$k^{th}$ \emph{column} of a random $\lambda\in Y_{n,m}$ under the measure
\eqref{qz1.8}, respectively. We also set $L_{n,m}^{\symmO}=L^{\symmO,
(1)}_{n,m}$ and $L_{n,m}^{\symmS}=L^{\symmS, (1)}$, the length of the longest
\emph{increasing} and \emph{decreasing} subsequences of a random $\pi\in
S_{n,m}$ under the uniform probability measure.

Set
\begin{equation}\label{qz1.9}
  \alpha= \frac{m}{\sqrt{2n}}.
\end{equation}
The limiting distribution of $L_n^{\symmO}$ differs depending on $\alpha$.
Indeed we prove in Theorem \ref{newthm1} and Theorem \ref{thmbigalpha1} below
that
\begin{align}
\label{qz1.10}    \lim_{n\to\infty} \Prob\biggl(
    \frac{L^\symmO_{n,[\alpha\sqrt{2n}]}-2\sqrt{2n+m}}{(2n+m)^{1/6}}
\le x \biggr) &= F_4(x),
\qquad 0\le \alpha<1, \\
\label{qz1.11}    \lim_{n\to\infty} \Prob\biggl(
    \frac{L^\symmO_{n,[\sqrt{2n}]}-2\sqrt{2n+m}}{(2n+m)^{1/6}}
\le x \biggr) &= F_1(x),
\qquad \alpha=1, \\
\label{qz1.12}  \lim_{n\to\infty} \Prob\biggl(
\frac{L^\symmO_{n,[\alpha\sqrt{2n}]}
-(\alpha+1/\alpha)\sqrt{2n+m}}{\sqrt{(1/\alpha-1/\alpha^3)}(2n+m)^{1/4}}
\biggr) &=\erf(x), \qquad \alpha>1,
\end{align}
where $F_4$ and $F_1$ are the distributions for the limiting fluctuations of
the largest eigenvalues of random GSE and GOE matrices respectively, and $\erf$
is the standard normal distribution. Again we also prove convergence of
moments. We note that $F_4=F_1^{(2)}$; the limiting distributions of the
largest eigenvalue of GSE and the second largest eigenvalue of GOE are the same
(see the discussion at the end of Section \ref{sec:results} below.)

The role of the number of fixed points for the limiting distribution can be
seen from the following point selection picture. Consider a unit square
$[0,1]\times [0,1]$ in the plane and set $\delta=\{(x,x): 0\le x\le 1\}$, the
diagonal. Suppose we select $n$ points at random in the lower triangle $0\le
x<y\le 1$ and take the mirror image of the points about the diagonal $\delta$.
We also select $m$ points at random on the diagonal $\delta$. Hence there is a
total of $2n+m$ points in the square.
\begin{figure}[ht]
 \centerline{\epsfig{file=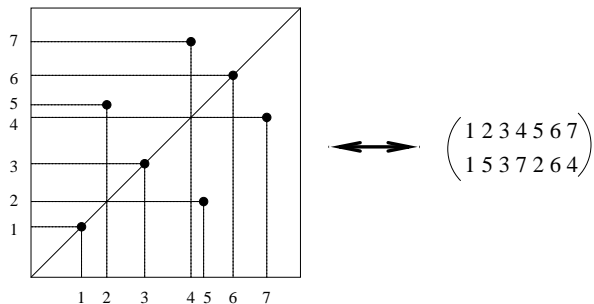, width=8cm}}
 \caption{Point selection process}
\label{fig-points}
\end{figure}
As illustrated in Figure \ref{fig-points}, one such choice of points gives
rise to a permutation $\pi$ satisfying $\pi^2=1$ with $m$ fixed points,
i.e., $\pi\in S_{n,m}$. The length $L_n^\symmO(\pi)$ of the longest
increasing subsequence of $\pi$ is then equal to the `length' of the
longest (piecewise linear) up/right path in the square from $(0,0)$ to
$(1,1)$, where the `length' of a path is defined by the number of points on
the path. The length of the longest up/right path in the above point
selection process has the same distribution as $L_n^\symmO$. Now note that
the points on $\delta$ form an increasing path.  When $m$ is large
comparing to $n$, there are many points on $\delta$ and we expect that the
longest path consists mostly of diagonal points. Hence we are in the linear
statistics situation, and thus the order of fluctuation of the length of
the longest path is expected to be $(\text{\it mean})^{1/2}$ by the usual central
limit theorem. On the other hand, when $m$ is small comparing to $n$, then
the longest path contains few diagonal points (none if $m=0$) and we are in
the situation of a $2$-dimensional maximization problem. In this case it has
been believed, and in a few cases (for example, \cite{BDJ, kurtj:disc,
GTW}) it is proved, that the fluctuation has order $(\text{\it mean})^{1/3}$. (For
random permutations, which have a similar interpretation as a point
selection process, one can see from the scaling in \eqref{qz1.4} that the
fluctuation is of order $(\text{\it mean})^{1/3}$.) Thus there must be a transition of
the limiting distribution as the size of $m$ varies. The results
\eqref{qz1.10}--\eqref{qz1.12} show that $\alpha=1$ is the transition
point. The fixed points play the role of adding a special line in the
$2$-dimensional maximization problem. (See \cite{BR4} for a relevant work
where two special lines are added to a $2$-dimensional maximization
problem.) We note that when $\alpha=1$, $L^\symmO_{n,[\sqrt{2n}]}$ in
\eqref{qz1.11} has the same limiting distribution as $\wt{L}^{(1)}_n$ (see
\eqref{qz1.5}). This is because the typical number of fixed points of a random
involution of size $k$ is $\sqrt{k}$. Indeed the result \eqref{qz1.5} is proved
by using \eqref{qz1.10}--\eqref{qz1.12} and taking a summation over the number
of fixed points (see Section \ref{as-asymp3} below) to which only $\alpha=1$
gives the main contribution.

Once the transition point $\alpha=1$ is known, it is of interest to investigate
the transition more carefully. We set
\begin{equation}\label{qz1.13}
  \alpha=1-\frac{2w}{(2n)^{1/6}},
\end{equation}
and take $n\to\infty$ while keeping $w$ fixed. We prove that (Theorem
\ref{newthm2} below) there is a one-parameter family of distribution functions
$F^\symmO(x;w)$, which is expressed in terms of the Riemann-Hilbert
representation for the Painlev\'e II equation (see Definition \ref{def3}
below), such that
\begin{equation}\label{qz1.14}
    \lim_{n\to\infty} \Prob\biggl(
    \frac{L^\symmO_{n,[\alpha\sqrt{2n}]}-2\sqrt{2n+m}}{(2n+m)^{1/6}}
\le x \biggr) = F^\symmO(x;w).
\end{equation}
The new class of distributions $F^\symmO(x;w)$ interpolates $F_4$ and $F_1$ as
$w\to\infty$ and $w=0$ respectively, and satisfies $\lim_{w\to-\infty}
F^\symmO(x;w)=0$, so \eqref{qz1.14} is consistent with
\eqref{qz1.10}--\eqref{qz1.12}. Alternatively, since $F_4=F_1^{(2)}$,
$F^\symmO(x;w)$ interpolates the limiting distributions of the second and the
first eigenvalues of a random GOE matrix.

The meaning of $F^\symmO(x;w)$ in terms of random matrices is not clear, but
there is a Coulomb gas interpretation. In (7.64)--(7.65) of \cite{PartI},
the following density function is introduced. Suppose $2N$ ordered particles on
the positive real line, $0< \xi_{2N}<\cdots<\xi_2<\xi_1$, are distributed
according to the density function
\begin{equation}\label{qz100}
  \frac{1}{Z_{2N}} e^{A\sum_{j=1}^{2N}(-1)^j\xi_j}
  \prod_{1\le i<j\le 2N} (\xi_i-\xi_j)
  \prod_{j=1}^{2N}e^{-\xi_j/2}d\xi_j,
\end{equation}
where $Z_{2N}$ is the normalization constant. Hence in addition to the usual
Coulomb gas interaction, there is a additional attraction between 
neighboring pairs $\xi_{2j-1}$ and $\xi_{2j}$, $j=1,\cdots,N$. When $A=0$, this
attraction vanishes, and one sees that \eqref{qz100} is the eigenvalue density
for the $2N\times 2N$ Laguerre orthogonal ensemble (LOE). On the other hand,
when $A\to\infty$, the neighboring particles $\xi_{2j-1}$, $\xi_{2j}$,
$j=1,2,\cdots,N$, coalesce, and \eqref{qz100} becomes the eigenvalue density
function for the $N\times N$ Laguerre symplectic ensemble (LSE). Thus
\eqref{qz100} interpolates LOE and LSE eigenvalue distributions. This density
function arises in a symmetric version of the growth model with a growth rule
given by the exponential distribution considered in Proposition 1.4 of
\cite{kurtj:shape}. A discrete version of the above density function was
considered in (4.27) of \cite{BR3}, and the limiting distribution of the
largest particle was precisely $F^\symmO(x;w)$. Now \emph{formally} taking the
exponential limit (set $q=1-1/L$ and $\alpha=1-A/L$ and take $L\to\infty$) of
(4.27) of \cite{BR3} (which is convincing, but is not justified yet), we see
that if we set
\begin{equation}\label{qz101}
  A= \frac{w}{N^{1/3}}
\end{equation}
in \eqref{qz100} and take $N\to\infty$, the scaled largest particle
$(\xi_1-4N)/(2N^{1/3})$ has the limiting distribution $F^\symmO(x;w)$. We plan
to exploit the justification of the exponential limit in a later publication.

\medskip
On the other hand, the longest decreasing subsequence corresponds to the
longest down/right path from $(0,1)$ to $(1,0)$ in the above point selection
process. Thus it is clear that the distribution of $L^\symmS_{n,m}$ is
insensitive to the number of fixed points. See Section \ref{sec:results} for
results and discussions.

The other ensemble we consider is the set of signed involutions. A signed
permutation $\pi$ is a bijection from $\{-n,\cdots, -1,1,\cdots, n\}$ onto
itself which satisfies $\pi(x)=-\pi(-x)$. The limiting distribution for a
random signed permutation is obtained by \cite{TW3, Bo}. In this paper we
consider random signed involutions with/without constraint on fixed points and
also negated points ($\pi(x)=-x$). See Section \ref{sec:results} for results.
Especially, we obtain another one-parameter family of distributions which now
interpolates $F_2$ and $F_1^2$. Here $F_1^2$ has the meaning as the limiting
distribution for the largest `eigenvalue' of the super-imposition of the
eigenvalues of two random GOE matrices. We note that in this case, $F_2$ is
equal to the second largest `eigenvalue' of such super-imposition (see the
discussion at the end of Section \ref{sec:results}).

\bigskip

For convenience of future reference, we summarize various definitions
introduced above. By the \emph{$k^{th}$ row/column of $\pi$} we mean the
$k^{th}$ row/column of the corresponding Young diagram under the
Robinson-Schensted map.

\begin{defn}\label{def-int}
Let $S_n$ be the symmetric group of $n$ letters and let
$S_n^\symmu$ be the set of the bijections from
$\{-n,\cdots,-2,-1,1,2,\cdots,n\}$ onto itself satisfying $\pi(x)=-\pi(-x)$.
We define
\begin{eqnarray}
\tilde{S}_n&=&\{\pi\in S_n: \pi=\pi^{-1}\}, \\
S_{n,m} &=& \{ \pi\in \tilde{S}_{2n+m} :\ |\{ x: \pi(x)=x\}|=m \},\\
\tilde{S}^\symmu_n &=&\{ \pi\in S_n^\symmu : \pi=\pi^{-1} \}, \\
S^{\symmu}_{n,m_+,m_-} &=& \{ \pi\in \tilde{S}^\symmu_{2n+m_++m_-} :\
|\{ x: \pi(x)=x\}|=2m_+ ,\ |\{ x: \pi(x)=-x \}|=2m_- \}, \\
\tilde{L}_n^{(k)}(\pi) &=& \text{the length of the $k^{th}$ row of $\pi
\in\tilde{S}_n$},
\qquad \tilde{L}_{n}=\tilde{L}^{(1)}_{n}, \\
\tilde{L}_n^{\symmu,(k)}(\pi) &=& \text{the length of the $k^{th}$ row
of $\pi\in\tilde{S}^\symmu_n$},
\qquad \tilde{L}^\symmu_{n}=\tilde{L}^{\symmu,(1)}_{n}, \\
L^{\symmO,(k)}_{n,m}(\pi) &=& \text{the length of the $k^{th}$ row
of $\pi\in S^\symmO_{n,m}$},
\qquad L^{\symmO}_{n,m}=L^{\symmO,(1)}_{n,m},\\
L^{\symmS,(k)}_{n,m}(\pi) &=& \text{the length of the $k^{th}$ column
of $\pi\in S^\symmS_{n,m}$},
\qquad L^{\symmS}_{n,m}=L^{\symmS,(1)}_{n,m},\\
L^{\symmu,(k)}_{n,m_+,m_-}(\pi) &=& \text{the length of the $k^{th}$ row
of $\pi\in S^\symmu_{n,m_+,m_-}$},
\qquad L^{\symmu}_{n,m_+,m_-}=L^{\symmu,(1)}_{n,m_+,m_-}.
\end{eqnarray}
\end{defn}

The results of this paper were announced in \cite{BR3}. Since we completed
this paper, there have been two applications. One is to random vicious walker
models \cite{forrester99, BR3}, and the other is to polynuclear growth models
\cite{SpohnP1, SpohnP2, BR4}. Indeed there are bijections between the above two
applications and various ensembles considered in this paper, and thus the
results in this paper can be employed to answer asymptotic questions in the
above applications.

\bigskip

The proofs of our theorems use the Poissonization and de-Poissonization
scheme of \cite{Jo2, BDJ}. We define the Poisson generating function, for
example, for $L^\symmO_{n,m}$ by (see Definition \ref{def-smml})
\begin{equation}
  Q^\symmO_l(\lambda_1,\lambda_2) :=
e^{-\lambda_1-\lambda_2} \sum_{n_1,n_2\ge 0}
\frac{\lambda_1^{n_1} \lambda_2^{n_2}}{n_1!n_2!}
\Prob\bigl( L^\symmO_{n_2,n_1}\le l\bigr).
\end{equation}
A generalization of the de-Poissonization lemma due to Johansson \cite{Jo2}
yields that $\Prob\bigl( L^\symmO_{n_2,n_1}\le l\bigr) \sim
Q^\symmO_l(n_1,n_2)$ as $n_1,n_2\to\infty$ ; see Section \ref{as-dePoi} for
the precise statement. Thus if we obtain the asymptotics of the generating
function, the asymptotics of the coefficients will follow. The point of the
scheme is that the Poisson generating functions can be expressed in terms
of Toeplitz and/or Hankel determinants. The necessary algebraic work for
this purpose was done in our earlier paper \cite{PartI}.  The general
theory of orthogonal polynomials then tells us that Toeplitz/Hankel
determinants can be expressed in terms of orthogonal polynomials. In the
paper \cite{PartI}, it turned out that for all the ensembles being
discussed in this paper, we need only one family of orthogonal polynomials
: $\pi_n(z;t)=z^n+\cdots$ which is orthogonal with respect to the weight
$e^{t\cos\theta} d\theta/(2\pi)$ on the unit circle. This orthogonal
polynomial is precisely the same orthogonal polynomial used in \cite{BDJ}
to analyze the random permutation problem. The authors in \cite{BDJ}
computed the uniform asymptotics of the normalization constant $N_n(t)$ of
$\pi_n$ as $n,t\to\infty$ using the steepest-descent analysis for the
corresponding Riemann-Hilbert problem (see \eqref{as4} below). The
difference between the present paper and
\cite{BDJ} is that we need $\pi_n(-\alpha;t)$ for all $\alpha\ge 0$, in
contrast to \cite{BDJ} where only one quantity, $N_n(t)$, is needed. But in
order to analyze $N_n(t)$, \cite{BDJ} controlled in a uniform way the
asymptotic behavior of the solution to the associated Riemann-Hilbert
problem. Therefore the asymptotics of $\pi_n(-\alpha;t)$ for $\alpha$
uniformly apart from $1$ can be (almost) directly read off from the
analysis of \cite{BDJ}, which eventually imply \eqref{qz1.10}. The point
$z=-1$ ($\alpha=1$) in the complex plane plays a special role in this
Riemann-Hilbert analysis as discussed in \cite{BDJ}---it is the point where
a gap starts to open up in the support of the associated equilibrium
measure as the relation of $t$ to $n$ varies. When $\alpha\to 1$ as
$n\to\infty$ according to the scaling \eqref{qz1.13} (which is required for
\eqref{qz1.14}) , we need more careful analysis of the Riemann-Hilbert problem
which is the new part of the asymptotic analysis of the orthogonal
polynomials and the Riemann-Hilbert problem. In this paper we establish
this goal by extending the analysis of \cite{BDJ}. In the analysis of the
Riemann-Hilbert problem in Section \ref{RHP}, we give a rather sketchy
presentation for the parts which overlap the work of \cite{BDJ}, but we
give a full proof for the new analysis required for the case $\alpha\to 1$.
We also rework portions of \cite{BDJ} as necessary for consistency of
presentation.

\bigskip

This paper is organized as follows. Section \ref{sec:dist} defines the
Tracy-Widom distribution functions as well as new classes of distribution
functions which will be used to state the main results, and discusses their
properties (Lemma \ref{lem12}). The main results of the paper are then
stated in Section \ref{sec:results}. Determinantal formulae and orthogonal
polynomial expressions for the Poisson generating functions are taken from
\cite{PartI} and summarized in Section \ref{sec:sum}. In Section
\ref{as-op}, we state the main estimates of the relevant quantities of
orthogonal polynomials; these estimates are key to the proofs of the
theorems of Section \ref{sec:results}. The de-Poissonization lemmas are
stated in Section \ref{as-dePoi}. Proofs of the main theorems are given in
Section
\ref{sec:pf} for involutions with a constraint on the number of fixed
points (Theorems \ref{newthm1},
\ref{newthm2}, \ref{newthm3} and \ref{thm-second1}), and in Section
\ref{as-asymp3} for general involutions and equivalently for $M^1_n$ (Theorems
\ref{thm22} and \ref{thm-second2}), respectively. The case when $\alpha>1$ is
considered in Section \ref{sec:bigalpha} (see Remark to Theorem \ref{newthm3}).
Finally the Riemann-Hilbert analysis is given in Section \ref{RHP}, which
proves the propositions in Section \ref{as-op}.

\bigskip
\noindent {\bf Notational remarks.} The ensemble $\tilde{S}_n$ in the present
paper is identical to $\tilde{S}^\symmO_n$ in \cite{PartI}. In \cite{PartI},
$\tilde{S}^\symmS_n$ is introduced to denote the ensemble of neginvolutions,
and the authors investigated the longest increasing subsequence of
$\pi\in\tilde{S}^\symmS_n$. But there is a bijection between
$\tilde{S}^\symmS_n$ and $\tilde{S}^\symmO_n$ and the longest increasing
subsequence of $\pi\in\tilde{S}^\symmS_n$ corresponds to the longest decreasing
subsequence of the image of $\pi$. In the present paper, we choose the view
point of considering both the increasing and decreasing subsequences of
involutions of the same ensemble rather than considering only the increasing
subsequences of involutions of the different ensembles.

\medskip
\noindent {\bf Acknowledgments.} The authors would like to thank Percy Deift
for helpful discussions and encouragement, especially for his help in
proving Lemma 2.1. We would also like to acknowledge many useful
conversations and communications with Peter Forrester, Kurt Johansson,
Charles Newman and Harold Widom. Special thanks is due to the referee who
gave us crucial advice, improving the exposition of the paper
significantly. The work of the first author was supported in part by a
Sloan Doctoral Dissertation Fellowship during the academic year 1998-1999
when he was a graduate student at Courant Institute of Mathematical
Sciences.

%
%

\section{Limiting distribution functions}\label{sec:dist}

Let $u(x)$ be the solution of the Painlev\'e II (PII) equation,
\begin{equation}\label{as10}
  u_{xx}=2u^3+xu,
\end{equation}
with the boundary condition
\begin{equation}\label{as11}
  u(x)\sim -Ai(x)\quad\text{as}\quad x\to +\infty,
\end{equation}
where $Ai$ is the Airy function.
The proof of the (global) existence and the uniqueness of the solution
was first established in \cite{HM}:
the asymptotics as $x\to \pm\infty$ are
(see e.g. \cite{HM,DZ2})
\begin{eqnarray}
\label{as17}  u(x) &=&
-Ai(x) + O\biggl( \frac{e^{-(4/3)x^{3/2}}}{x^{1/4}}\biggr),
\qquad\text{as $x\to +\infty$,}\\
\label{as18}  u(x) &=&
-\sqrt{\frac{-x}{2}}\biggl( 1+O\bigl(\frac1{x^2}\bigr)\biggr),
\qquad\text{as $x\to -\infty$.}
\end{eqnarray}
Recall that $Ai(x)\sim \frac{e^{-(2/3)x^{3/2}}}{2\sqrt{\pi}x^{1/4}}$
as $x\to +\infty$.
Define
\begin{equation}\label{as12}
  v(x):= \int_{\infty}^{x} (u(s))^2 ds,
\end{equation}
so that $v'(x)= (u(x))^2$.

We can now introduce the Tracy-Widom (TW) distributions.
(Note that $q:=-u$, which Tracy and Widom used in their papers,
solves the same
differential equation with the boundary condition
$q(x)\sim +Ai(x)$ as $x\to\infty$.)

\begin{defn}[TW distribution functions]\label{def1}
Set
\begin{eqnarray}
\label{as46-1}   F(x) &:=& \exp\biggl(\frac12\int_x^{\infty} v(s)ds\biggr)
= \exp\biggl(-\frac12\int_x^{\infty} (s-x)(u(s))^2ds\biggr),\\
\label{as47-1}   E(x) &:=& \exp\biggl(\frac12\int_x^{\infty} u(s)ds\biggr),
\end{eqnarray}
and set
\begin{eqnarray}
   F_2(x) :=& F(x)^2 &=\exp\biggl(-\int_x^{\infty} (s-x)(u(s))^2ds\biggr),\\
   F_1(x) :=& F(x)E(x) &= \bigl(F_2(x)\bigr)^{1/2}
e^{\frac12\int_x^{\infty} u(s)ds} ,\\
   F_4(x) :=& F(x)\bigl[E(x)^{-1}+E(x)\bigr]/2
&= \bigl(F_2(x)\bigr)^{1/2}
\biggl[ e^{-\frac12\int_x^{\infty} u(s)ds} +
e^{\frac12\int_x^{\infty} u(s)ds} \biggr] /2 .
\end{eqnarray}
\end{defn}

In \cite{TW1} and \cite{TW2},
Tracy and Widom proved that under proper centering
and scaling, the distribution of the largest eigenvalue of
a random GUE/GOE/GSE matrix converges to $F_2(x)$ / $F_1(x)$ / $F_4(x)$
as the size of the matrix becomes large.
We note that from the asymptotics \eqref{as17} and \eqref{as18},
for some positive constant $c$,
\begin{eqnarray}
\label{as46}   F(x) &=& 1+ O\bigl( e^{-cx^{3/2}}\bigr)
\qquad \text{as $x\to+\infty$,}\\
\label{as47}   E(x) &=& 1+ O\bigl( e^{-cx^{3/2}}\bigr)
\qquad \text{as $x\to+\infty$,}\\
\label{as48}   F(x) &=& O\bigl( e^{-c|x|^{3}}\bigr)
\qquad\qquad \text{as $x\to-\infty$,}\\
\label{as49}   E(x) &=& O\bigl( e^{-c|x|^{3/2}}\bigr)
\quad\qquad \text{as $x\to-\infty$.}
\end{eqnarray}
Hence in particular, $\lim_{x\to+\infty} F_\beta(x) = 1$ and
$\lim_{x\to-\infty} F_\beta(x) = 0$, $\beta=1,2,4$.
Monotonicity of $F_\beta(x)$ follows from the fact that $F_\beta(x)$
is the limit of a sequence of distribution functions.
Therefore $F_\beta(x)$ is indeed a distribution function.

\begin{defn}\label{def2}
  Define $\chi_{\text{GOE}}$, $\chi_{\text{GUE}}$ and
$\chi_{\text{GSE}}$ to be random variables whose distribution
functions are given by $F_1(x)$, $F_2(x)$ and $F_4(x)$, respectively.
Define $\chi_{\text{GOE}^2}$
to be a random variable with the distribution function $F_1(x)^2$.
\end{defn}

\bigskip
As indicated in Introduction, we need new classes of distribution functions
to describe the transitions from $\chi_{\text{GSE}}$ to $\chi_{\text{GOE}}$
and from $\chi_{\text{GUE}}$ to $\chi_{\text{GOE}^2}$. First we consider
the Riemann-Hilbert problem (RHP) for the Painlev\'e II equation \cite{FN,
JMU}.  Let $\Gamma$ be the real line $\R$, oriented from $+\infty$ to
$-\infty$, and let $m(\cdot\thinspace;x)$ be the solution of the following
RHP:
\begin{equation}\label{as20}
 \begin{cases}
    m(z;x) \qquad \text{is analytic in $z\in\C\setminus\Gamma$,}\\
    m_+(z;x)=m_-(z;x) \begin{pmatrix} 1& -e^{-2i(\frac43z^3+xz)}\\
e^{2i(\frac43z^3+xz)}& 0 \end{pmatrix} \quad \text{for $z\in\Gamma$,}\\
    m(z;x) = I+O\bigl(\frac1{z}\bigr) \qquad \text{as $z\to\infty$.}
 \end{cases}
\end{equation}
Here $m_+(z;x)$ (resp,. $m_-$) is the limit of $m(z';x)$ as $z'\to z$
from the left (resp., right) of the contour $\Gamma$:
$m_\pm(z;x)=\lim_{\epsilon\downarrow 0}m(z\mp i\epsilon;x)$.
Relation \eqref{as20} corresponds to
the RHP for the PII equation with the special monodromy data
$p=-q=1, r=0$ (see \cite{FN, JMU}, also \cite{FZ,DZ2}).
In particular if the solution is expanded at $z=\infty$,
\begin{equation}\label{as7.27}
   m(z;x) = I+ \frac{m_1(x)}{z} + O\bigl(\frac1{z^2}\bigr),
\qquad \text{as $z\to\infty$},
\end{equation}
we have
\begin{eqnarray}
\label{as22}   2i(m_1(x))_{12} = -2i(m_1(x))_{21} &=& u(x), \\
\label{as23}   2i(m_1(x))_{22} = -2i(m_1(x))_{11}&=& v(x),
\end{eqnarray}
where $u(x)$ and $v(x)$ are defined in \eqref{as10}-\eqref{as12}.

Now we define two one-parameter families of distribution functions.

\begin{defn}\label{def3}
  Let $m(z;x)$ be the solution of RHP \eqref{as20} and
denote by $m_{jk}(z;x)$ the $(jk)$-entry of $m(z;x)$.
For $w>0$, define
\begin{equation}\label{as61}
   F^\symmO(x;w) :=
F(x) \biggl\{\bigl[m_{22}(-iw;x)-m_{12}(-iw;x) \bigr]E(x)^{-1}
+ \bigl[m_{22}(-iw;x)+m_{12}(-iw;x)\bigr]E(x) \biggr\}/2,
\end{equation}
and for $w<0$, define
\begin{equation}\label{as62}
 \begin{split}
   F^\symmO(x;w) &:=
e^{\frac83w^3-2xw}F(x) \\
&\times \biggl\{\bigl[-m_{21}(-iw;x)+m_{11}(-iw;x) \bigr]E(x)^{-1}
- \bigl[m_{21}(-iw;x)+m_{11}(-iw;x)\bigr]E(x) \biggr\}/2.
 \end{split}
\end{equation}
Also define
\begin{eqnarray}
   F^\symmu(x;w) &:=& m_{22}(-iw;x)F_2(x),
\quad\qquad\qquad\qquad w>0,\\
   F^\symmu(x;w) &:=& -e^{\frac83w^3-2xw}m_{21}(-iw;x)F_2(x),
\qquad w<0.
\end{eqnarray}
\end{defn}

\medskip
First, $F^\symmO(x;w)$ and $F^\symmu(x;w)$ are real from Lemma \ref{lem12} (i)
below. Note that $F^\symmO(x;w)$ and $F^\symmu(x;w)$ are continuous at $w=0$
since at $z=0$, the jump condition of the RHP \eqref{as20} implies
$(m_{12})_+(0;x)= -(m_{11})_-(0;x)$ and $(m_{22})_+(0;x)= -(m_{21})_-(0;x)$. In
fact, $F^\symmO(x;w)$ and $F^\symmu(x;w)$ are entire in $w\in\C$ from the RHP
\eqref{as20}.

From \eqref{as46}-\eqref{as49} and \eqref{as94}-\eqref{as97} below,
we see that
\begin{equation}
  \lim_{x\to+\infty} F^\symmO(x;w), F^\symmu(x;w) =1, \qquad
\lim_{x\to-\infty} F^\symmO(x;w), F^\symmu(x;w) = 0
\end{equation}
for any fixed $w\in\R$.
Also Theorem \ref{newthm2} below shows that $F^\symmO(x;w)$ and $F^\symmu(x;w)$
are limits of distribution functions, implying that
they are monotone in $x$.
Therefore, $F^\symmO(x;w)$ and $F^\symmu(x;w)$ are indeed
distribution functions for each $w\in\R$.

\begin{defn}
  Define $\chi^\symmO_w$ and $\chi^\symmu_w$ to be random variables
with distribution functions $F^\symmO(x;w)$ and
$F^\symmu(x;w)$, respectively.
\end{defn}

We close this section by summarizing some properties of $m(-iw;x)$ in the
following lemma. In particular the lemma implies that $F^\symmO(x;w)$
interpolates between $F_4(x)$ and $F_1(x)$, and $F^\symmu(x;w)$
interpolates between $F_2(x)$ and $F_1(x)^2$ (see Corollary \ref{cor16}).

\begin{lem}\label{lem12}
  Let $\sigma_3=\bigl(\begin{smallmatrix} 1&0\\0&-1 \end{smallmatrix}
\bigr)$, $\sigma_1=\bigl(\begin{smallmatrix} 0&1\\1&0 \end{smallmatrix}
\bigr)$, and set $[a,b]=ab-ba$.
\begin{enumerate}
\item
  For real $w$, $m(-iw;x)$ is real.
\item
For fixed $w\in\R$, we have
\begin{eqnarray}
\label{as94}  m(-iw;x) &=& \bigl( I+\bigl( e^{-cx^{3/2}} \bigr) \bigr)
\begin{pmatrix}
1 & -e^{\frac83w^3-2xw} \\ 0 & 1
\end{pmatrix}, \qquad  \text{$w>0$, $x\to+\infty$,}\\
\label{as95}  m(-iw;x) &=& \bigl( I+\bigl( e^{-cx^{3/2}} \bigr) \bigr)
\begin{pmatrix}
1 & 0 \\ -e^{-\frac83w^3+2xw} & 1
\end{pmatrix}, \qquad  \text{$w<0$, $x\to+\infty$,}\\
\label{as96}  m(-iw;x) &\sim& \frac1{\sqrt{2}} \biggl( \begin{smallmatrix}
1&-1\\ 1&1 \end{smallmatrix} \biggr)
e^{(-\frac43w^3+xw)\sigma_3}
e^{(\frac{\sqrt{2}}3(-x)^{3/2}+\sqrt{2}w^2(-x)^{1/2})\sigma_3},
\quad  \text{$w>0$, $x\to-\infty$,}\\
\label{as97}  m(-iw;x) &\sim& \frac1{\sqrt{2}} \biggl( \begin{smallmatrix}
1&1\\ -1&1 \end{smallmatrix} \biggr)
e^{(-\frac43w^3+xw)\sigma_3}
e^{(-\frac{\sqrt{2}}3(-x)^{3/2}-\sqrt{2}w^2(-x)^{1/2})\sigma_3},
\quad  \text{$w<0$, $x\to-\infty$.}
\end{eqnarray}
\item
  For any $x$, we have
\begin{equation}\label{as98}
  \lim_{w\to 0^+} m(-iw;x) =
\lim_{w\to 0^-} \sigma_1 m(-iw;x)\sigma_1  = \begin{pmatrix}
\frac12\bigl(E(x)^2+E(x)^{-2}\bigr) & -E(x)^2 \\
\frac12\bigl(-E(x)^2+E(x)^{-2}\bigr) & E(x)^2
\end{pmatrix}.
\end{equation}
\item
   For fixed $w\in\R\setminus\{0\}$,
$m(-iw;x)$ solves the differential equation
\begin{equation}\label{as62.5}
   \frac{d}{dx}m = w[m,\sigma_3]+ u(x)\sigma_1 m,
\end{equation}
where $u(x)$ is the solution of the PII equation \eqref{as10}, \eqref{as11}.
\item
   For fixed $x$, $m(-iw;x)$ solves
\begin{equation}
   \frac{\partial}{\partial w}m = (-4w^2+x)[m,\sigma_3] -4w u(x)\sigma_1 m
-2\biggl(\begin{smallmatrix}u^2&-u'\\u'&-u^2\end{smallmatrix}\biggr) m.
\end{equation}
\item
For any $x$, we have
\begin{equation}
  m(z;x)=\sigma_1m(-z;x)\sigma_1.
\end{equation}
\end{enumerate}
\end{lem}

\begin{cor}\label{cor16}
   We have
\begin{eqnarray}
   F^\symmO(x;0) &=& F_1(x),\\
   \lim_{w\to\infty} F^\symmO(x;w) &=& F_4(x), \\
   \lim_{w\to-\infty} F^\symmO(x;w) &=& 0, \\
\label{as116}   F^\symmu(x;0) &=& F_1(x)^2, \\
   \lim_{w\to\infty} F^\symmu(x;w) &=& F_2(x), \\
   \lim_{w\to-\infty} F^\symmu(x;w) &=& 0.
\end{eqnarray}
\end{cor}

\begin{proof}
   The values at $w=0$ follow from \eqref{as98}.
For $w\to\pm\infty$, note from the RHP \eqref{as20} that we have
$\lim_{z\to\infty} m(z;x)= I$.
\end{proof}

\begin{proof}[Proof of Lemma \ref{lem12}]
  Let $v(z)=v(z;x)$ denote the jump matrix of the RHP \eqref{as20}.
  Since $\overline{v(-z)}=v(z)$ for $z\in\R$,
$M(z):=\overline{m(-\overline{z};x)}$ also solves the same RHP.
By the uniqueness of the solution of the RHP \eqref{as20}, we have
\begin{equation}
  \overline{m(-\overline{z};x)} = m(z;x), \qquad z\in\C\setminus\R.
\end{equation}
Thus,  $m(-iw;x)$ is real for $w\in\R$, thus proving (i).

\smallskip
 By the symmetry of the jump matrix,
$\sigma_1 v(-z)^{-1} \sigma_1 = v(z)$, we obtain, by a argument similar to (i),
\begin{equation}\label{as69}
  \sigma_1 m(-z;x) \sigma_1 = m(z;x),
\end{equation}
which is (vi).

\smallskip
  The asymptotics results (ii) as $x\to\pm\infty$ follow from the
calculations in Section 6, pp.329-333, of \cite{DZ2}.

\smallskip
For the proof of (iv),
define a new matrix function
\begin{equation}\label{as70}
   f(z;x):=m(z;x) e^{-i\theta\sigma_3}, \qquad \theta:= \frac43z^3+xz.
\end{equation}
Then $f(\cdot\thinspace;x)$ satisfies the jump condition $f_+(z;x)=f_-(z;x)
\bigl(\begin{smallmatrix} 1&-1\\1&0 \end{smallmatrix} \bigr)$
for $z\in\R$, and $f(z;x)e^{i\theta\sigma_3}
\to I$ as $z\to\infty$.
Since the jump matrix for $f(z;x)$ is independent of $x$,
$f'(z;x)$, the derivative with respect to $x$, satisfies
$f'_+(z;x)=f'_-(z;x)
\bigl(\begin{smallmatrix} 1&-1\\1&0 \end{smallmatrix} \bigr)$, and
$f'e^{i\theta\sigma_3}+i\theta' f\sigma_3e^{i\theta\sigma_3}
\to 0$ as $z\to\infty$.
Hence $f'f^{-1}$ has no jump across $\R$, and
satisfies $f'f^{-1} +i\theta'f\sigma_3 f^{-1} \to 0$ as $z\to\infty$.
If we write $m(z;x)= I + \frac{m_1(x)}{z} + O(1/z^2)$ as $z\to\infty$,
we have $i\theta'f\sigma_3 f^{-1}= iz\sigma_3 +i [m_1,\sigma_3]+O(z^{-1})$
as $z\to\infty$.
Thus $f'f^{-1}$ is entire and as $z\to\infty$,
$f'f^{-1}\sim -iz\sigma_3 -i [m_1,\sigma_3]$.
Therefore by Liouville's theorem, we obtain
\begin{equation}\label{as102}
   f'(z;x)(f(z;x))^{-1} = -iz\sigma_3 - i[m_1,\sigma_3].
\end{equation}
Recalling that
$u(x)= 2i(m_1(x))_{12} = -2i(m_1(x))_{21}$
in \eqref{as22}, we have $[m_1,\sigma_3]=iu(x)\sigma_1$.
Changing $f$ to $m$ from \eqref{as70}, \eqref{as102} is
\begin{equation}\label{e-q}
   \frac{d}{dx} m(z;x)= iz[m(z;x),\sigma_3] + u(x)\sigma_1 m(z;x).
\end{equation}
This is \eqref{as62.5} when $z=-iw$.

\smallskip
The proof of (v) is very similar to that of (iv), and the detail is
left to the reader.
We only note that in the derivation of (v), we need the identity
\begin{equation}
  \frac{d}{dx} m_1= i[m_2,\sigma_3]-i[m_1,\sigma_3]m_1,
\end{equation}
which can be obtained from \eqref{e-q} by setting $m(z;x)=I+\frac{m_1(x)}{z}
+\frac{m_2(x)}{z^2} +O(1/z^3)$ as $z\to\infty$.

\smallskip
Finally we prove (iii).
  Note that $\lim_{w\to 0^\pm} m(-iw;x)=m_\pm(0;x)$.
 From the jump condition at $z=0$, we have
\begin{equation}\label{as63}
  m_+(0;x)=m_-(0;x)\biggl(\begin{smallmatrix} 1&-1\\1&0
\end{smallmatrix} \biggr).
\end{equation}
Letting $z\to 0$, $Im z>0$, in (vi), we have $\sigma_1
 m_+(0;x) \sigma_1 = m_-(0;x)$, which together with \eqref{as63}
implies that $m_+(0;x)=\sigma_1 m_+(0;x)\sigma_1
\bigl(\begin{smallmatrix} 1&-1\\1&0
\end{smallmatrix} \bigr)$.
Thus we have
\begin{equation}\label{as65}
   m_+(0;x)= \biggl(\begin{smallmatrix} a(x)&b(x)\\a(x)+b(x)&-b(x)
\end{smallmatrix} \biggr)
\end{equation}
for some $a(x)$, $b(x)$.
Also the condition $\det v(z)=1$ for all $z\in\R$ implies that
$\det m(z;x)=1$ for all $z\in\C\setminus\R$, and hence we have
\begin{equation}\label{as66}
   b^2+2ab+1=0.
\end{equation}
Now letting $z\to 0$, $Im z<0$ in \eqref{e-q}, we obtain
\begin{equation}
   m'_+(0;x)(m_+(0;x))^{-1} = \biggl( \begin{smallmatrix}
0&u(x)\\u(x)&0 \end{smallmatrix} \biggr).
\end{equation}
Thus from \eqref{as65} and \eqref{as66}, $b'/b=-u$,
which has the solution
\begin{equation}
   b(x)=b(y) e^{-\int_y^x u(s)ds}.
\end{equation}
 From \eqref{as94} with $w=0^+$, we have
\begin{equation}
   b(x)= (m_{12})_+(0;x) \to -1, \qquad \text{as $x\to+\infty$}.
\end{equation}
Therefore, $b(x)= -e^{\int_x^\infty u(s)ds}$,
which is $-E(x)^2$ from \eqref{as47-1}.
Now \eqref{as66} gives $a(x)=\frac12(E(x)^2+E(x)^{-2})$,
proving \eqref{as98}.
\end{proof}

%
%

\section{Statement of results}\label{sec:results}

\subsection{Involutions with constraint on the number of
fixed points}

Recall (Definition \ref{def-int} in Introduction) the ensembles
$\tilde{S}_{n,m}$, $\tilde{S}^\symmu_{n,m_+,m_-}$ of (signed) involutions with
constraint on the number of fixed (and negated) points. We scale the random
variables:
\begin{eqnarray}
\label{res1}  \chi^\symmO_{n,m} &:=&
\frac{L^{\symmO}_{n,m}-2\sqrt{2n+m}}{(2n+m)^{1/6}}, \\
\label{res2}  \chi^\symmS_{n,m} &:=&
\frac{L^{\symmS}_{n,m}-2\sqrt{2n+m}}{(2n+m)^{1/6}}, \\
\label{res3}  \chi^\symmu_{n,m_+,m_-} &:=&
\frac{L^{\symmu}_{n,m_+,m_-}-
2\sqrt{4n+2m_++2m_-}}{2^{2/3}(4n+2m_++2m_-)^{1/6}}.
\end{eqnarray}

\begin{thm}\label{newthm1}
 For fixed $\alpha$ and $\beta$,
we have:
\begin{itemize}
\item
\begin{align}
\label{qz3.4}
  \lim_{n\to\infty} \Prob\bigl( \chi^\symmO_{n,[\sqrt{2n}\alpha]}
\le x \bigr) &= F_4(x),
\qquad 0\le \alpha<1, \\
  \lim_{n\to\infty} \Prob\bigl( \chi^\symmO_{n,[\sqrt{2n}]}
\le x \bigr) &= F_1(x), \\
  \lim_{n\to\infty} \Prob\bigl( \chi^\symmO_{n,[\sqrt{2n}\alpha]}
\le x \bigr) &= 0, \qquad\qquad \alpha>1.
\end{align}
\item
\begin{equation}\label{qz3.7}
  \lim_{n\to\infty} \Prob\bigl( \chi^\symmS_{n,[\sqrt{2n}\beta]}
\le x \bigr) = F_1(x), \qquad \beta\ge 0.
\end{equation}
\item
\begin{align}
  \lim_{n\to\infty}
\Prob\bigl( \chi^\symmu_{n,[\sqrt{n}\alpha],[\sqrt{n}\beta]}
\le x \bigr) &= F_2(x),
\qquad 0\le\alpha<1, \beta\ge 0, \\
  \lim_{n\to\infty}
\Prob\bigl( \chi^\symmu_{n,[\sqrt{n}],[\sqrt{n}\beta]}
\le x \bigr) &= F_1(x)^2,
\qquad \beta\ge 0, \\
  \lim_{n\to\infty}
\Prob\bigl( \chi^\symmu_{n,[\sqrt{n}\alpha],[\sqrt{n}\beta]}
\le x \bigr) &= 0,
\qquad\qquad \alpha>1, \beta\ge 0.
\end{align}
\end{itemize}
\end{thm}

As indicated earlier in the introduction, as $\alpha\to 1$ at a certain rate,
we see smooth transitions.

\begin{thm}\label{newthm2}
For fixed $w\in\R$ and $\beta\ge 0$,
we have
\begin{align}
  \lim_{n\to\infty} \Prob\bigl( \chi^\symmO_{n,m}
\le x \bigr) &= F^\symmO(x;w),
\qquad m=[\sqrt{2n} - 2w(2n)^{1/3}],\\
     \lim_{n\to\infty}
\Prob\bigl( \chi^\symmu_{n,m_+,m_-}
\le x \bigr) &= F^\symmu(x;w),
\quad\qquad m_+=[\sqrt{n} - 2wn^{1/3}],\ m_-=[\sqrt{n}\beta].
\end{align}
\end{thm}

From Corollary \ref{cor16},
this result is consistent with Theorem \ref{newthm1}.
We also have convergence of moments.

\begin{thm}\label{newthm3}
For any $p=1,2,3,\cdots$, the followings hold.
For fixed $\alpha$ and $\beta$,
\begin{align}
  \lim_{n\to\infty} \Exp\bigl( (\chi^\symmO_{n,[\sqrt{2n}\alpha]})^p \bigr)
&= \Exp\bigl( (\chi_{\text{GSE}})^p \bigr),
\qquad 0\le \alpha<1, \\
  \lim_{n\to\infty} \Exp\bigl( (\chi^\symmO_{n,[\sqrt{2n}]})^p \bigr)
&= \Exp\bigl( (\chi_{\text{GOE}})^p \bigr),\\
  \lim_{n\to\infty} \Exp\bigl( (\chi^\symmS_{n,[\sqrt{2n}\beta]})^p \bigr)
&= \Exp\bigl( (\chi_{\text{GOE}})^p \bigr),
\qquad 0\le \beta, \\
  \lim_{n\to\infty} \Exp\bigl(
(\chi^\symmu_{n,[\sqrt{n}\alpha],[\sqrt{n}\beta]})^p \bigr)
&= \Exp\bigl( (\chi_{\text{GUE}})^p \bigr),
\qquad 0\le\alpha<1, \beta\ge 0,\\
  \lim_{n\to\infty} \Exp\bigl(
(\chi^\symmu_{n,[\sqrt{n}],[\sqrt{n}\beta]})^p \bigr)
&= \Exp\bigl( (\chi_{\text{GOE}^2})^p \bigr),
\qquad \beta\ge 0.
\end{align}
Also for fixed $w\in\R$ and $\beta\ge0$,
\begin{align}
  \lim_{n\to\infty} \Exp\bigl( (\chi^\symmO_{n,m})^p \bigr)
&= \Exp\bigl( (\chi^\symmO_w)^p \bigr),
\qquad m=[\sqrt{2n} - 2w(2n)^{1/3}],\\
  \lim_{n\to\infty} \Exp\bigl( (\chi^\symmu_{n,m_+,m_-})^p \bigr)
&= \Exp\bigl( (\chi^\symmu_w)^p \bigr),
\qquad m_+=[\sqrt{n} - 2wn^{1/3}],\ m_-=[\sqrt{n}\beta].
\end{align}
\end{thm}

\begin{rem}
Theorem \ref{newthm1} shows that when $\alpha>1$ is fixed, we have used
incorrect scaling. When properly scaled, the resulting limiting distribution is
Gaussian. See Section \ref{sec:bigalpha} for the statement and the proof.
\end{rem}

The proof of Theorem \ref{newthm1}, \ref{newthm2}, \ref{newthm3} are provided
in Section \ref{sec:pf}.

In terms of the points selection process, which is a version of (directed site)
percolation, mentioned in the Introduction, the above results show that the
limiting distribution of the longest path depends on the geometry of the
domain, while the order of fluctuation is the same: $(\text{\it mean})^{1/3}$. From
\eqref{qz1.4}, the longest up/right path in a rectangle $0\le x,y\le 1$ has
$F_2$ in the limit, while the longest up/right path in a lower triangle $0\le
x<y\le 1$ has $F_4$, \eqref{qz3.4}, in the limit if there is no points on the
edge $0\le x=y\le 1$. If there are points on $0\le x=y\le 1$, they affect the
length of the longest up/right path. On the other hand, the longest down/right
path corresponding to $L^\symmS_{n,m}$ can be thought as the longest path from
the point $(0,1)$ to the \emph{line} $0\le x=y\le 1$. Thus the result
\eqref{qz3.7} shows that the point-to-line maximizing path has different
limiting distribution from the point-to-point maximizing path, $F_2$ from
\eqref{qz1.4}, though the fluctuation order is identical. One can also state
similar results for Poisson process (see Theorem \ref{thm11} below) and
also certain directed site percolation processes considered by
\cite{kurtj:shape} (see \cite{BR3}). This observation came from discussions
of one of us (J.B.) with Charles Newman to whom we are specially grateful.

\subsection{General involutions}

Now we consider general involutions and signed involutions without any
conditions on the number of fixed or negated points.

\begin{thm}\label{thm22}
For any fixed $x\in\R$, we have
\begin{align}
\label{as110}
  \lim_{n\to\infty} \Prob\biggl( \tilde{\chi}_n :=
\frac{\tilde{L}_n - 2\sqrt{n}}{n^{1/6}}\le x\biggr) &= F_1(x),\\
\label{as13-3}
   \lim_{n\to\infty} \Prob\biggl( \tilde{\chi}^\symmu_n :=
\frac{\tilde{L}^\symmu_n - 2\sqrt{2n}}{2^{2/3}(2n)^{1/6}}\le x\biggr)
&= F_1(x)^2.
\end{align}
Also for any $p=1,2,3,\cdots$,
\begin{align}
\label{as111}
  \lim_{n\to\infty} \Exp\bigl( (\tilde{\chi}_n)^p\bigr)
&= \Exp\bigl( (\chi_{\text{GOE}})^p \bigr),\\
\label{as13-4}
  \lim_{n\to\infty}\Exp\bigl( (\tilde{\chi}^\symmu_n)^p \bigr)
&= \Exp\bigl( (\chi_{\text{GOE}^2})^p \bigr).
\end{align}
\end{thm}

As mentioned in the introduction, this result proves that the first row of a
random Young diagram under the $1$-Plancherel measure $M^1_n$ behaves
statistically like the largest eigenvalue of a random GOE matrix as
$n\to\infty$. The proof of the theorem is given in Section \ref{as-asymp3}.

\subsection{Second rows}

For the second row, we scale the same way
as in \eqref{res1}--\eqref{res3}, and denote
the scaled random variables by $\chi^{\symmO,(2)}_{n,m}$,
$\chi^{\symmS,(2)}_{n,m}$
and $\chi^{\symmu,(2)}_{n,m_+,m_-}$, respectively.

\begin{thm}\label{thm-second1}
   Let $\alpha, \beta\ge 0$ be fixed. Then
 \begin{align}
   \lim_{n\to\infty}
\Prob\bigl(\chi^{\symmO,(2)}_{n,\sqrt{2n}\alpha}\bigr) &= F_4(x),\\
   \lim_{n\to\infty}
\Prob\bigl(\chi^{\symmS,(2)}_{n,[\sqrt{2n}\beta]}\bigr)
&= F_4(x),\\
   \lim_{n\to\infty}
\Prob\bigl(\chi^{\symmu,(2)}_{n,[\sqrt{n}\alpha],[\sqrt{n}\beta]}\bigr)
&= F_2(x),
 \end{align}
and for any $p=1,2,3,\cdots$,
 \begin{align}
   \lim_{n\to\infty}
\Exp\bigl((\chi^{\symmO,(2)}_{n,\sqrt{2n}\alpha})^p\bigr)
&= \Exp\bigl((\chi_{GSE})^p\bigr),\\
   \lim_{n\to\infty}
\Exp\bigl((\chi^{\symmS,(2)}_{n,[\sqrt{2n}\beta]})^p \bigr)
&= \Exp\bigl((\chi_{\text{GSE}})^p \bigr),\\
   \lim_{n\to\infty}
\Exp\bigl((\chi^{\symmu,(2)}_{n,[\sqrt{n}\alpha],[\sqrt{n}\beta]})^p \bigr)
&= \Exp \bigl((\chi_{\text{GUE}})^p\bigr).
 \end{align}
\end{thm}

This theorem is proved in Section \ref{sec:pf}.

As in the first row, these results yield the following theorem on the second
rows of general (signed) involutions. The proof is very similar to the proof of
Theorem \ref{thm22} and we skip the details.

\begin{thm}\label{thm-second2}
   For any fixed $x\in\R$, we have
 \begin{align}\label{e-t1}
   \lim_{n\to\infty}
\Prob\biggl(\tilde{\chi}^{(2)}_{n}:=
\frac{\tilde{L}^{(2)}_{n}-2\sqrt{n}}{n^{1/6}} \le x \biggr)
&= F_4(x), \\
\label{qz3.31}   \lim_{n\to\infty} \Prob\biggl(\tilde{\chi}^{\symmu,(2)}_{n}:=
\frac{\tilde{L}^{\symmu,(2)}_{n}-2\sqrt{2n}}{2^{2/3}(2n)^{1/6}}\le x \biggr) &=
F_2(x).
 \end{align}
Also for any $p=1,2,3,\cdots$,
 \begin{align}
      \lim_{n\to\infty}
\Exp\bigl((\tilde{\chi}^{(2)}_{n})^p \bigr)
&= \Exp\bigl((\chi_{\text{GSE}})^p \bigr),\\
   \lim_{n\to\infty}
\Exp\bigl((\tilde{\chi}^{\symmu,(2)}_{n})^p \bigr)
&= \Exp \bigl((\chi_{\text{GUE}})^p\bigr).
 \end{align}
\end{thm}

We conclude this section with some remarks on GOE and GSE. If the
conjecture given in the introduction that the $k^{th}$ row of a random
involution behaves in the limit like the $k^{th}$ largest eigenvalue of a
random GOE matrix were true, the result \eqref{e-t1} suggests that the
limiting distribution, $F_1^{(2)}$, of the second largest eigenvalue of GOE
is equal to the limiting distribution, $F_4$, of the largest eigenvalue of
GSE. Equivalently, since a GSE matrix has double eigenvalues, the second
eigenvalues of GOE and GSE are expected to have the same limiting
distribution. An indication for this is Theorem 10.6.1 of \cite{Mehta},
which says that the distributions of $N$ alternate angles of the
eigenvalues of a random $2N\times 2N$ matrix taken from the circular
orthogonal ensemble (COE) are identical to those of the $N$ angles of the
eigenvalues of a random $N\times N$ matrix taken from the circular
symplectic ensemble (CSE). Indeed for $2N\times 2N$ Laguerre ensembles, we
have proved that the joint distributions of the second, fourth, sixth,\dots
largest eigenvalues of the Laguerre orthogonal ensemble (LOE) and the
Laguerre symplectic ensemble (LSE) are identical (see Remark 1 to Corollary
7.6 of
\cite{PartI}). In particular, since the $k^{th}$ largest eigenvalue of a
Laguerre ensemble has the same limiting distribution as the corresponding
quantity for a Gaussian ensemble, the above remark implies that
\begin{equation}
  F_4^{((2k-1)}=F_4^{(2k)}=F_1^{(2k)}, \qquad k=1,2,\cdots.
\end{equation}
Thus \eqref{as110} and \eqref{e-t1} imply that the first and the second row of
a random involution has the same limiting distribution as the first and the
second eigenvalue of GOE, respectively.

Recently in \cite{ForRai99}, the authors proved that the
same property holds true for GOE and GSE. They also proved that, among many
things, the $(2k)^{th}$ `eigenvalue' of a superimposition of two random GOE
matrices has the same distribution as the $k^{th}$ eigenvalue of a random
GUE matrix. In particular, when $k=1$, this implies that
\begin{equation}
  (F_1^2)^{(2)}=F_2,
\end{equation}
and hence \eqref{as13-3} and \eqref{qz3.31} states that the first and the
second rows of a random signed involution have the same limiting
distribution as the first and the second `eigenvalues' of the
superimposition of two random GOE matrices, respectively.

%
%

\section{Poisson generating functions}\label{sec:sum}

We review the results from \cite{PartI} which we will need in the proof of
the theorems of Section \ref{sec:results}. As in \cite{PartI}, throughout
the paper, the notation $\tsymmg$ indicates an arbitrary member of the set
$\{\tsymmO,\tsymmS,\tsymmu\}$.

\begin{defn}\label{def-smml}
 We define the Poisson generating functions for the distributions
introduced above:
\begin{eqnarray}
\label{as12-38}
  Q^\symmO_l(\lambda_1,\lambda_2) &:=&
e^{-\lambda_1-\lambda_2} \sum_{n_1,n_2\ge 0}
\frac{\lambda_1^{n_1} \lambda_2^{n_2}}{n_1!n_2!}
\Prob\bigl( L^\symmO_{n_2,n_1}\le l\bigr), \\
\label{as12-39}
  Q^\symmS_l(\lambda_1,\lambda_2) &:=&
e^{-\lambda_1-\lambda_2} \sum_{n_1,n_2\ge 0}
\frac{\lambda_1^{n_1} \lambda_2^{n_2}}{n_1!n_2!}
\Prob\bigl( L^\symmS_{n_2,n_1}\le l\bigr), \\
\label{as12-40}
  Q^\symmu_l(\lambda_1,\lambda_2,\lambda_3) &:=&
e^{-\lambda_1-\lambda_2-\lambda_3} \sum_{n_1,n_2,n_3\ge 0}
\frac{\lambda_1^{n_1} \lambda_2^{n_2} \lambda_3^{n_3}}{n_1!n_2!n_3!}
\Prob\bigl( L^\symmu_{n_3,n_1,n_2}\le l\bigr).
\end{eqnarray}
\end{defn}

As in \cite{PartI}, let $\tilde{f}^\symmO_{nml}$ (resp.,
$\tilde{f}^\symmS_{nml}$) be the number of involutions on $n$ numbers with $m$
fixed points with no increasing (resp., decreasing) subsequence of length
greater than $l$. Thus
$\tilde{f}^\symmO_{(2n_2+n_1)n_1l}=\Prob(L^\symmO_{n_2,n_1}\le l) \cdot
|S_{n_2,n_1}|$, etc. Also let $\tilde{f}^\symmu_{nm_+m_-l}$ be the number of
signed involutions on $2n$ letters with $2m_+$ fixed points and $2m_-$ negated
points with no increasing subsequence of length greater than $l$ :
$\tilde{f}^\symmu_{(2n_3+n_1+n_2)n_1n_2l}=\Prob(L^\symmu_{n_3,n_1,n_2}\le
l)\cdot |S^\symmu_{n_3,n_1,n_2}|$. We also define
\begin{align}
\label{3.4}
P^\symmO_l(t;\alpha)
&:=
e^{-\alpha t-t^2/2}
\sum_{0\le n} \frac{t^n}{n!}
\sum_{0\le m}
\alpha^m \tilde{f}^\symmO_{nml},\\
\label{3.5}
P^\symmS_l(t;\beta)
&:=
e^{-\beta t-t^2/2}
\sum_{0\le n} \frac{t^n}{n!}
\sum_{0\le m}
\beta^m \tilde{f}^\symmS_{nml},\\
\label{3.6}
P^\symmu_l(t;\alpha,\beta)
&:=
e^{-\alpha t-\beta t-t^2}
\sum_{0\le n} \frac{t^n}{n!}
\sum_{0\le m_+,m_-}
\alpha^{m_+} \beta^{m_-} \tilde{f}^\symmu_{nm_+m_-l}.
\end{align}
Using $|S_{n,m}|=\frac{(2n+m)!}{n!m!2^n}$ and
$|S^\symmu_{n,m_+,m_-}|=\frac{(2n+m_++m_-)!}{n!m_+!m_-!}$, it is easy to check
that
\begin{eqnarray}
\label{as12-41}
  P^\symmO_l(t;\alpha) &=& Q^\symmO_l(\alpha t, t^2/2), \\
\label{as12-42}
  P^\symmS_l(t;\beta) &=& Q^\symmS_l(\beta t, t^2/2), \\
\label{as12-43}
  P^\symmu_l(t;\alpha,\beta) &=&
Q^\symmu_l(\alpha t, \beta t, t^2).
\end{eqnarray}
It turns out that the $P$-formulae in \eqref{3.4}-\eqref{3.6}
are useful for algebraic manipulations (see \cite{PartI}),
while the $Q$-formulae
\eqref{as12-38}-\eqref{as12-40} are well adapted to asymptotic analysis.

\bigskip
The following results from \cite{PartI} provide the starting point for our
analysis in this paper.
For
a nonnegative integer $k$, define $\pi_k(z;t)=z^k+\cdots$
to be the monic orthogonal
polynomial of degree $k$ with respect to the weight function
$\exp(t(z+1/z))dz/(2\pi i)$ on the unit circle.
Let the norm of $\pi_k(z;t)$ be
$N_k(t)$ :
\begin{equation}
   \int_\Sigma \pi_n(z;t)\overline{\pi_m(z;t)}
e^{t(z+1/z)}\frac{dz}{2\pi iz}
= N_n(t)\delta_{nm}.
\end{equation}
We note that all the coefficients of $\pi_n(z;t)$ are real.
Define
\begin{equation}
   \pi^*_n(z;t):= z^n\pi_n(z^{-1};t).
\end{equation}
Then

\begin{thm}[Collorary 4.3 and Collorary 2.7 of \cite{PartI}]\label{cor12}
For $\alpha,\beta \geq 0$,
\begin{align}
   P_{2l}^\symmO(t;\alpha) &= e^{-\alpha t-t^2/2}
{1\over 2}
\biggl\{ \bigl[\pi^*_{2l-1}(-\alpha;t)
-\alpha\pi_{2l-1}(-\alpha;t)\bigr] D_l^{--}(t)
+ \bigl[\pi^*_{2l-1}(-\alpha;t)
+\alpha\pi_{2l-1}(-\alpha;t)\bigr] D_{l-1}^{++}(t)
\biggr\}, \\
   P_{2l+1}^\symmO(t;\alpha) &= e^{-\alpha t-t^2/2}
{1\over 2}
\biggl\{ \bigl[\pi^*_{2l}(-\alpha;t)+\alpha\pi_{2l}(-\alpha;t)\bigr]
e^{t}D_l^{+-}(t)
+  \bigl[\pi^*_{2l}(-\alpha;t)-\alpha\pi_{2l}(-\alpha;t)\bigr]
e^{-t}D_{l}^{-+}(t) \biggr\},\\
\label{qz4.14}   P_{2l+1}^\symmS(t;\beta) &= e^{-t^2/2} D_l^{++}(t),\\
\label{qz4.15}   P_{2l+1}^\symmu(t;\alpha,\beta) &= e^{-\alpha t-t^2}
\pi^*_l(-\alpha;t) D_l(t),
\end{align}
where for any real $t\ge 0$, $D_l(t)$ and $D^{\pm\pm}_l(t)$ are certain
Toeplitz and Hankel determinants which in turn can be written as
\begin{eqnarray}
e^{-t^2}D_l(t) &=&\prod_{j\geq l} N_j(t)^{-1}, \\
e^{-t^2/2}D^{--}_l(t) &=&
\prod_{j\geq l} N_{2j+2}(t)^{-1}(1+\pi_{2j+2}(0;t)),\\
e^{-t^2/2}D^{++}_l(t) &=&
\prod_{j\geq l} N_{2j+2}(t)^{-1}(1-\pi_{2j+2}(0;t)), \\
e^{-t^2/2+t}D^{+-}_l(t) &=&
\prod_{j\geq l} N_{2j+1}(t)^{-1}(1-\pi_{2j+1}(0;t)),\\
e^{-t^2/2-t}D^{-+}_l(t) &=&
\prod_{j\geq l} N_{2j+1}(t)^{-1}(1+\pi_{2j+1}(0;t)).
\end{eqnarray}
\end{thm}

\begin{rem}
  The absence of $\beta$ on the right hand side of \eqref{qz4.14} is fairly
  simple to explain.  Observe that in the point selection model, the
  longest decreasing subsequence can always be chosen to be symmetric about
  the diagonal; moreover, any decreasing subsequence can contain at most
  one diagonal point.  Thus if the longest decreasing subsequence has $l$
  points, then removing the diagonal points will result in a longest
  decreasing subsequence with $2[l/2]$ points.  The independence from
  $\beta$ is thus special to $P^\symmS_l$ for $l$ odd; for $l$ even, we do
  indeed have $\beta$-dependence.  But by the monotonicity of $P^\symmS_l$
  in $l$, we only need \eqref{qz4.14} to compute the limiting distribution;
  in particular, the limiting distribution will not depend on $\beta$.
  A similar remark applies to \eqref{qz4.15}.
\end{rem}

As a special case:

\begin{thm}[Theorem 2.5 and Collorary 4.3 of \cite{PartI}]\label{thm2.5}
For $l\ge 0$, we have the following formulae:
{
\allowdisplaybreaks
\begin{align}
P^\symmO_{2l+2}(t;0) &= e^{-t^2/2}
\left[
D^{--}_{l+1}(t)+
D^{++}_{l}(t)\right]/2,\\
P^\symmO_{2l+1}(t;0) &= e^{-t^2/2}
\left[
e^tD^{+-}_l(t)+
e^{-t}D^{-+}_l(t)\right]/2,\\
P^\symmS_{2l}(t;0) &=
e^{-t^2/2} D^{++}_l(t),\\
   P_{2l}^\symmO(t;1) &= P_{2l}^\symmS(t;1)
= e^{-t-t^2/2} D_l^{-+}(t),\\
   P_{2l+1}^\symmO(t;1) &= P_{2l+1}^\symmS(t;1)
= e^{-t^2/2} D_l^{++}(t),\\
P^\symmu_{2l}(t;0,0) &= e^{-t^2} D_l(t),\\
   P_{4l+1}^\symmu(t;1,\beta) &= e^{-t-t^2}
D_l^{++}(t)D_l^{-+}(t),\\
   P_{4l+3}^\symmu(t;1,\beta) &= e^{-t-t^2}
D_l^{++}(t)D_{l+1}^{-+}(t).
\end{align}
}
Also
$
P^\symmO_0(t;0) = e^{-t^2/2} D^{--}_0(t)
= e^{-t^2/2}
$.
\end{thm}

For the second row, we define the Poisson generating functions
in a similar manner.
Then we have

\begin{thm}[Theorem 5.8 and Collorary 5.12 of \cite{PartI}]\label{thm:sum1}
For $\alpha,\beta\ge 0$,
\begin{align}
P_l^{\symmO,(2)}(t;\alpha) &= P^\symmO_l(t;0),\\
P_{2l+1}^{\symmS,(2)}(t;\beta) &= P^\symmO_{2l}(t;0),\\
P_{2l+1}^{\symmu,(2)}(t;\alpha,\beta) &= P^\symmu_{2l}(t;0,0).
\end{align}
\end{thm}

%
%

\section{Asymptotics of orthogonal polynomials}\label{as-op}

Let $\Sigma=\{z\in\C : |z|=1\}$ be the unit circle in the complex plane,
oriented counterclockwise. Set
\begin{equation}
   \psi(z;t):= e^{t(z+z^{-1})}.
\end{equation}
Let $\pi_n(z;t)=z^n+\cdots$ be the $n$-th monic orthogonal polynomial with
respect to the measure $\psi(z;t)dz/(2\pi iz)$ on the unit circle. From
Theorems \ref{cor12} and \ref{thm2.5}, in order to obtain the asymptotics of
the Poisson generating functions, we need the asymptotics, as $k,t\to\infty$,
of
\begin{equation}\label{qz5.2}
   N_{k}(t), \qquad \pi_k(z;t), \qquad \pi^{*}_{k}(z;t).
\end{equation}
In this section, we summarize the asymptotic results for these quantities.

Define the $2\times 2$ matrix-valued function of $z$ in $\C\setminus\Sigma$
by
\begin{equation}\label{as4}
    Y(z;k;t):=\begin{pmatrix}
    \pi_{k}(z;t)&
\int_{\Sigma}\frac{\pi_{k}(s;t)}{s-z}\frac{\psi(s;t)ds}{2\pi i s^{k}}\\
    -N_{k-1}(t)^{-1}\pi^{*}_{k-1}(z;t)& -N_{k-1}(t)^{-1}
\int_{\Sigma}\frac{\pi^{*}_{k-1}(s;t)}{s-z}\frac{\psi(s;t)ds}{2\pi i s^{k}}
       \end{pmatrix}, \qquad k\geq 1.
\end{equation}
Then $Y(\cdot\thinspace;k;t)$ solves the following
Riemann-Hilbert problem (RHP) (see Lemma 4.1 in \cite{BDJ}):
\begin{equation}\label{as-RHP}
 \begin{cases}
     Y(z;k;t) \quad\text{is analytic in}\quad z\in\C\setminus\Sigma,\\
     Y_+(z;k;t)=Y_-(z;k;t) \begin{pmatrix}
1&\frac{1}{z^{k}}\psi(z;t)\\0&1 \end{pmatrix},
\quad\text{on}\quad z\in\Sigma,\\
     Y(z;k;t) \biggl(\begin{smallmatrix} z^{-k}&0\\0&z^k
\end{smallmatrix} \biggr)
=I+O(\frac1{z}) \quad\text{as}\quad z\to\infty.
 \end{cases}
\end{equation}
Here the notation $Y_+(z;k)$ (resp., $Y_-$) denotes the limiting value
$\lim_{z'\to z}Y(z';k)$ with $|z'|<1$ (resp., $|z'|>1$). Note that $k$ and $t$
play the role of external parameters in the above RHP ; in particular, the term
$O(\frac1{z})$ does not imply a uniform bound in $k$ and $t$. One can easily
show that the solution of the above RHP is unique, hence \eqref{as4} is the
unique solution of the above RHP. This RHP formulation of orthogonal
polynomials on the unit circle given in \cite{BDJ} is an adaptation of a result
of Fokas, Its and Kitaev in \cite{FIK} where they considered orthogonal
polynomials on the real line.

From \eqref{as4}, the quantities in \eqref{qz5.2} are equal to
\begin{eqnarray}
\label{as7.7}   N_{k-1}(t)^{-1} &=& -Y_{21}(0;k;t), \\
\label{as7.8}  \pi_k(z;t) &=& Y_{11}(z;k;t), \\
\label{as7.9}  \pi^{*}_{k}(z;t) &=& z^kY_{11}(z^{-1};k;t) =
Y_{21}(z;k+1;t)(Y_{21}(0;k+1;t))^{-1}.
\end{eqnarray}
(For the other entries of $Y$, one can check directly from \eqref{as4} that
$Y_{12}(0;k;t)= N_k(t)$, $Y_{22}(0;k;t)=\pi_k(0;t)$.) Thus the asymptotic
analysis of the RHP \eqref{as-RHP} would yield the asymptotics of the above
quantities, and hence eventually the theorems in Section \ref{sec:results}. The
asymptotic analysis of the RHP \eqref{as-RHP} was conducted in \cite{BDJ} with
special interest on $Y_{21}(0;k;t)$. But as mentioned in the Introduction,
\cite{BDJ} controlled the solution $Y(z)$ to the RHP \eqref{as-RHP} in a
uniform way. In \cite{BDJ} and the theorem below, it is natural to distinguish
five different regimes of $k$ and $t$. From the analysis of \cite{BDJ}, the
following results for $Y(0;k;t)$, except for $\pi_k(0;t)$ in the case (ii)
below, can be directly read off. For example, (5.34)--(5.35) of \cite{BDJ} yield
the estimates for the case (iii) below when $x\ge 0$. For the case (ii), we
need to improve the $L^1$ norm bound, (5.23) in \cite{BDJ}, of the associated
jump matrix. If one is interested only on $N_{k-1}(t)$, the first integral
involving $w^{(3)}(s)$ in the displayed equation before (5.19) of \cite{BDJ}
vanishes, and hence the bound (5.23) of \cite{BDJ} was enough. But for
$\pi_k(0;t)$, this integral does not vanish, and we need an improved bound. See
the discussion in \eqref{as12.39}--\eqref{as12.41} below.

\begin{prop}[\cite{BDJ}]\label{thm1}
There exists $M_0>0$ such that
as $k,t\to\infty$, we have the following
asymptotic results for $N_{k-1}(t)$ and $\pi_k(0;t)$
in each different region of $k$ and $t$.
\begin{enumerate}
 \item If $0\le 2t \le ak$ with $0<a<1$, then
\begin{equation}
   \bigl|N_{k-1}(t)^{-1}-1\bigr| , \ \
\bigl|\pi_k(0;t)\bigr| \le Ce^{-ck},
\end{equation}
for some constants $C,c>0$.
 \item If $a k \le 2t \le k-Mk^{1/3}$ with some $M>M_0$
and $0<a<1$, then
\begin{equation}
   \bigl|N_{k-1}(t)^{-1}-1\bigr| , \ \
\bigl|\pi_k(0;t)\bigr|  \le
\frac{C}{k^{1/3}}e^{-\frac{2\sqrt{2}}{3}k
\bigl(1-\frac{2t}{k}\bigr)^{3/2}},
\end{equation}
for some constant $C>0$.
 \item If $2t = k-\frac{x}{2^{1/3}}k^{1/3}$ with
$-M \le x\le M$ for some constant $M>0$, then
\begin{equation}
      \biggl|N_{k-1}(t)^{-1}-1-
\frac{2^{1/3}}{k^{1/3}}v(x) \biggr| , \ \
\biggl|\pi_k(0;t)+ (-1)^k \frac{2^{1/3}}{k^{1/3}}u(x) \biggr|
\le \frac{C}{k^{2/3}},
\end{equation}
for some constant $C>0$, where $u(x)$ and $v(x)$ are defined in \eqref{as10}
and \eqref{as12} respectively.
 \item If $k+Mk^{1/3} \le 2t \le ak$ with some $M>M_0$ and $a>1$, then
\begin{equation}
   \biggl|\sqrt{\frac{2t}{k}}e^{k(2t/k-\log(2t/k)-1)}N_{k-1}(t)^{-1}- 1
\biggr|, \ \ \biggl| (-1)^k\sqrt{\frac{2t}{2t-k}}\pi_k(0;t) -1 \biggr|
\le \frac{C}{2t-k},
\end{equation}
for some constant $C>0$.
 \item If $ak \le 2t \le bk$ with $1<a<b$, then
\begin{equation}
   \biggl|\sqrt{\frac{2t}{k}}e^{k(2t/k-\log(2t/k)-1)}N_{k-1}(t)^{-1}- 1
\biggr|, \ \ \biggl| (-1)^k\sqrt{\frac{2t}{2t-k}}\pi_k(0;t) -1 \biggr|
\le \frac{C}{k},
\end{equation}
for some constant $C>0$.
\end{enumerate}
\end{prop}

Now we are interested in $\pi_{k}(z;t)$. If $z$ is apart from $-1$ and is
fixed, then similar estimates for $Y(z;k;t)$ can be obtained from the
analysis of \cite{BDJ}. The result below when $x\ge 0$ is (almost) direct
from the work of
\cite{BDJ}. For the case when $x<0$, the analysis of \cite{BDJ} expresses the
bound in terms of the so-called $g$-function, and we need further analysis for
this $g$-function. When $z=0$, this $g$-function becomes very simple:
$g(0)=\pi i$. See \eqref{as12-115}-\eqref{as13.136} below.

\begin{prop}\label{thm2}
  For  $2t=k-x(k/2)^{1/3}$, $x$ fixed,
and for each fixed $z\in\C\setminus\Sigma$, we have
\begin{eqnarray}
\label{as27}  \lim_{k\to\infty} e^{tz}\pi_k(z;t) = 0,
&\lim_{k\to\infty} e^{tz}\pi^*_k(z;t) = 1,& \quad |z|<1,\\
\label{as28}  \lim_{k\to\infty} z^{-k}e^{tz^{-1}}\pi_k(z;t) = 1,
&\lim_{k\to\infty} z^{-k}e^{tz^{-1}}\pi^*_k(z;t) = 0,&\quad  |z|>1.
\end{eqnarray}
\end{prop}

\begin{cor}\label{cor4}
  For $2t=k-x(k/2)^{1/3}$, $x$ fixed,
we have for fixed $\alpha>1$,
\begin{equation}
   \lim_{k\to\infty} e^{-\alpha t}\pi_k(-\alpha;t)=0, \qquad
\lim_{k\to\infty} e^{-\alpha t}\pi^*_k(-\alpha;t)=0.
\end{equation}
\end{cor}

\begin{proof}
  Write
\begin{equation}
 \begin{split}
    e^{-\alpha t}\pi_k(-\alpha;t)
&= \alpha^ke^{t(-\alpha+\alpha^{-1})}
\alpha^{-k}e^{-t\alpha^{-1}} \pi_k(-\alpha;t) \\
&= e^{kf(\alpha;2t/k)} \alpha^{-k}e^{-t\alpha^{-1}} \pi_k(-\alpha;t),
 \end{split}
\end{equation}
where
\begin{equation}
  f(\alpha;\gamma) = \frac{\gamma}2(-\alpha+\alpha^{-1}) + \log\alpha.
\end{equation}
The function $f(\alpha;1)$ is strictly decreasing for $\alpha>0$, and
$f(1;1)=0$.
Hence $f(\alpha;1)<0$ for $\alpha>1$.
Note that
\begin{equation}
   f(\alpha;\gamma) = f(\alpha;1) + \frac{\gamma-1}{2} (-\alpha+\alpha^{-1}).
\end{equation}
When $x\le 0$, $2t/k\geq 1$, hence $f(\alpha;2t/k)\le f(\alpha;1)$. On the
other hand, when $x>0$, since $2t/k-1=-x/(2^{1/3}k^{2/3})$, we have
$f(\alpha;2t/k)\le \frac12f(\alpha;1)$ if $k>\biggl(
\frac{2^{2/3}x(-\alpha+\alpha^{-1})}{f(\alpha;1)} \biggr)^{3/2}$. Therefore
\eqref{as28} implies that
\begin{equation}
  \bigl| e^{-\alpha t}\pi_k(-\alpha;t) \bigr| \le
e^{\frac{k}2f(\alpha;1)} \bigl|
(-\alpha)^{-k}e^{-t\alpha^{-1}} \pi_k(-\alpha;t)\bigr|
\to 0,
\end{equation}
as $k\to\infty$.
Similar calculations give the desired result for $\pi^*_k$.
\end{proof}

When $z\to -1$, (which is required for the proof of Theorem \ref{newthm2}), the
estimates for $Y(z;k;t)$ can not be directly read off from the result of
\cite{BDJ}. However, with more detailed estimates, the same procedure as in
\cite{BDJ} gives us the following restuls. See Subsubsection \ref{subsub3}
\eqref{as13.72} - \eqref{as13.82} for the case $x\ge 0$, and Subsubsection
\ref{subsub6} \eqref{as12-109} - \eqref{as13.153} for the case $x<0$. Recall
from Section \ref{sec:dist} that $m(\cdot\thinspace,x)$ solves the RHP for the
PII equation \eqref{as20}.

\begin{prop}\label{thm3}
  Let $2t=k-x(k/2)^{1/3}$ where $x$ is a fixed number.
Set
\begin{equation}
  \alpha = 1- \frac{2^{4/3}w}{k^{1/3}}.
\end{equation}
We have for $w>0$ fixed,
\begin{eqnarray}
\label{as7.32}
   \lim_{k\to\infty} (-1)^ke^{-t\alpha}\pi_{k}(-\alpha;t)
&=& -m_{12}(-iw;x) ,\\
\label{as7.33}
   \lim_{k\to\infty} e^{-t\alpha}\pi^*_{k}(-\alpha;t)
&=& m_{22}(-iw;x),
\end{eqnarray}
and for $w<0$ fixed,
\begin{eqnarray}
\label{as7.34}
   \lim_{k\to\infty} (-\alpha)^{-k}e^{-t\alpha^{-1}}\pi_{k}(-\alpha;t)
&=& m_{11}(-iw;x) ,\\
\label{as7.35}
   \lim_{k\to\infty} \alpha^{-k}e^{-t\alpha^{-1}}\pi^*_{k}(-\alpha;t)
&=& -m_{21}(-iw;x).
\end{eqnarray}
\end{prop}

\begin{cor}\label{cor5}
   For $w<0$, under the same condition as the above proposition,
we have
\begin{eqnarray}
   \lim_{k\to\infty} (-1)^{-k}e^{-t\alpha}\pi_{k}(-\alpha;t)
&=& m_{11}(-iw;x)e^{(8/3)w^3-2xw}, \\
   \lim_{k\to\infty} e^{-t\alpha}\pi^*_{k}(-\alpha;t)
&=& -m_{21}(-iw;x)e^{(8/3)w^3-2xw}.
\end{eqnarray}
\end{cor}

\begin{proof}
  Note that under the stated conditions, we have
\begin{equation}
  e^{t(\alpha^{-1}-\alpha)}\alpha^k
= e^{(8/3)w^3-2xw + O( k^{-1/3})}.
\end{equation}
\end{proof}

\begin{rem}
  As noted in Section \ref{sec:dist}, it follows from the RHP \eqref{as20}
that $(m_{12})_+(0;x)= -(m_{11})_-(0;x)$ and $(m_{22})_+(0;x)=
-(m_{21})_-(0;x)$, and hence by the above Corollary, the limits in
\eqref{as7.32}, \eqref{as7.33} are in fact continuous across $w=0$.
\end{rem}

For convergence of moments, we need uniform bound of $\pi_k(z;t)$ for
$|x|\ge M$ for a fixed number $M>0$. The results \eqref{as7.42} and
\eqref{as7.43} are essentially in the analysis of \cite{BDJ}, while
\eqref{as7.44}-\eqref{as7.47} are new estimates. We again need to extend the
method of \cite{BDJ} to obtain the results below. The proof is provided in
Section \ref{RHP}. See Subsubsection \ref{subsub1} and \ref{subsub2} for the
case $x\ge M$, and Subsubsection \ref{subsub4} and \ref{subsub5} for the case
$x\le -M$.

\begin{prop}\label{thm4}
Define $x$ through the relation
\begin{equation}
   \frac{2t}{k}= 1-\frac{x}{2^{1/3}k^{2/3}}.
\end{equation}
Then there exists $M_0$ such that the following holds for any fixed $M>M_0$.
Let $0<b<1$ and $0<L<2^{-3/2}\sqrt{M}$ be fixed. Then as $k,t\to\infty$, we
have for $x\ge M$,
\begin{eqnarray}
\label{as7.42}
  \bigl| e^{tz}\pi_k(z;t) \bigr| \le Ce^{-c|x|^{3/2}},  && |z|\le b, \\
\label{as7.43}
  \bigl| e^{tz^{-1}}z^{-k}\pi_k(z;t)-1 \bigr| \le Ce^{-c|x|^{3/2}}, &&
|z|\ge b^{-1}, \\
\label{as7.44}
  \bigl| e^{-t\alpha}\pi_k(-\alpha;t) \bigr| \le Ce^{c|x|},&&
\alpha=1-2^{4/3}k^{-1/3}w, \quad  -L\le w\le L, \\
\label{as7.45}  \bigl| e^{-t\alpha^{-1}}(-\alpha)^{-k}
\pi_k(-\alpha;t) -1 \bigr|
\le Ce^{-c|x|^{3/2}},&&  \alpha=1-2^{4/3}k^{-1/3}w, \quad -L\le w\le L,
\end{eqnarray}
and for $x\le -M$,
\begin{eqnarray}
\label{as7.46}
  \bigl| e^{-t\alpha}\pi_k(-\alpha;t) \bigr| \le C,&&
0< \alpha\le 1,\\
\label{as7.47}
   \bigl| e^{-t\alpha^{-1}}(-\alpha)^{-k}\pi_k(-\alpha;t) \bigr|
\le C,&&  \alpha\ge 1.
\end{eqnarray}
\end{prop}

\begin{cor}\label{cor-new}
  Let $\alpha=1-2^{4/3}wk^{-1/3}$ and $-L\le w\le L$ for fixed $L>0$.
  Under the assumption of the above proposition, for $x\le -M$,
we have
\begin{eqnarray}
\label{as7.48}
  \bigl| e^{-t\alpha}\pi_k(-\alpha;t) \bigr| &\le& Ce^{c|x|},\\
\label{as7.49}
  \bigl| e^{-t\alpha^{-1}}(-\alpha)^{-k}\pi_k(-\alpha;t) \bigr|
&\le& Ce^{c|x|}.
\end{eqnarray}
for some positive constants $C$, $c$.
\end{cor}

\begin{proof}
   We have
\begin{eqnarray}
\label{as7-48}
  |e^{-t(\alpha-\alpha^{-1})}\alpha^{k}| &=&e^{2xw+\frac83w^3+O(k^{-1})}, \\
\label{as7-49}
  |e^{t(\alpha-\alpha^{-1})}\alpha^{-k}| &=&e^{-2xw-\frac83w^3+O(k^{-1})}.
\end{eqnarray}
Above proposition shows that \eqref{as7.48} is true for $w\ge 0$.
For $w<0$, write
\begin{equation}
   e^{-t\alpha}\pi_k(-\alpha;t)
= \biggl[e^{-t(\alpha-\alpha^{-1})}(-\alpha)^{k} \biggr]
\biggl[ e^{-t\alpha^{-1}}(-\alpha)^{-k}\pi_k(-\alpha;t)\biggr].
\end{equation}
Now \eqref{as7.48} follows from \eqref{as7.47} and
\eqref{as7-48}.
The estimate \eqref{as7-49} is proved similarly.
\end{proof}

The result below is new and is used for the asymptotics of
$L^\symmO_{n,[\sqrt{2n}\alpha]}$ when $\alpha>1$. See Subsection \ref{sub3} for
the proof.

\begin{prop}\label{prop7.7}
  Let $\alpha>1$ be fixed.
When
\begin{equation}
   \frac{t}{k} = \frac{\alpha}{\alpha^2+1}
-\frac{\alpha(\alpha^2-1)^{1/2}}{(\alpha^2+1)^{3/2}}
\cdot\frac{x}{\sqrt{k}}, \qquad \text{$x$ fixed,}
\end{equation}
we have
\begin{equation}
  \lim_{k\to\infty} e^{-\alpha t}(-\alpha)^{k}\pi_{k}(-\alpha^{-1};t)
= \frac1{\sqrt{2\pi}} \int_{-\infty}^x e^{-\frac12y^2}dy.
\end{equation}
\end{prop}

%
%

\section{De-Poissonization lemmas}\label{as-dePoi}

In this section, we describe a series of Tauberian type de-Poissonization
lemmas, which enable us to extract the asymptotics of the coefficient from the
knowledge of the asymptotics of its generating function. Lemma \ref{lem10}
below is due to Johansson \cite{Jo2} and Lemma \ref{lem11} is taken from
Section 8 in \cite{BDJ}. Lemma \ref{lem19} and Lemma \ref{lem20} are
multi-index versions. Lemma \ref{lem10} and \ref{lem19} are enough for both
convergence in distribution and the convergence of moments, but for convenience
of computations, we use Lemma \ref{lem11} and \ref{lem20} for the convergence
of moments in the subsequent sections.

For a sequence $q=\{q_n\}_{n\ge 0}$, we define its Poisson generating
function by
\begin{equation}
\phi(\lambda) = e^{-\lambda} \sum_{0\le n} q_n
\frac{\lambda^n}{n!}.
\end{equation}

\begin{lem}\label{lem10}
For any fixed real number $d>0$, set
\begin{eqnarray}
\mu^{(d)}_n &=& n+(2\sqrt{d+1}+1)\sqrt{n\log n},\\
\nu^{(d)}_n &=& n-(2\sqrt{d+1}+1)\sqrt{n\log n}.
\end{eqnarray}
Then there are constants $C$ and $n_0$ such that
for any sequence $q=\{q_n\}_{n\ge 0}$ satisfying
(i) $q_n\ge q_{n+1}$ (ii) $0\le q_n\le 1$, for all $n\ge 0$,
\begin{equation}
\phi(\mu^{(d)}_n) - Cn^{-d} \le q_n \le \phi(\nu^{(d)}_n)+Cn^{-d}
\end{equation}
for all $n\ge n_0$.
\end{lem}

\begin{lem}\label{lem11}
  For any fixed real number $d>0$, there exist constants
$C$ and $n_0$ such that
for any sequence $q=\{q_n\}_{n\ge 0}$ satisfying (i) and (ii) above,
\begin{eqnarray}
   q_n &\le& C\phi(n-d\sqrt{n}), \\
  1-q_n &\le& C(1-\phi(n+d\sqrt{n}))
\end{eqnarray}
for all $n\ge n_0$.
\end{lem}

For multi-indexed sequences, there are similar results.
For $q=\{q_{n_1,n_2}\}_{n_1,n_2\ge 0}$, define
\begin{equation}
\phi(\lambda_1,\lambda_2) =
e^{-\lambda_1-\lambda_2}
\sum_{n_1,n_2\ge 0}
q_{n_1n_2}
\frac{\lambda_1^{n_1} \lambda_2^{n_2}}{n_1!n_2!}.
\end{equation}
From the above two lemmas, we easily obtain the following lemmas.

\begin{lem}\label{lem19}
For any real number $d>0$, define $\mu^{(d)}_n$ and $\nu^{(d)}_n$
as in Lemma \ref{lem10}.
Then there exist constants
$C$ and $n_0$ such that
for any $q=\{q_{n_1,n_2}\}_{n_1,n_2\ge 0}$ satisfying
(i) $q_{n_1,n_2}\ge q_{n_1+1,n_2}$, $q_{n_1,n_2}\ge q_{n_1.n_2+1}$
(ii) $0\le q_{n_1,n_2} \le 1$, for all $n_1,n_2\ge 0$,
\begin{equation}
\phi(\mu^{(d)}_{n_1},\mu^{(d)}_{n_2})
-
C(n_1^{-d}+n_2^{-d})
\le
q_{n_1n_2}
\le
\phi(\nu^{(d)}_{n_1},\nu^{(d)}_{n_2})
+
C(n_1^{-d}+n_2^{-d})
\end{equation}
for all $n_1,n_2\ge n_0$.
\end{lem}

Similarly,
\begin{lem}\label{lem20}
For any fixed real number $d>0$, there exist constants
$C$ and $n_0$ such that for any $q=\{q_{n_1,n_2}\}_{n_1,n_2\ge 0}$
satisfying two condition in Lemma \ref{lem19},
\begin{eqnarray}
  q_{n_1n_2} &\le& C\phi(n_1-d\sqrt{n_1}, n_2-d\sqrt{n_2}),\\
  1-q_{n_1n_2} &\le& C (1-\phi(n_1-d\sqrt{n_1}, n_2+d\sqrt{n_2}))
\end{eqnarray}
for $n_1,n_2\ge n_0$.
\end{lem}

\begin{rem}
Similar lemmas hold true for sequences of arbitrarily many indices.
\end{rem}

\section{Proofs of Theorems \ref{newthm1}, \ref{newthm2}, \ref{newthm3}
and \ref{thm-second1}}\label{sec:pf}

The following results follow from Proposition \ref{thm1} above. The result
\eqref{as50} is derived in Lemma 7.1 (iii) of \cite{BDJ}, and the other cases
are similar. We omit the details.

\begin{cor}\label{cor6}
  Let $M>M_0$, where $M_0$ is given in Proposition \ref{thm4}.
Then there exist positive constants $C$ and $c$ which
are independent of $M$, and a positive constant $C(M)$ which
may depend on $M$, such that the following results hold for large $l$.

(i)\cite{BDJ} Define $x$ by $2t=l-x(l/2)^{1/3}$. For $-M< x <M$,
\begin{equation}\label{as50}
  \biggl| \sum_{j\geq l} \log N_j(t)^{-1} - 2\log F(x) \biggr|
\le \frac{C(M)}{l^{1/3}} + Ce^{-cM^{3/2}}.
\end{equation}

(ii) Define $x$ by $t=l-(x/2)l^{1/3}$.
For $-M<x<M$,
\begin{align}
  \biggl| \sum_{j\geq l} \log N_{2j+2}(t)^{-1} - \log F(x) \biggr| , \qquad
\biggl| \sum_{j\geq l} \log N_{2j+1}(t)^{-1} - \log F(x) \biggr|
&\le& \frac{C(M)}{l^{1/3}} + Ce^{-cM^{3/2}}, \\
  \biggl| \sum_{j\geq l} \log (1-\pi_{2j+2}(0;t)) - \log E(x)  \biggr| ,\ \
\biggl| \sum_{j\geq l} \log (1+\pi_{2j+1}(0;t)) - \log E(x)  \biggr|
&\le& \frac{C(M)}{l^{1/3}} + Ce^{-cM^{3/2}}, \\
   \biggl| \sum_{j\geq l} \log (1+\pi_{2j+2}(0;t)) + \log E(x)  \biggr| , \ \
\biggl| \sum_{j\geq l} \log (1-\pi_{2j+1}(0;t)) + \log E(x)  \biggr|
&\le& \frac{C(M)}{l^{1/3}} + Ce^{-cM^{3/2}}.
\end{align}
\end{cor}

These results yield the asymptotics of the determinants
in Theorem \ref{cor12}.

\begin{cor}\label{cor7}
  There exits $M_1$ such that for $M>M_1$,
there exist positive constants $C$ and $c$ which
are independent of $M$, and a positive constant $C(M)$ which
may depend on $M$ such that the following results hold for large $l$.

(i) Define $x$ by $2t=l-x(l/2)^{1/3}$.  For $-M< x <M$,
\begin{equation}
\label{as11-14}
   \bigl| e^{-t^2}D_l(t) - F(x)^2 \bigr| \le
\frac{C(M)}{l^{1/3}} + Ce^{-cM^{3/2}}.
\end{equation}

(ii) Define $x$ by $t=l-(x/2)l^{1/3}$. For $-M< x <M$,
\begin{eqnarray}
\label{as11-15}
  \bigl| e^{-t^2/2}D^{--}_l(t) - F(x)E(x)^{-1} \bigr| &\le&
\frac{C(M)}{l^{1/3}} + Ce^{-cM^{3/2}}, \\
\label{as11-16}
  \bigl| e^{-t^2/2}D^{++}_{l-1}(t) - F(x)E(x) \bigr| &\le&
\frac{C(M)}{l^{1/3}} + Ce^{-cM^{3/2}}, \\
\label{as11-17}
  \bigl| e^{-t^2/2+t}D^{+-}_l(t) - F(x)E(x)^{-1} \bigr| &\le&
\frac{C(M)}{l^{1/3}} + Ce^{-cM^{3/2}}, \\
\label{as11-18}
  \bigl| e^{-t^2/2-t}D^{-+}_l(t) - F(x)E(x) \bigr| &\le&
\frac{C(M)}{l^{1/3}} + Ce^{-cM^{3/2}}.
\end{eqnarray}
\end{cor}

\begin{proof}
  For $C$ and $c$ in the above corollary,
take $M_1>M_0$ such that $Ce^{-cM_1^{3/2}}\le \frac12$. Once we fix $M>M_1$,
then for $l$ is large, $\frac{C(M)}{l^{1/3}} + Ce^{-cM^{3/2}} <1$, and hence by
\eqref{as50} above, $\bigl|\sum_{j\geq l} \log N_j(t)^{-1} - 2\log F(x)\bigr|
\le 1$. Using $|e^x -1 | \le (e-1)|x|$ for $|x|\le 1$,
\begin{equation}
 \begin{split}
   \bigl| e^{-t^2}D_l(t) - F(x)^2  \bigr|
&=  F(x)^2 \bigl|
e^{\bigl(\sum_{j\geq l} \log N_j(t)^{-1} - 2\log F(x)\bigr)} -1 \bigr|\\
&\le (e-1)F(x)^2 \biggl|\sum_{j\geq l} \log N_j(t)^{-1} - 2\log F(x)\biggr|.
 \end{split}
\end{equation}
But from \eqref{as46} and \eqref{as48}, $F(x)$ is bounded for $x\in\R$.
Hence using \eqref{as50}, we obtain the result for (i)
with new constants $C$, $c$ and $C(M)$.
For (ii), we note that $F(x)E(x)$ and $F(x)E(x)^{-1}$ are bounded
for $x\in\R$ from \eqref{as46}-\eqref{as49}.
\end{proof}

From Proposition \ref{thm2}, Corollary \ref{cor4} and
Theorems \ref{cor12}, \ref{thm2.5},
this Corollary
immediately yields the following asymptotics for
Poisson generating functions.

\begin{prop}\label{thm11}
  Let $2t=l-x(l/2)^{1/3}$ where $x$ is fixed.
As $l\to\infty$, for each fixed $\alpha,\beta$,
\begin{eqnarray}
\label{as57}   P_{l}^\symmO(t;\alpha) &\to& F_4(x), \qquad 0\le\alpha<1,\\
\label{as81}   P_{l}^\symmO(t;1) &\to& F_1(x),\\
\label{as82}   P_{l}^\symmO(t;\alpha) &\to& 0, \qquad\qquad \alpha>1,\\
\label{as83}   P_{l}^\symmS(t;\beta) &\to& F_1(x), \qquad \beta\geq 0.
\end{eqnarray}
Let $4t=l-x(2l)^{1/3}$ where $x$ is fixed.
As $l\to\infty$, for each fixed $\alpha,\beta$,
\begin{eqnarray}
\label{as84}   P_{l}^\symmu(t;\alpha,\beta) &\to& F_2(x), \qquad
0\le\alpha<1,\medspace \beta\geq 0,\\
\label{as84.5}   P_{l}^\symmu(t;1,\beta) &\to& F_1(x)^2, \qquad
\beta\geq 0,\\
\label{as85}   P_{l}^\symmu(t;\alpha,\beta) &\to& 0, \qquad\qquad
\alpha>1,\ \ \beta\geq 0.
\end{eqnarray}
\end{prop}

Similarly, using
Proposition \ref{thm3} and Corollary \ref{cor5}
and Theorem \ref{cor12},
we have:

\begin{thm}\label{thm12}
Let $2t=l-x(l/2)^{1/3}$ where $x$ is fixed.
As $l\to\infty$, we have for any fixed $w\in\R$,
\begin{equation}
   P_{l}^\symmO(t;\alpha) \to F^\symmO(x;w),
\qquad \alpha = 1- \frac{2^{4/3}w}{l^{1/3}}.
\end{equation}
Let $4t=l-x(2l)^{1/3}$ where $x$ is fixed.
As $l\to\infty$, we have for each fixed $\beta$ and $w\in\R$,
\begin{equation}
   P_{l}^\symmu(t;\alpha,\beta) \to F^\symmu(x;w),
\qquad \alpha = 1- \frac{2^{5/3}w}{l^{1/3}}, \medspace \beta\ge 0.
\end{equation}
\end{thm}

Recall the relation between $Q^\symmg_l(\lambda)$ and $P^\symmg_l(t)$ in
\eqref{as12-41}-\eqref{as12-43}. We now use the de-Poissonization Lemma
\ref{lem19} to obtain the asymptotic results of Theorem \ref{newthm1} and
\ref{newthm2}. In order to apply the de-Poissonization Lemma, we need the
following monotonicity results.

\begin{lem}[Monotonicity]
  For any $l$, $\Prob(L^\symmO_{k,m}\le l)$, $\Prob(L^\symmS_{k,m}\le l)$
and $\Prob(L^\symmu_{k,m_+,m_-}\le l)$ are monotone decreasing
in $k$, $m$, $m_+$ and $m_-$.
\end{lem}

\begin{proof}
  We first consider $\Prob(L^\symmO_{k,m}\le l)$.
Let $f_{km}:=\Prob(L^\symmO_{k,m}\le l)\cdot |S_{k,m}|$ be the number of
elements in $S_{k,m}$ with no increasing subsequence greater than $l$. Consider
the map $h : S_{k,m-1}\times \{1,2,\cdots,2k+m\} \to S_{k,m}$ defined as
follows: for $(\pi,j)\in S_{k,m-1}\times \{1,2,\cdots, 2k+m\}$, set
$h(\pi,j)(x)=\pi(x)$ for $1\le x<j-1$, $h(\pi,j)(j)=j$, and
$h(\pi,j)(x)=\pi(x-1)$ for $j<x\le 2k+m$. Then it is easy to see that
$h^{-1}(\sigma)$ consists of $m$ elements, hence
$(2k+m)|S_{k,m-1}|=m|S_{k,m}|$. Moreover if $\pi\in S_{k,m-1}$ has an
increasing subsequence of length greater than $l$, then $h(\pi,j)$ has an
increasing subsequence of length greater than $l$. Thus $(2k+m)f_{k(m-1)}\ge
mf_{km}$. But since $|S_{k,m}|=\frac{(2k+m)!}{2^kk!m!}$, we obtain
$\Prob(L^\symmO_{k,m-1}\le l)\ge \Prob(L^\symmO_{k,m}\le l)$.

A similar argument works for the other cases.
Note that $|S^\symmu_{k,m_+,m_-}|=\frac{(2k+m_+,m_-)!}{k!m_+!m_-!}$.
\end{proof}

Thus Lemma \eqref{lem19} can be applied to obtain the asymptotics results in
Theorems \ref{newthm1} and \ref{newthm2}. The proofs are similar to that in
Section 9 of \cite{BDJ}.

\bigskip

Now we consider convergence of moments.
For this, we first obtain the following estimates which follow
from Proposition \ref{thm1} (i),(ii),(iv),(v) above.
The proof is very similar to that of Lemma 7.1 (i),(ii),(iv),(v)
of \cite{BDJ}.
Compare the results with \eqref{as46}-\eqref{as49}
noting Corollary \ref{cor6}.

\begin{cor}\label{cor11}
   Set
\begin{equation}
   t = l-\frac{x}{2}l^{1/3}.
\end{equation}
There exits $M_2$ such that for a fixed $M>M_2$,
there are positive constants $C=C(M)$ and $c=c(M)$
such that the following results hold.
\begin{enumerate}
  \item For $x\ge M$,
\begin{eqnarray}
1-\prod_{j\geq l} N_{2j+2}(t)^{-1},\ \
1-\prod_{j\geq l} N_{2j+1}(t)^{-1} &\le& Ce^{-c|x|^{3/2}},\\
1-\prod_{j\geq l} (1- \pi_{2j+2}(0;t)),\ \
1-\prod_{j\geq l} (1+ \pi_{2j+1}(0;t)) &\le& Ce^{-c|x|^{3/2}},\\
1-\prod_{j\geq l} (1+ \pi_{2j+2}(0;t)),\ \
1-\prod_{j\geq l} (1- \pi_{2j+1}(0;t)) &\le& Ce^{-c|x|^{3/2}},
\end{eqnarray}
  \item For $x\le -M$,
\begin{eqnarray}
\prod_{j\geq l} N_{2j+2}(t)^{-1},\ \
\prod_{j\geq l} N_{2j+1}(t)^{-1} &\le& Ce^{-c|x|^{3}}, \\
\prod_{j\geq l} (1- \pi_{2j+2}(0;t)),\ \
\prod_{j\geq l} (1+ \pi_{2j+1}(0;t)) &\le& Ce^{- c|x|^{3/2}}, \\
\prod_{j\geq l} (1+ \pi_{2j+2}(0;t)),\ \
\prod_{j\geq l} (1- \pi_{2j+1}(0;t)) &\le& Ce^{+ c|x|^{3/2}}.
\end{eqnarray}
\end{enumerate}
\end{cor}

\begin{rem}
 From the definitions of $P^\symmg_l(t)$ and the equalities of Theorem
\ref{cor12}, we know that all the infinite products above are
between $0$ and $1$.
\end{rem}

Now as in Section 9 of \cite{BDJ}, using
Lemma \ref{lem20} and Theorems \ref{newthm1}, \ref{newthm2},
this implies Theorem \ref{newthm3}.

\medskip
Theorem \ref{thm-second1} follows from
Theorem \ref{thm:sum1}.

\section{Proofs of Theorems \ref{thm22} and \ref{thm-second2}}
\label{as-asymp3}

In this section, we prove Theorem \ref{thm22} by summing up the asymptotic
results of Theorems \ref{newthm1}, \ref{newthm2} and \ref{newthm3}. Theorem
\ref{thm-second2} can be proved in a similar way from Theorem \ref{thm:sum1}.

\begin{proof}[{\bf Proof of \eqref{as110}}]


Note that we have a disjoint union
\begin{equation}
  \tilde{S}_n = \bigcup_{2k+m=n} S_{k,m}.
\end{equation}
Set $p^\symmS_{kml}=\Prob( L^\symmS_{km}\le l)$, the probability that
the length of the longest decreasing subsequence of $\pi\in S_{k,m}$
is less than or equal to $l$.
As the first row and the first column of $\pi$ in $\tilde{S}_n$
have the same statistics, we have
\begin{equation}\label{as112}
  \Prob( \tilde{L}_n \le l)
= \frac1{|\tilde{S}_n|}
\sum_{2k+m=n} p^\symmS_{kml} |S_{k,m}|.
\end{equation}
Note that
\begin{equation}\label{as13.7}
  |S_{k,m}|
= \binom{2k+m}{2k} \frac{(2k)!}{2^{k}k!}.
\end{equation}
As $n\to\infty$ (see pp.66-67 of \cite{Kn}),  we have
\begin{equation}\label{as113}
  |\tilde{S}_n|
= \sum_{2k+m=n} |S_{k,m}|
= \frac1{\sqrt{2}} n^{n/2} e^{-n/2+\sqrt{n}-1/4}
\biggl( 1+\frac{7}{24}n^{-1/2} + O\bigl(n^{-3/4}\bigr)\biggr),
\end{equation}
and the main contribution to the sum comes from
$\sqrt{n}-n^{\epsilon+1/4} \le m \le \sqrt{n}+n^{\epsilon+1/4}$.

Fix $0<a<1<b$.
We split the sum in \eqref{as112} into two pieces :
\begin{equation}\label{as114}
   \Prob( \tilde{L}_n \le l) =
\frac1{|\tilde{S}_n|} \biggl[
\sum_{(*)}p^\symmS_{kml} |S_{k,m}|
+ \sum_{(**)}p^\symmS_{kml} |S_{k,m}| \biggr],
\end{equation}
where $(*)$ is the region $a\sqrt{n}\le m\le b\sqrt{n}$
and $(**)$ is the rest.

For $2k+m=n$, the quantity $|S_{k,m}|=\binom{n}{2k}\frac{(2k)!}{2^kk!}$
is unimodal for $0\le k\le n$,
and the maximum is achieved when $k\sim (n-\sqrt{n})/2$ as $n\to\infty$.
Hence
\begin{equation}
    \sum_{(**)}p^\symmS_{kml} |S_{k,m}|
\le n\cdot \max \bigl( |S_{k,[a\sqrt{n}]}|, \ \
|S_{k,[b\sqrt{n}]}|
\bigr).
\end{equation}
Using Stirling's formula for \eqref{as13.7},
for any fixed $c$, when $2k+[c\sqrt{n}]=n$,
\begin{equation}
   |S_{k,[c\sqrt{n}]}| \sim
n^{n/2}e^{-n/2+\sqrt{n}(c-c\log c)}
\frac{e^{-1/2+c^2/4}}{\sqrt{\pi c}n^{1/4}}.
\end{equation}
Hence using \eqref{as113}, we have
\begin{equation}
   \frac1{|\tilde{S}_n|}
\sum_{(**)}p^\symmS_{kml} |S_{k,m}|
\le Cn^{3/4}\cdot \max \bigl( e^{\sqrt{n}(a-1-a\log a)}, \ \
e^{\sqrt{n}(b-1-b\log b)} \bigr).
\end{equation}
But $f(x)=x-1-x\log x$ is increasing in $0<x<1$, is decreasing in $x>1$, and
$f(1)=0$.
Therefore there are positive constants $C$ and $c$ such that
for large $n$,
\begin{equation}\label{as119}
   \frac1{|\tilde{S}_n|}
\sum_{(**)}p^\symmS_{kml} |S_{k,m}|
\le C e^{-c\sqrt{n}}.
\end{equation}

On the other hand, Lemma \ref{lem19} says that (recall \eqref{as12-42})
for any fixed real number $d>0$, there is a constant $C$ such that
for $a\sqrt{n}\le m\le b\sqrt{n}$,
\begin{equation}
 \begin{split}
   &P^\symmS_l\biggl((2\mu^{(d)}_{k})^{1/2} ;
\mu^{(d)}_m(2\mu^{(d)}_{k})^{-1/2}\biggr)
-Cn^{-d/2}\\
&\qquad\qquad\qquad \le
p^\symmS_{kml}
\le
P^\symmS_l\biggl((2\nu^{(d)}_{k})^{1/2} ;
\nu^{(d)}_m(2\nu^{(d)}_{k})^{-1/2}\biggr)
+Cn^{-d/2},
 \end{split}
\end{equation}
for sufficiently large $n$.
Since $P^\symmS_l(t;\beta)\le P^\symmS_{l+1}(t;\beta)$,
Theorem \ref{cor12} for $P^\symmS_{2l+1}(t;\beta)$ yields
\begin{equation}
   e^{-\mu^{(d)}_{k}}D^{++}_{[(l-1)/2]}((2\mu^{(d)}_{k})^{1/2})
- Cn^{-d/2}
\le p^\symmS_{nml}
\le e^{-\nu^{(d)}_{k}}D^{++}_{[l/2]}((2\nu^{(d)}_{k})^{1/2})
+ Cn^{-d/2}.
\end{equation}

Let $l=[2\sqrt{n}+xn^{1/6}]$.
For $a\sqrt{n}\le m\le b\sqrt{n}$,
hence for $\frac{n-b\sqrt{n}}2 \le k\le \frac{n-a\sqrt{n}}2$,
\begin{equation}\label{as149}
  \bigl(l/2-(2\mu^{(d)}_{k})^{1/2}\bigr)2
(l/2)^{-1/3}, \ \
\bigl(l/2-(2\nu^{(d)}_{k})^{1/2}\bigr)2
(l/2)^{-1/3} =
x +O\bigl( n^{-1/6}\sqrt{\log n} \bigr).
\end{equation}
Also note that
from the asymptotics \eqref{as17}, \eqref{as18}  and
\eqref{as46}-\eqref{as49},
\begin{equation}\label{as150}
   (F(x)E(x))' = -\frac12 (v(x)+u(x))F(x)E(x)
\end{equation}
is bounded for $x\in\R$.
Hence using \eqref{as11-16} in Corollary \ref{cor7},
\eqref{as149} and \eqref{as150}, we obtain
\begin{equation}\label{as154}
 \begin{split}
   &\bigl|e^{-\nu^{(d)}_{k}}D^{++}_{[l/2]}((2\nu^{(d)}_{k})^{1/2})
-(FE)(x)\bigr|\\
&\qquad \le
\bigl|e^{-\nu^{(d)}_{k}}D^{++}_{[l/2]}((2\nu^{(d)}_{k})^{1/2})
-(FE)(\bigl(l/2-(2\nu^{(d)}_{k})^{1/2}\bigr)2
(l/2)^{-1/3})\bigr| \\
&\qquad \qquad+ \bigl|(FE)(\bigl(l/2-(2\nu^{(d)}_{k})^{1/2}\bigr)2
(l/2)^{-1/3}) - (FE)(x)\bigr|\\
&\qquad
\le C(M)n^{-1/6} + Ce^{-cM^{3/2}} + Cn^{-1/6}\sqrt{\log n}.
 \end{split}
\end{equation}
Therefore we have
\begin{equation}\label{as125}
   \sum_{(*)} p^\symmS_{nml} |S_{n,m}|
\le \biggl( F(x)E(x) +
C(M) n^{-1/6} + Ce^{-cM^{3/2}} + Cn^{-1/6}\sqrt{\log n} \biggr)
\sum_{(*)} |S_{n,m}|.
\end{equation}
Similarly,
\begin{equation}\label{as126}
   \sum_{(*)} p^\symmS_{nml} |S_{n,m}|
\ge \biggl( F(x)E(x) -
C(M) n^{-1/6} - Ce^{-cM^{3/2}} - Cn^{-1/6}\sqrt{\log n} \biggr)
\sum_{(*)} |S_{n,m}|.
\end{equation}
But from \eqref{as119},
\begin{equation}\label{as127}
   \frac1{|\tilde{S}_n|} \sum_{(*)} |S_{n,m}|
= 1 - \frac1{|\tilde{S}_n|} \sum_{(**)} |S_{n,m}|
= 1 + O(e^{-c\sqrt{n}}).
\end{equation}
Thus using \eqref{as114}, \eqref{as119}, \eqref{as125}, \eqref{as126}
and  \eqref{as127}, we obtain \eqref{as110}.
\end{proof}

\medskip
\begin{proof}[{\bf Proof of \eqref{as111}}]

As in Section 9 in \cite{BDJ}, integrating by parts,
\begin{equation}
  \Exp\bigl( (\tilde{\chi}_n)^p\bigr)
= \int_{-\infty}^\infty x^p dF_n(x)
= -\int_{-\infty}^0 px^{p-1}F_n(x) dx
+ \int_0^\infty px^{p-1}(1-F_n(x))dx,
\end{equation}
where $F_n(x):= \Prob(\tilde{\chi}_n \le x)
=\Prob(\tilde{L}_n\le 2\sqrt{n}+xn^{1/6})$.
From Theorem  \ref{cor12} and Corollary \ref{cor11}, we have
\begin{eqnarray}
\label{as163}  1- e^{-t^2/2} D^{++}_l(t) &\le& Ce^{-c|x|^{3/2}},
\qquad x\ge M, \\
\label{as164}   e^{-t^2/2} D^{++}_l(t)&\le& Ce^{-c|x|^{3}},
\qquad x\le -M,
\end{eqnarray}
for a fixed $M>M_2$ where $t=l-(x/2)l^{1/3}$.
Noting that $P_{2l+1}^\symmS(t;\beta)=e^{-t^2/2} D^{++}_l(t)$
for all $\beta\ge 0$, from
\eqref{as112}, Lemma \ref{lem20} and \eqref{as163}, \eqref{as164},
we obtain
\begin{eqnarray}
  1- F_n(x) &\le& Ce^{-c|x|^{3/2}},
\qquad x\ge M, \\
  F_n(x) &\le& Ce^{-c|x|^{3}},
\qquad x\le -M.
\end{eqnarray}
Now using convergence in distribution,
the dominated convergence theorem gives \eqref{as111}.
\end{proof}

\begin{rem}
We could also proceed using
\begin{equation}
   \Prob( \tilde{L}_n \le l)
= \frac1{|\tilde{S}_n|}
\sum_{2k+m=n} p^\symmO_{kml} |S_{k,m}|.
\end{equation}
The main contribution to the sum from $|S_{k,m}|$ comes from the region
$|m-\sqrt{n}|\le n^{1/4+\epsilon}$.  On the other hand, from Theorem
\ref{newthm2}, when $m= \sqrt{n} -2wn^{1/3}$, the quantity $p^\symmO_{kml}$
converges to $F(x;w)$.  Since the region $m= \sqrt{n}+ cn^{1/4+\epsilon}$
is much narrower than the region $m= \sqrt{n} +cn^{1/3}$, the main
contribution to the sum comes from when $w=0$, implying that
\begin{equation}
  \Prob( \tilde{L}_n \le l)
\sim \frac1{|\tilde{S}_n|}
\sum_{m=0}^{n} F(x;0) |S_{n,m}| = F_1(x).
\end{equation}
In the following proof for signed involutions, we make this
argument rigorous.
\end{rem}

\medskip
\begin{proof}[{\bf Proof of \eqref{as13-3}}]


   We have a disjoint union
\begin{equation}
   \tilde{S}^\symmu_n = \bigcup_{2k+m_++m_- =n}
S^\symmu_{k,m_+,m_-}.
\end{equation}
Hence again
\begin{equation}
   \Prob(\tilde{L}^\symmu_n \le l) =
\frac1{\bigl| \tilde{S}^\symmu_n \bigr|}
\sum_{2k+ m_++m_- =n} p^\symmu_{km_+m_-l}
\bigl| S^\symmu_{k,m_+,m_-} \bigr|.
\end{equation}
One can check that
\begin{equation}
   \bigl| S^\symmu_{k,m_+,m_-} \bigr| =
\frac{(2k+m_++m_-)!}{k!m_+!m_-!}.
\end{equation}
Hence, we have
\begin{equation}\label{as173}
   \bigl| \tilde{S}^\symmu_n \bigr| =
\sum_{2k+ m_++m_- =n} \bigl| S^\symmu_{k,m_+,m_-} \bigr|
= \sum_{0\le k\le [\frac{n}2]}
\sum_{0\le m_+\le n-2k} f(m_+,k)
\end{equation}
where
\begin{equation}
   f(m_+,k) := \frac{n!}{m_+!(n-m_+-2k)!k!}.
\end{equation}
For fixed $0\le k\le [\frac{n}2]$, $f(m_+,k)$ is unimodal in $m_+$
and achieves its maximum when $m_+\sim n/2-k$.
And $f(n/2-k,k)$ is unimodal in $k$ and the maximum is attained
when $k\sim n/2-\sqrt{n/2}$.
Hence $f(m_+,k)$ has its maximum when
$(m_+,k)\sim (\sqrt{\frac{n}2}, \frac{n}2-\sqrt{\frac{n}2})$.
Consider the disc $D$ of radius $n^{1/4+\epsilon}$ centered at
$(\sqrt{\frac{n}2}, \frac{n}2-\sqrt{\frac{n}2})$.
We will show that the main contribution to the sum in
\eqref{as173} comes from $D$.
Set
\begin{equation}
  m_+=\sqrt{\frac{n}2}+x, \quad k= \frac{n}2-\sqrt{\frac{n}2}+y,
\qquad |x|,|y|\le n^{1/4+\epsilon},
\end{equation}
By Stirling's formula,
\begin{equation}\label{as176}
  f(m_+,k) = \frac1{\sqrt{en}\pi} (2n)^{n/2}
e^{-n/2+\sqrt{2n}}e^{-\frac{x^2+(x+2y)^2}{\sqrt{2n}}}
\bigl( 1+O(n^{-1/4+3\epsilon/2}) \bigr).
\end{equation}
Hence from the unimodality discussed above,
\begin{equation}\label{as177}
   \sum_{(m_+,k)\notin D} f(m_+,k) \le n^2
\max_{(m_+,k)\in\partial D} f(m_+,k)
\le \frac{n^2}{\sqrt{en}\pi} (2n)^{n/2}
e^{-n/2+\sqrt{2n}} e^{-5\sqrt{2}n^{2\epsilon}},
\end{equation}
and by summing up using \eqref{as176},
\begin{equation}
  \sum_{(m_+,k)\in D} f(m_+,k)
= \frac1{\sqrt{2e}} (2n)^{n/2}e^{-n/2+\sqrt{2n}}
\bigl( 1+O(n^{-1/4+3\epsilon/2}) \bigr).
\end{equation}
Hence we have
\begin{equation}\label{as179}
   \bigl| \tilde{S}^\symmu_n \bigr| =
\frac1{\sqrt{2e}} (2n)^{n/2}e^{-n/2+\sqrt{2n}}
\bigl( 1+O(n^{-1/4+3\epsilon/2}) \bigr),
\end{equation}
and the main contribution to the sum in \eqref{as173}
comes from $D$.

As in Theorem \ref{thm22}, we write
\begin{equation}\label{as201}
   \Prob(\tilde{L}^\symmu_n\le l)
= \frac1{\bigl| \tilde{S}^\symmu_n \bigr|}
\biggl[ \sum_{D} p_{km_+m_-l}
\bigl| \tilde{S}^\symmu_{k,m_+,m_-} \bigr|
+ \sum_{D^c} p_{km_+m_-l}
\bigl| \tilde{S}^\symmu_{k,m_+,m_-} \bigr| \biggr].
\end{equation}
From \eqref{as177} and \eqref{as179},
\begin{equation}\label{as202}
   \frac1{\bigl| \tilde{S}^\symmu_n \bigr|}
\sum_{D^c} p_{nm_+m_-l}
\bigl| \tilde{S}^\symmu_{n,m_+,m_-} \bigr|
\le e^{-10n^{2\epsilon}}.
\end{equation}

On the other hand, by the remark to Lemma \ref{lem19}
and Theorem \ref{cor12}, (recall \eqref{as12-43})
\begin{equation}
   p_{km_+m_-l}
\le e^{-\nu^{(d)}_{m_+}-\nu^{(d)}_{k}}\pi^*_{[l/2]}
\biggl(-\frac{\nu^{(d)}_{m_+}}{(\nu^{(d)}_{k})^{1/2}};
(\nu^{(d)}_{k})^{1/2}\biggr)
D_{[\frac{l}2]} ((\nu^{(d)}_{k})^{1/2}) + Cn^{-d/2},
\end{equation}
for large $n$.
We have a similar inequality of the other direction with
$\nu$, $l$ and $+Cn^{-d/2}$ replaced by $\mu$, $l-1$
and $-Cn^{-d/2}$.

Let $l=[2\sqrt{2n}+x2^{2/3}(2n)^{1/6}]$.
In the region $D$,
\begin{equation}
   (l/2-4(\nu^{(d)}_{k})^{1/2})(l/4)^{-1/3} = x+O(n^{-1/6}\sqrt{\log n}),
\end{equation}
and
\begin{equation}
   (1-\nu^{(d)}_{m_+}/(\nu^{(d)}_{k})^{1/2})2^{-4/3}(l/2)^{1/3}
= O(n^{-1/12+\epsilon}).
\end{equation}
Hence as in \eqref{as154}, using Corollary \ref{cor7} \eqref{as11-14},
Proposition \ref{thm3} and Corollary \ref{cor5}
\begin{equation}
 \begin{split}
   &\biggl| e^{-\nu^{(d)}_{m_+}-\nu^{(d)}_{k}}\pi^*_{[l/2]}
\biggl(-\frac{\nu^{(d)}_{m_+}}{(\nu^{(d)}_{k})^{1/2}};
(\nu^{(d)}_{k})^{1/2}\biggr)
D_{[\frac{l}2]} ((\nu^{(d)}_{k})^{1/2})
- F^\symmu(x;0) \biggr| \\
&\qquad\qquad \le
Cn^{-1/6}\sqrt{\log n} + Cn^{-1/12+\epsilon}
+ C(M)n^{-1/6} + Ce^{-cM^{3/2}}
 \end{split}
\end{equation}
for a constant $C(M)$ which may depend on $M$ and
constants $C$ and $c$ which is independent of $M$.
Thus for large $n$,
\begin{equation}
   \Prob\biggl(
\frac{\tilde{L}^\symmu_n - 2\sqrt{2n}}{2^{2/3}(2n)^{1/6}}\le x\biggr)
\le F^\symmu(x;0) + e(n,M)
\end{equation}
with some error $e(n,M)$ such that
$\lim_{M\to\infty} \lim_{n\to\infty} e(n,M) =0$.
Similarly we have an inequality for the other direction.
Recalling $F^\symmu(x;0)=F_1(x)^2$ from \eqref{as116},
we obtain \eqref{as13-3}.
\end{proof}

\medskip
\begin{proof}[{\bf Proof of \eqref{as13-4}}]

  Integrating by parts,
\begin{equation}
  \Exp\bigl( (\tilde{\chi}^\symmu_n)^p\bigr)
= \int_{-\infty}^\infty x^p dF_n(x)
= -\int_{-\infty}^0 px^{p-1}F_n(x) dx
+ \int_0^\infty px^{p-1}(1-F_n(x))dx,
\end{equation}
where $F_n(x):= \Prob(\tilde{\chi}^\symmu_n \le x)
= \Prob(\tilde{L}^\symmu_n \le 2\sqrt{2n}+x2^{2/3}(2n)^{1/6})$.
Note that when $x< -(4n)^{1/3}$, $F_n(x)=0$,
and when $x> 2^{1/6}n^{5/6}-(4n)^{1/3}$, $F_n(x)=1$.

Let $M>M_0$ fixed.
Consider the case when $-(4n)^{1/3} \le x\le -M$.
 From \eqref{as201} and \eqref{as202},
\begin{equation}
  F_n(x) \le
\frac1{\bigl| \tilde{S}^\symmu_n \bigr|}
 \sum_{D} p_{km_+m_-l}
\bigl| \tilde{S}^\symmu_{k,m_+,m_-} \bigr|
+ Ce^{-10n^{2\epsilon}},
\end{equation}
where $l=[2\sqrt{2n}+x2^{2/3}(2n)^{1/6}]$.
We apply Lemma \ref{lem20}, Corollary \ref{cor11} and Corollary \ref{cor-new}
\eqref{as7.49}.
Note that we are in the region $\alpha \to 1$ faster than $k^{-1/3}$,
hence $w$ is bounded, say $-1\le w\le 1$.
So we can apply Corollary \ref{cor-new} \eqref{as7.49}.
Then we obtain
\begin{equation}
  F_n(x) \le Ce^{-c|x|^3} + Ce^{-10n^{2\epsilon}}.
\end{equation}
Since $-(4n)^{1/3} \le x\le -M$, we have
\begin{equation}
  e^{-10n^{2\epsilon}} \le e^{-\frac{10}{2^{4\epsilon}}|x|^{6\epsilon}},
\end{equation}
thus,
\begin{equation}
  F_n(x) \le Ce^{-c|x|^3} + Ce^{-\frac{10}{2^{4\epsilon}}|x|^{6\epsilon}}.
\end{equation}

On the other hand, when
$M\le x\le 2^{1/6}n^{5/6}-(4n)^{1/3}$, similarly we have
\begin{equation}
  1-F_n(x) \le
\frac1{\bigl| \tilde{S}^\symmu_n \bigr|}
 \sum_{D} (1-p_{km_+m_-l})
\bigl| \tilde{S}^\symmu_{k,m_+,m_-} \bigr|
+ Ce^{-10n^{2\epsilon}},
\end{equation}
Using Lemma \ref{lem20}, Corollary \ref{cor11} and Proposition \ref{thm4}
\eqref{as7.45}, we obtain
\begin{equation}
  1-F_n(x) \le Ce^{-c|x|^{3/2}} + Ce^{-10n^{2\epsilon}}.
\end{equation}
Since $M\le x\le 2^{1/6}n^{5/6}-(4n)^{1/3}$,
\begin{equation}
   e^{-10n^{2\epsilon}} \le
e^{-\frac{10}{2^{2\epsilon/5}}|x|^{12\epsilon/5}},
\end{equation}
thus
\begin{equation}
  1-F_n(x) \le Ce^{-c|x|^{3/2}}
+ Ce^{-\frac{10}{2^{2\epsilon/5}}|x|^{12\epsilon/5}}.
\end{equation}
Therefore, using dominated convergence theorem, we obtain \eqref{as13-4}.
\end{proof}

%
%

\section{Asymptotics for $\alpha>1$}\label{sec:bigalpha}

As we remarked after Theorem \ref{newthm3}, when $\alpha>1$, we must use a
different scaling to obtain useful results.

Let $L^\symmO(t;\alpha)$ and $L^\symmu(t;\alpha,\beta)$ be random variables
with the distribution functions given by $\Prob(L^\symmO(t;\alpha)\le l)=
P^\symmO_l(t;\alpha)$ and $\Prob(L^\symmu(t;\alpha,\beta)\le l)=
P^\symmu_l(t;\alpha,\beta)$, respectively: the Poissonized version of
$L^\symmO$ and $L^\symmu$. Under appropriate scalings, we obtain the Gaussian
distribution in the limit.

\begin{thm}
For $\alpha>1$ and $\beta\ge 0$ fixed,
\begin{eqnarray}
\label{b1}
  \lim_{t\to\infty} \Prob\biggl(
\frac{L^\symmO(t;\alpha)-(\alpha+\alpha^{-1})t}
{\sqrt{(\alpha-\alpha^{-1})t}}
\le x \biggr)
&=& \frac1{\sqrt{2\pi}} \int_{-\infty}^{x} e^{-\frac12 y^2} dy,\\
\label{b2}
  \lim_{t\to\infty} \Prob\biggl(
\frac{L^\symmu(t;\alpha,\beta)-2(\alpha+\alpha^{-1})t}
{\sqrt{2(\alpha-\alpha^{-1})t}}
\le x \biggr)
&=& \frac1{\sqrt{2\pi}} \int_{-\infty}^{x} e^{-\frac12 y^2} dy.
\end{eqnarray}
\end{thm}

\begin{proof}
  Let $l=(\alpha+\alpha^{-1})t+\sqrt{(\alpha-\alpha^{-1})t}$
for $L^\symmO$.
For large $t$, $2t/l\le c<1$ for some $c>0$.
  From Proposition \ref{thm1} (i), using Theorem \ref{cor12},
it is easy to see that
$e^{-t^2/2}D_l^{\pm\pm}(t), e^{-t^2\pm t}D_l^{\pm\mp}(t) \to 1$
exponentially as $l\to\infty$.
Now Theorem \ref{cor12}, \eqref{as7.42}
and Proposition \ref{prop7.7} imply \eqref{b1}.
 For $L^\symmu$, let $l=2(\alpha+\alpha^{-1})t+\sqrt{2(\alpha-\alpha^{-1})t}$.
Similarly,
$e^{-t^2}D_l(t) \to 1$ exponentially as $l\to\infty$, and
we obtain \eqref{b2}.
\end{proof}

Unfortunately, we can no longer apply the de-Poissonization technique; the
difficulty is that $(\alpha+\alpha^{-1})t$ depends too strongly on small
perturbations in $\alpha$.  Indeed, as we shall see, the asymptotics
of the non-Poisson processes are different.

Consider the case of involutions with $[\alpha t]$ fixed points and
$[t^2/2]$ 2-cycles; the case of signed involutions is analogous.  By
symmetry, this is the same as the largest increasing subset distribution
for the point selection process in the triangle $0\le y\le x\le 1$ with
$[t^2/2]$ generic points and $[\alpha t]$ diagonal points.  As was observed
in Remark 2 to Corollary 7.6 of \cite{PartI}, it is equivalent to consider
weakly increasing subsets where the extra points are added to the line
$y=0$ instead of to the diagonal.

As in \eqref{res1}, let
\[
\chi^\symmO_{[t^2/2],[\alpha t]} =
\frac{L^\symmO_{[t^2/2],[\alpha t]}
-(\alpha+1/\alpha)t}{\sqrt{(1/\alpha-1/\alpha^3)t}}.
\]

\begin{thm}\label{thmbigalpha1}
As $t\to\infty$, the variable $\chi^\symmO_{[t^2/2],[\alpha t]}$
converges in distribution and
moments to $N(0,1)$.
\end{thm}

\begin{proof}
Let $S(t)$ be the set of points at time $t$, and let $I$ be a largest
increasing subset of $S(t)$.  Then there will exist some number $0\le s^+\le
1$ (not unique) such that
\[
(S(t)\cap \{y=0,0\le x\le s^+\}) \subset I,
\]
and such that every other point of $I$ has $x>s^+$ and $y>0$.  For any
$0\le s\le 1$, we thus have
\[
f_1(s)+f_2(s) \le |I|,
\]
where $f_1(s)$ is the number of points of $S(t)$ with $y=0$ and $0\le x\le
s$, and where $f_2(s)$ is the largest increasing subset of $S(t)$ lying
entirely in the (part-open) trapezoid with $x\ge s$, $y>0$.

Since $f_1(s)$ is binomial with parameters $[\alpha t]$ and $s$,

\begin{lem}
Let $M>0$ be sufficiently large and fixed.  For all $0\le s<1$, there exist
positive constants $C$, $c$ independent of $s$ such that for $w\ge M$,
\[
\Prob(f_1(s)>\alpha s t+w t^{1/2}) \le C e^{-c|w|^2},
\]
while for $w\le -M$,
\[
\Prob(f_1(s)<\alpha s t+w t^{1/2}) \le C e^{-c|w|^2}.
\]
\end{lem}

For $f_2$, we have:

\begin{lem}
Let $M>0$ be sufficiently large and fixed.  For $0\le s<1$, there exist
positive constants $C$ and $c$ independent of $s$ such that for all $w\ge M$,
\[
\Prob(f_2(s) > 2\sqrt{(1-s)}t + w t^{1/3}) \le C e^{-c|w|^{3/2}},
\]
and for all $w\le -M$,
\[
\Prob(f_2(s) < 2\sqrt{(1-s)}t + w t^{1/3}) \le C e^{-c|w|^3}.
\]
\end{lem}

\begin{proof}
We first show the corresponding large-deviation result for the
Poissonization.  Define $f'_2(s,t)$ to be the length of the longest
increasing subsequence when the number of points in the trapezoid is
Poisson with parameter $t^2(1-s^2)/2$.  Then $f'_2(s,t)$ is bounded between
the corresponding processes for the rectangle $s\le x\le 1$, $0\le y\le 1$
and for the triangle $0\le (x-s)/(1-s)\le y\le 1$.  In particular, if
$f'_2(s,t)$ deviates significantly from $\sqrt{1-s}t$, so must the
appropriate bounding process; the result follows immediately from the
corresponding results for rectangles and triangles.

The corresponding large-deviation result when the number of points is fixed
then follows from Lemma \ref{lem11}.  In our case, the number of points in
the trapezoid is binomial with parameters $t^2/2$ and $(1-s^2)$; the lemma
follows via essentially the same argument used to prove Lemma \ref{lem11}.
\end{proof}

As we will see, the value $s=1-\alpha^{-2}$ deserves special attention:

\begin{lem}
Let $M>0$ be sufficiently large and fixed.  There exist positive constants
$C$, $c$ such that for $w\ge M$,
\[
\Prob(f_1(1-\alpha^{-2})+f_2(1-\alpha^{-2}) >
(\alpha+1/\alpha) t + w t^{1/2}) \le C e^{-c\min(|w|^2,|w|^{3/2}t^{1/4})},
\]
and for $w\le -M$,
\[
\Prob(f_1(1-\alpha^{-2})+f_2(1-\alpha^{-2}) <
(\alpha+1/\alpha) t + w t^{1/2}) \le C e^{-c|w|^2}.
\]
Moreover, if we define
\[
\chi_0(t)=
\frac{f_1(1-\alpha^{-2})+f_2(1-\alpha^{-2}) - (\alpha+1/\alpha)t}
{\sqrt{(1/\alpha-1/\alpha^3)t}},
\]
then $\chi_0(t)$ converges to a standard normal distribution, both
in distribution and in moments.
\end{lem}

\begin{proof}
That $\chi_0(t)$ converges as stated follows from the fact that if we
write $\chi_0(t)=\chi_1(t)+\chi_2(t)$, with
\begin{align}
\chi_1(t) &= \frac{f_1(1-\alpha^{-2})-(\alpha-1/\alpha)t}{\sqrt{(1/\alpha-1/\alpha^3)t}}\\
\chi_2(t) &=
\frac{f_2(1-\alpha^{-2})-(2/\alpha)t}{\sqrt{(1/\alpha-1/\alpha^3)t}},
\end{align}
then $\chi_1(t)$ converges in distribution and moments to a standard normal
distribution, and $\chi_2(t)$ converges in distribution and moments to 0.

For the large deviation bounds, we note that if $x+y>z+w$, then either
$x>z$ or $y>w$.  Thus for any $0\le b\le 1$, we have
\begin{align}
\Prob(f_1(1-\alpha^{-2}&)+f_2(1-\alpha^{-2}) >
(\alpha+1/\alpha) t + w t^{1/2})\\
&\le
\Prob(f_1(1-\alpha^{-2})>(\alpha-1/\alpha) t + b w t^{1/2})
+
\Prob(f_2(1-\alpha^{-2})>(2/\alpha) t + (1-b) w t^{1/2})\\
&
\le
C e^{-c |bw|^2}+C e^{-c |(1-b)w|^{3/2} t^{1/4}};
\end{align}
the result follows by balancing the two terms.  In the other case, the $C
e^{-c|w|^2}$ term always dominates.
\end{proof}

\begin{lem}
For any sufficiently small $\epsilon>0$, there exist positive constants
$C$, $c$ such that
\[
\Pr(s^+-(1-1/\alpha^2)>t^{\epsilon/3-1/3})
<
C e^{-c t^\epsilon}
\]
and
\[
\Pr((1-1/\alpha^2)-s^+>t^{\epsilon/2-1/2})
<
C e^{-c t^\epsilon}
\]
for all sufficiently large $t$.
\end{lem}

\begin{proof}
Define a sequence $s_i$ by taking
\[
s_i  = 1-(1-2/(i+2))^2/\alpha^2
\]
for all $i\ge 0$.  Similarly define a sequence $s'_i$ by
\[
s'_i = \max(t_i,0),
\]
with
\[
t_i=1-(1+2 e^{2^{-1-i}})^2/\alpha^2
\]
for $i<0$ and
\[
t_i=1-(1+4/(i+1))^2/\alpha^2
\]
for $i\ge 0$.  Note that $s_i$ is strictly decreasing and $t_i$ is
strictly increasing.

\begin{lem}
For all $i\ge 0$,
\[
\alpha s_i + 2 \sqrt{1-s_{i+1}} < \alpha+1/\alpha.
\]
For all $i$,
\[
\alpha s'_{i+1} + 2\sqrt{1-s'_{i}} < \alpha+1/\alpha.
\]
\end{lem}

\begin{proof}
In the first case, we have
\[
\alpha+1/\alpha - (\alpha s_i + 2 \sqrt{1-s_{i+1}})
=
\frac{4}{(i+2)^2(i+3)\alpha}.
\]
In the second case, it suffices to verify the formula with $s'$ replaced
by $t$.  For $i<-1$,
\[
\alpha+1/\alpha - (\alpha t_{i+1} +2\sqrt{1-t_{i}})
=
4 e^{2^{-i}/4}/\alpha.
\]
For $i=-1$,
\[
\alpha+1/\alpha - (\alpha t_{i+1} +2\sqrt{1-t_{i}})
=
(24-4e)/\alpha
\]
Finally, for $i\ge 0$,
\[
\alpha+1/\alpha - (\alpha t_{i+1} +2\sqrt{1-t_{i}})
=
\frac{8i}{(i+1) (i+2)^2 \alpha}.
\]
\end{proof}

Let $i_1=t^{1/6-\epsilon/6}$, $i_2=t^{1/4-\epsilon/4}$.  Then
there exist constants $C$, $c$ such that for $0\le i\le i_1$,
\[
\Prob( f_1(s_i)+f_2(s_{i+1}) > f_1(1-\alpha^{-2})+f_2(1-\alpha^{-2}))
<
C e^{-c t^{\epsilon}}.
\]
Since $f_1(s_i)+f_2(s_{i+1})$ is an upper bound on $f_1(s)+f_2(s)$
with $s_{i+1}\le s \le s_i$, it follows that
\[
\Prob( s^+\in [s_{i+1},s_i] )<C e^{-c t^{\epsilon}}.
\]
Similarly, for $i\le i_2$,
\[
\Prob( f_1(s'_{i+1})+f_2(s'_i) > f_1(1-\alpha^{-2})+f_2(1-\alpha^{-2}))
<
C e^{-c t^{\epsilon}},
\]
and thus
\[
\Prob( s^+\in [s'_i,s'_{i+1}] )<C e^{-c t^{\epsilon}}.
\]
Since there are only $i_1+i_2+\log\log\alpha$ such events to
consider, the result follows.
\end{proof}

In particular, with probability
$1-C e^{-c t^{\epsilon}}$, we have
\[
f_1(1-\alpha^{-2}-t^{\epsilon/2-1/2})
+
f_2(1-\alpha^{-2}+t^{\epsilon/3-1/3})
\le
L(t)
\le
f_1(1-\alpha^{-2}+t^{\epsilon/3-1/3})
+
f_2(1-\alpha^{-2}-t^{\epsilon/2-1/2}).
\]

But then using the fact that
\[
f_1(1-\alpha^{-2})-f_1(1-\alpha^{-2}-t^{\epsilon/2-1/2})
\]
and
\[
f_1(1-\alpha^{-2}-t^{\epsilon/2-1/2})-f_1(1-\alpha^{-2})
\]
are Poisson, and using the large deviation behavior of $f_2(s)$,
we find that
\[
\Prob(
L(t)-(
f_1(1-\alpha^{-2}-t^{\epsilon/2-1/2})
+
f_2(1-\alpha^{-2}+t^{\epsilon/3-1/3})
)
\ge t^{1/2-\epsilon'}
)
\le C e^{-c t^\epsilon},
\]
and
\[
\Prob(
(
f_1(1-\alpha^{-2}-t^{\epsilon/3-1/3})
+
f_2(1-\alpha^{-2}+t^{\epsilon/2-1/2})
)
-
L(t)
\ge t^{1/2-\epsilon'}
)
\le C e^{-c t^\epsilon}.
\]
So $\chi(t)-\chi_0(t)$ converges to 0 in a fairly strong sense;
in particular, they must have the same limiting distribution and
limiting moments.
\end{proof}

\begin{rem}
The above proof could equally well be applied to the Poisson process;
the beta distribution would then be replaced by a Poisson distribution.
\end{rem}

For the signed involution case, with $[2\alpha t]$ fixed points,
$[2\beta t]$ negated points, and $[2t^2]$ 2-cycles, again we let
\[
\chi^\symmu_{[t^2],[2\alpha t],[2\beta t]} =
\frac{L^\symmu_{[t^2/2],[\alpha t],[\beta t]}
-2(\alpha+1/\alpha)t}{\sqrt{2(1/\alpha-1/\alpha^3)t}}.
\]
Then the analogous argument proves that $\chi'(t)$ also converges in
moments and distribution to a standard normal.

%
%

\section{Steepest descent type analysis for Riemann-Hilbert problems}
\label{RHP}

In this section, we prove the asymptotics of orthogonal polynomials results
stated in Section \ref{as-op} by applying the steepest descent type method
to the Riemann-Hilbert problem (RHP) \eqref{as-RHP}. The steepest descent
method for RHP's, the Deift-Zhou method, was introduced by Deift and Zhou
in \cite{DZ1} and is developed further in \cite{DZ2}, \cite{DVZ1}, and
finally placed in a systematic form by Deift, Venakides and Zhou in
\cite{DVZ}. The steepest descent analysis of the RHP \eqref{as-RHP} was
first conducted in \cite{BDJ}.  The analysis of \cite{BDJ} has many
similarities with \cite{DKMVZ, DKMVZ2, DKMVZ3} where the asymptotics of
orthogonal polynomials on the real line with respect to a general weight is
obtained, leading to a proof of universality conjectures in random matrix
theory.

As mentioned in the Introduction and Section \ref{as-op}, in this section we
extend the analysis of \cite{BDJ} and obtain new estimates on the orthogonal
polynomials $\pi_k(z;t)$. The extension is done roughly in two categories. In
\cite{BDJ}, the quantity of interest was $Y_{21}(0;k;t)$, and so the
$z$-dependence of the error bound of $Y(z;k;t)$ was not considered carefully.
But in the present paper we need the asymptotics of $Y(z)$ for general $z\in
\C$, and also for the case when $z\to -1$ as $k,t\to\infty$. Hence the first
category of our extension is to investigate how the error estimate
depends on $z$. This task sometimes requires improved estimates of the
solution $Y(z)$.  See for example, \eqref{qz10.41} below where an improved
$L^1$ norm bound of the jump matrix is needed. On the other hand, as we
will see, the asymptotic solution $Y(z)$ is expressed in terms of the
so-called the $g$-function. Thus we need detailed analysis of the
$g$-function to obtain the asymptotics of the orthogonal polynomials. In
the special case $z=0$, we have $g(0)=\pi i$ (Lemma 4.2 in
\cite{BDJ}). Hence in \cite{BDJ}, the analysis of the $g$-function was
quite simple. But in the present paper, we need general values of $g(z)$,
and this in some cases requires further analysis. Hence the analysis of the
$g$-function is the second category of our extension. For example, see
\eqref{qz10.109}-\eqref{qz10.121} below where we need a further analysis of the
$g$-function.

Again the analysis in this section relies heavily on the analysis of
\cite{BDJ} and we extend the method of \cite{BDJ}. Thus the analysis below
overlaps in many parts with the analysis of \cite{BDJ}. However, for
continuity of presentation and also for the convenience of readers, we
include some calculations which overlap \cite{BDJ}. When the analysis
overlaps that of \cite{BDJ}, we only sketch the method, and instead we
focus on new features to indicate how to prove the propositions in Section
\ref{as-op}.

\medskip
We say that a RHP is normalized at $\infty$ if the solution $m$
satisfies the condition $m\to I$ as $z\to\infty$.
Thus for instance, the RHP's \eqref{as20}, \eqref{as12.1} are normalized
at $\infty$, while the RHP \eqref{as-RHP} is not.

In \cite{BDJ}, it turned out that the asymptotic analysis differs critically
when $(2t)/k\le 1$ and $(2t)/k>1$, due to the difference of (the support of)
the associated equilibrium measure (see Lemma 4.3 \cite{BDJ}). Hence we discuss
those two cases separately in Subsections \ref{sub1} and \ref{sub2}, which
extend Section 5 and Section 6 of \cite{BDJ}, respectively. Each section is
also divided into three subcases. In each subcase the corresponding case of the
propositions in Section \ref{as-op} (except Proposition \ref{prop7.7}) are
proved. Subsection \ref{sub3} is new, and Proposition \ref{prop7.7} is proved
in this subsection.

%
%

\subsection{When $(2t)/k \le 1$.}\label{sub1}

The following algebraic transformations \eqref{as12.1}--\eqref{as12.11} of
RHP's are taken from (5.1)--(5.3) of \cite{BDJ}.

Define
\begin{equation}\label{as12.1}
\begin{split}
  m^{(1)}(z;k;t) := Y(z;k;t)
\begin{pmatrix} (-1)^ke^{tz}&0 \\ 0&(-1)^ke^{-tz} \end{pmatrix},
\qquad &|z|<1, \\
  m^{(1)}(z;k;t) := Y(z;k;t)
\begin{pmatrix} z^{-k}e^{tz^{-1}}&0 \\ 0&z^ke^{-tz^{-1}} \end{pmatrix},
\qquad &|z|>1.
\end{split}
\end{equation}
Then $m^{(1)}$ solves a new RHP which is equivalent to the RHP
\eqref{as-RHP} in the sense that a solution of one RHP yields algebraically
a solution of the other RHP, and vice versa:
\begin{equation}
\begin{cases}
m^{(1)}(z;k;t) \qquad \text{is analytic in $\C\setminus\Sigma$,}\\
m_+^{(1)}(z;k;t)= m_-^{(1)}(z;k;t) \begin{pmatrix}
(-1)^kz^ke^{t(z-z^{-1})}& (-1)^k\\
0& (-1)^kz^{-k}e^{-t(z-z^{-1})}
\end{pmatrix} \quad\text{on $\Sigma$},\\
m^{(1)}(z;k;t)=I+O(\frac1z) \qquad \text{as $z\rightarrow \infty$},\\
\end{cases}
\end{equation}
where $\Sigma$ is the unit circle oriented counterclockwise as before.
Here and in the sequel, $m_+(z)$ (resp. $m_-(z)$) is understood as the limit
from the left (reps. right) side of the contour as one goes along
the orientation of the contour.
Now we define $m^{(2)}(z;k;t)$ in terms of $m^{(1)}(z;k;t)$ as follows:
\begin{equation}\label{as12.4}
 \begin{split}
  &\text{for even $k$,} \\
  &\quad \begin{cases}
    m^{(2)} \equiv m^{(1)}& \quad |z|>1,\\
    m^{(2)} \equiv m^{(1)} \bigl( \begin{smallmatrix} 0& -1\\ 1&0
\end{smallmatrix} \bigr)& \quad |z|<1. \end{cases}\\
  &\text{for odd $k$,}\\
  &\quad\begin{cases}
    m^{(2)} \equiv \bigl( \begin{smallmatrix} 1&0\\0&-1
\end{smallmatrix} \bigr)  m^{(1)} \bigl(
\begin{smallmatrix} 1&0\\0&-1 \end{smallmatrix} \bigr)& \quad |z|>1,\\
    m^{(2)} \equiv \bigl( \begin{smallmatrix} 1&0\\0&-1
\end{smallmatrix} \bigr) m^{(1)} \bigl(
\begin{smallmatrix} 0&-1\\-1&0 \end{smallmatrix} \bigr)&
\quad |z|<1. \end{cases}
 \end{split}
\end{equation}
Then $m^{(2)}(\cdot;k;t)$ solves another RHP
\begin{equation}
\begin{cases}
m^{(2)}_+(z;k;t)=m^{(2)}_-(z;k;t) v^{(2)}(z;k;t)\quad\text{on} \quad\Sigma,\\
m^{(2)}(z;k;t)=I+O(\frac1z) \ \ \text{as}\ \ z\rightarrow \infty\\
\end{cases}
\end{equation}
where
\begin{equation}
   v^{(2)}(z;k;t)= \begin{pmatrix}
1& -(-1)^kz^ke^{t(z-z^{-1})} \\
(-1)^kz^{-k}e^{-t(z-z^{-1})}&0 \end{pmatrix}.
\end{equation}
The jump matrix has the following factorization,
\begin{equation}\label{as12.7}
    v^{(2)}= \begin{pmatrix} 1&0\\
(-1)^kz^{-k}e^{-t(z-z^{-1})}&1\end{pmatrix}
\begin{pmatrix} 1&-(-1)^kz^ke^{t(z-z^{-1})}\\0&1\end{pmatrix}
=: (b_-^{(2)})^{-1}b_+^{(2)},
\end{equation}
We note that through the changes $Y\to m^{(1)} \to m^{(2)}$,
we have
\begin{eqnarray}
\label{as12.8}   Y_{11}(z;k;t) &=& -(-1)^ke^{-tz}m^{(2)}_{12}(z;k;t),
\qquad |z|<1,\\
\label{as12.9}   Y_{21}(z;k;t) &=& -e^{-tz}m^{(2)}_{22}(z;k;t),
\qquad\qquad |z|<1,\\
\label{as12.10}   Y_{11}(z;k;t) &=& z^ke^{-tz^{-1}}m^{(2)}_{11}(z;k;t),
\quad \qquad |z|>1,\\
\label{as12.11}   Y_{21}(z;k;t) &=& (-z)^ke^{-tz^{-1}}m^{(2)}_{21}(z;k;t),
\qquad |z|>1.
\end{eqnarray}

As in (5.4) of \cite{BDJ}, the absolute value of the (12)-entry of the jump
matrix $v^{(2)}$ is $e^{kF(\rho,\theta;\frac{2t}{k})}$ where
\begin{equation}\label{as13.11}
  F(z;\gamma)=F(\rho e^{i\theta};\gamma)
:= \frac{\gamma}{2}(\rho-\rho^{-1})\cos\theta +\log\rho,
\quad z=\rho e^{i\theta}.
\end{equation}
The absolute value of the (21)-entry of $v^{(2)}$ is
$e^{-kF(\rho e^{i\theta};\frac{2t}{k})}$.
Note that
\begin{equation}\label{as12.12}
   F(\rho,\theta;\gamma)=-F(\rho^{-1},\theta;\gamma).
\end{equation}
\begin{figure}[ht]
 \centerline{\epsfig{file=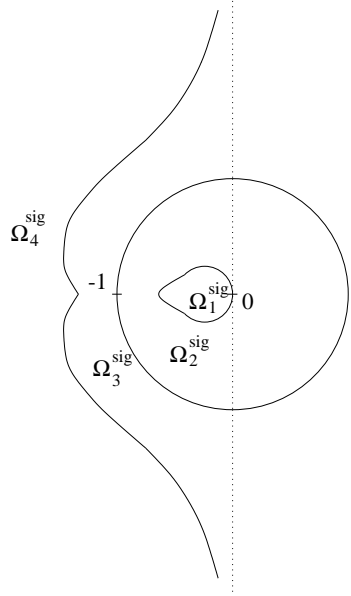, width=4.5cm}}
 \caption{curves of $F(z;\gamma)=0$ when $0<\gamma<1$}
\label{fig-sig}
\end{figure}
Figure \ref{fig-sig} shows the curves $F(z;\gamma)=0$.
In $\Omega^{sig}_1\cup\Omega^{sig}_3$, $F>0$, and
in $\Omega^{sig}_2\cup\Omega^{sig}_4$, $F<0$.
The region $\Omega^{sig}_2$ becomes smaller as $\gamma$ increases,
and when $\gamma=1$, the curve $F(z;\gamma)=0$ contacts the unit circle
$\Sigma$ at $z=-1$ with the angle $\pi/3$.

We distinguish three subcases as in Section 5 of \cite{BDJ}.

\subsubsection{The case $0\le 2t\le ak$ for some $0<a<1$.}\label{subsub1}

It is possible to fix $\rho_a<1$ such that the circle $\{ z: |z|=\rho_a\}$ is
in the region $\Omega^{sig}_2$ for all such $t$ and $k$. Define
$m^{(3)}(z;k;t)$ by (see (5.9) of \cite{BDJ})
\begin{equation}\label{as12.13}
  \begin{cases}
     m^{(3)}=m^{(2)} (b_+^{(2)})^{-1},
\qquad &\rho_a<|z|<1 , \\
     m^{(3)}=m^{(2)} (b_-^{(2)})^{-1},
\qquad &1<|z|<\rho_a^{-1}, \\
     m^{(3)}=m^{(2)}
\qquad\qquad &|z|<\rho_a,\ \ |z|>\rho_a^{-1}.
  \end{cases}
\end{equation}
Then $m^{(3)}$ satisfies a new jump condition $m^{(3)}_+=m^{(3)}_-v^{(3)}$ on
$\Sigma^{(3)}:=\{z :|z|=\rho_a,\rho_a^{-1}\}$, where $v^{(3)}=b_+^{(2)}$,
$|z|=\rho_a$ and $v^{(3)}=(b_-^{(2)})^{-1}$, $|z|=\rho_a^{-1}$. This
$\Sigma^{(3)}$ is not the best choice (see the next subsubsection). But for a
simple and direct estimate, we use this choice in this subsubsection. From the
choice of $\rho_a$, we have (see (5.13)--(5.14) of \cite{BDJ})
\begin{equation}
   |v^{(3)}(z;k;t)-I|\le e^{-ck}
\qquad \text{for all $z\in\Sigma^{(3)}$},
\end{equation}
which implies that $I-C_{w^{(3)}}$ is invertible
and the norm of the inverse is uniformly bounded,
where $w^{(3)}:=v^{(3)}-I$ and $C_{w^{(3)}}(f):=  C_-(fw^{(3)})$
on $L^2(\Sigma^{(3)},|dz|)$, $C_\pm$ being Cauchy operators
(see (2.5)--(2.9) in \cite{BDJ} and references therein).
From the general theory of RHP's, we have
\begin{equation}\label{as12.15}
      m^{(3)}(z)=I+\frac1{2\pi i} \int_{\Sigma^{(3)}}
\frac{((I-C_{w^{(3)}})^{-1}I)(s)w^{(3)}(s)}{s-z}ds,
\quad z\notin\Sigma^{(3)}.
\end{equation}
This implies the estimates (see (5.16) of \cite{BDJ})
\begin{equation}
   |m^{(3)}_{22}(0;k;t) -1|, \quad |m^{(3)}_{12}(0;k;t)| \le Ce^{-ck},
\end{equation}
which using \eqref{as12.13}, \eqref{as12.8}, \eqref{as12.9} and \eqref{as7.7},
\eqref{as7.8}, yield Proposition \ref{thm1} (i). This is precisely the result
contained in (5.17) of \cite{BDJ}.

\medskip
From \eqref{as12.13}, \eqref{as12.8}, \eqref{as12.10} and \eqref{as7.8},
we have
\begin{eqnarray}
\label{as12.17}
   \pi_k(z;t) &=& -(-1)^ke^{-tz}m^{(3)}_{12}(z;k;t),
\quad\quad\qquad\qquad\qquad\qquad |z|<\rho_a, \\
\label{as12.20}
   \pi_k(z;t) &=& z^{k}e^{-tz^{-1}}m^{(3)}_{11}(z;k;t),
\quad\qquad\qquad\qquad\quad\qquad\qquad |z|>\rho_a^{-1},\\
\label{as12.18}
   \pi_k(z;t) &=& z^ke^{-tz^{-1}}m^{(3)}_{11}(z;k;t)
-(-1)^ke^{-tz}m^{(3)}_{12}(z;k;t),
\quad \rho_a<|z|<\rho_a^{-1}.
\end{eqnarray}
Let $0<b<1$ be a fixed number. From Figure \ref{fig-sig}, we could have
chosen $\rho_a$ such that $\rho_a>b$. When $|z|\le b$ and $|z|\ge b^{-1}$,
we have $\dist(z,\Sigma^{(3)})\ge c>0$. Since $z$ is uniformly bounded away
from the contour, we can extend the argument leading to (5.17) in
\cite{BDJ} where the uniform boundedness of $0$ from the contour is
used. Hence using
\eqref{as12.15},
\eqref{as12.17} and \eqref{as12.20} imply that
\begin{eqnarray}
\label{as12.21}   |e^{tz}\pi_k(z;t)| &\le& Ce^{-ck},
\qquad |z|\le b, \\
\label{as12.21.5}   |e^{tz^{-1}}z^{-k}\pi_k(z;t)| &\le& Ce^{-ck},
\qquad |z|\ge b^{-1}.
\end{eqnarray}
These are \eqref{as7.42}, \eqref{as7.43} in Proposition \ref{thm4}
of the special case $x\ge 2^{1/3}(1-a)k^{2/3}$.

\medskip
On the other hand, let $L>0$ be a fixed number.
Set $\alpha=1-2^{4/3}k^{-1/3}w$ with $-L\le w\le L$ as
in Proposition \ref{thm4}.
Since $\rho_a$ is fixed, when $k$ is large,
$\dist(-\alpha,\Sigma^{(3)})\ge c>0$.
Then from \eqref{as12.15},
\begin{equation}\label{as12.22}
   |m^{(3)}_{11}(-\alpha;k;t)-1|, \quad |m^{(3)}_{12}(-\alpha;k;t)|
\le Ce^{-ck}.
\end{equation}
Note that
\begin{eqnarray}
  \frac12(s-s^{-1}) &\le& s-1, \quad s>0,\\
  -\frac12(s-s^{-1}) +\log s &\le& \frac23(1-s)^3,
\quad \frac12\le s\le1, \\
  -\frac12(s-s^{-1}) +\log s &\le& 0,
\quad s\ge 1.
\end{eqnarray}
Thus for $\gamma\le 1$, $s\ge\frac12$,
\begin{equation}\label{as12.24}
 \begin{split}
   F(-s;\gamma) &= -\frac{\gamma}2(s-s^{-1})+\log s
=\frac{1-\gamma}2(s-s^{-1})-\frac12(s-s^{-1}) +\log s\\
&\le (1-\gamma)(s-1) +\frac23|s-1|^3.
 \end{split}
\end{equation}
For large $k$, $\alpha\ge \frac12$ for all $-L\le w\le L$,
and hence
\begin{equation}\label{as12.25}
  \bigl|(-\alpha)^ke^{-t(\alpha-\alpha^{-1})}\bigr|
= e^{kF(-\alpha;\frac{2t}{k})}
\le e^{\frac{32}2|w|^3-k(1-\frac{2t}{k})\frac{2^{4/3}w}{k^{1/3}}}
\le e^{\frac{32}2L^3}e^{2^{4/3}Lk^{2/3}} = Ce^{ck^{2/3}}.
\end{equation}
Similarly, since $\alpha^{-1}\ge \frac12$ for large $k$,
\begin{equation}\label{as12.26}
  \bigl|(-\alpha)^{-k}e^{t(\alpha-\alpha^{-1})}\bigr|
= e^{kF(-\alpha^{-1};\frac{2t}{k})}
\le e^{k(1-\frac{2t}{k})\frac{2^{4/3}w}{k^{1/3}\alpha}
+ \frac{32}2|w|^3\alpha^{-3}}
\le Ce^{ck^{2/3}}.
\end{equation}
Therefore from \eqref{as12.18} and \eqref{as12.22},
\begin{eqnarray}
  \bigl|e^{-t\alpha}\pi_k(-\alpha;t)\bigr|
= \bigl|(-\alpha)^ke^{-t(\alpha-\alpha^{-1})}m^{(3)}_{11}(-\alpha;k;t)
-(-1)^km^{(3)}_{12}(-\alpha;k;t)\bigr|
&\le& Ce^{ck^{2/3}}, \\
  \bigl|e^{-t\alpha^{-1}}(-\alpha)^{-k}\pi_k(-\alpha;t) -1 \bigr|
= \bigl|m^{(3)}_{11}(-\alpha;k;t) -1
-\alpha^{-k}e^{t(\alpha-\alpha^{-1})}m^{(3)}_{12}(-\alpha;k;t)\bigr|
&\le& Ce^{-ck}.
\end{eqnarray}
Noting $x\sim k^{2/3}$, these are \eqref{as7.44}, \eqref{as7.45} in Proposition
\ref{thm4} of the special case $x\ge 2^{1/3}(1-a)k^{2/3}$. Thus we have
extended the argument of (5.17) of \cite{BDJ} to the case when $z\to -1$. This
is an example of the extension of the first category mentioned at the beginning
of Section \ref{RHP} (though it is straightforward to extend in this case).

\subsubsection{The case $ak\le 2t\le k-M2^{-1/3}k^{1/3}$ for some $0<a<1$
and $M> M_0$.}
\label{subsub2}

In the previous section, the contour $\Omega^{(3)}$
was not the best choice.
We could have chosen the steepest descent curve for
$F(z;\gamma)$.
For the previous case, it was not necessary to use the steepest descent
curve to obtain the desired results,
but for the case at hand and in the future calculations, we need to
use the steepest descent curve.

For fixed $\theta$ satisfying $0\le\theta<\pi/2$ or $(3\pi)/2<\theta<2\pi$,
$F(\rho,\theta;\gamma)$ is always negative for $0<\rho<1$, and as
$\rho\downarrow 0$, it decreases to minus infinity. On the other hand, one can
check that (see (5.5) in \cite{BDJ}) when $\gamma\le 1$, the minimum of
$F(\rho,\theta;\gamma)$, $0<\rho\le 1$ is attained, for fixed $\pi/2\le
\theta\le (3\pi)/2$, at
\begin{equation}\label{as12.28}
  \rho=\rho_{\theta} := \frac{1-\sqrt{1-\gamma^2\cos^2 \theta}}
{-\gamma\cos{\theta}},
\end{equation}
and $F(\rho_\theta,\theta;\gamma)<0$. And it is straightforward to check that
for $0\le\gamma\le 1$, $\pi/2\le\theta\le 3\pi/2$,
\begin{equation}
\label{as12.29}
   F(\rho_\theta,\theta;\gamma) = \sqrt{1-\gamma^2\cos^2\theta}
+ \log \biggl[ \frac{1-\sqrt{1-\gamma^2\cos^2\theta}}{-\gamma\cos\theta}
\biggr] \le -\frac{2\sqrt{2}}3(1+\gamma\cos\theta)^{3/2}.
\end{equation}
This is an extension of (5.13) in \cite{BDJ} where only the case $\theta=\pi$
was considered. We need this improved version to obtain a better $L^1$ estimate
of $v^{(3)}-I$ in the sequel. Also $F(\rho_\theta,\theta;\gamma)$ is increasing
in $\pi/2\le\theta\le\pi$ and is decreasing in $\pi\le\theta\le 3\pi/2$. In
fact, the saddle points for $\frac{\gamma}2(z-z^{-1})+\log z$ are $z=-\rho_\pi$
and $z=-\rho_\pi^{-1}$.

\begin{figure}[ht]
 \centerline{\epsfig{file=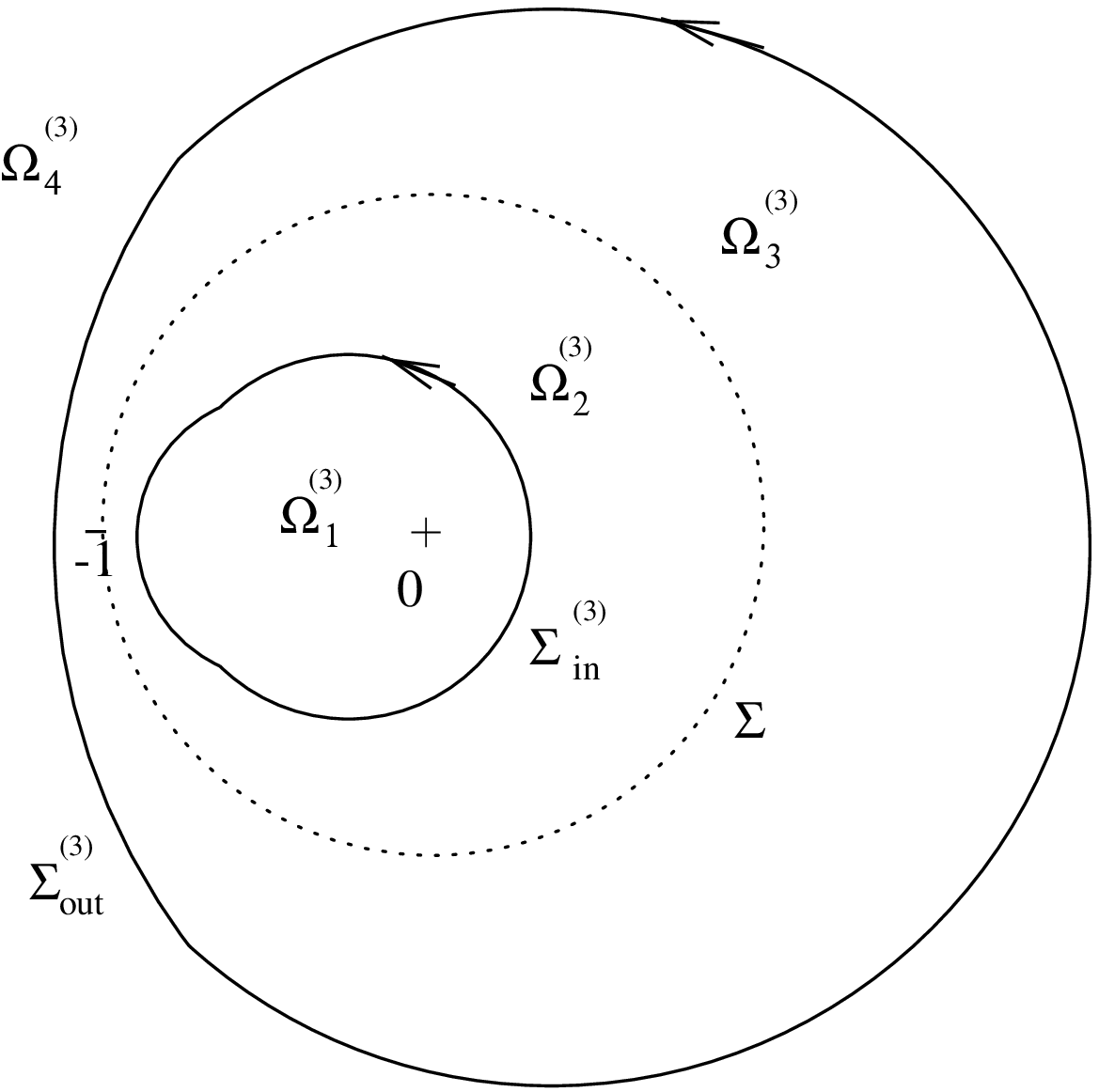, width=6cm}}
 \caption{$\Sigma^{(3)}$ and $\Omega^{(3)}$ when $\gamma<1$}
\label{fig-gammale1}
\end{figure}
In this time, define $\Sigma^{(3)}:=\Sigma^{(3)}_{in}\cup\Sigma^{(3)}_{out}$
by, as in (5.6) of \cite{BDJ},
\begin{equation}\label{as12.30}
  \begin{split}
   &\Sigma^{(3)}_{in}=\{\rho_{\theta}e^{i\theta} :
3\pi/4 \le \theta\le 5\pi/4 \} \cup\{\rho_{3\pi/4}e^{i\theta} :
0\le \theta\le 3\pi/4 , 5\pi/4\le\theta<2\pi \},\\
   &\Sigma^{(3)}_{out}=\{\rho_{\theta}^{-1}e^{i\theta} :
3\pi/4 \le \theta\le 5\pi/4 \} \cup\{\rho_{3\pi/4}^{-1}e^{i\theta} :
0\le \theta\le 3\pi/4 , 5\pi/4\le\theta<2\pi \},
  \end{split}
\end{equation}
where $\rho_\theta$ is defined in \eqref{as12.28} with
$\gamma=(2t)/k$. Orient $\Sigma^{(3)}$ as in Figure
\ref{fig-gammale1}. Note that $\Sigma^{(3)}$ lies in $\Omega_2^{sig}$ and
for $3\pi/4 \le \theta\le 5\pi/4$, it is the steepest descent curve. The
reason why we choose a part of the circle as the contour for the remaining
angles is to ensure the uniform boundedness of the Cauchy operators (see
Section 5 of \cite{BDJ}.) This does not affect the asymptotics since as
we will see the
main contribution to the asymptotics comes from neighborhood of $z=-1$.

Define the regions $\Omega^{(3)}_j$, $j=1,\cdots,4$ as in Figure
\ref{fig-gammale1}. Define $m^{(3)}(z;k;t)$ by, as in (5.9) of \cite{BDJ},
\begin{equation}\label{as12.31}
  \begin{cases}
     m^{(3)}=m^{(2)} (b_+^{(2)})^{-1},
\qquad &\text{in $\Omega^{(3)}_2$,}\\
     m^{(3)}=m^{(2)} (b_-^{(2)})^{-1},
\qquad &\text{in $\Omega^{(3)}_3$,}\\
     m^{(3)}=m^{(2)}
\qquad\qquad &\text{in $\Omega^{(3)}_1$, $\Omega^{(3)}_4$,}
  \end{cases}
\end{equation}
where $b_\pm^{(2)}$ are defined in \eqref{as12.7}.
Then $m^{(3)}$ solves a new RHP with the jump matrix $v^{(3)}(z;k;t)$
where
\begin{equation}\label{as12.36}
  \begin{cases}
    v^{(3)} =
\begin{pmatrix} 1&-(-1)^kz^ke^{t(z-z^{-1})}\\0&1
\end{pmatrix} \quad &\text{on}\quad\Sigma_{in}^{(3)},\\
    v^{(3)} =
\begin{pmatrix} 1&0\\(-1)^kz^{-k}e^{-t(z-z^{-1})}&1
\end{pmatrix} \quad &\text{on}\quad\Sigma_{out}^{(3)}.
  \end{cases}
\end{equation}
Set $w^{(3)}:=v^{(3)}-I$.
For $z\in\Sigma^{(3)}_{in}$, from the choice of $\Sigma^{(3)}$ and \eqref{as12.29},
the (12)-entry of the jump matrix satisfies
for $3\pi/4\le \arg z \le 5\pi/4$,
\begin{equation}\label{as12.33}
  |z^ke^{t(z-z^{-1})}| = e^{kF(\rho_\theta,\theta;\frac{2t}{k})}
\le e^{-\frac{2\sqrt{2}}{3}k(1+\frac{2t}{k}\cos\theta)^{3/2}}
\le e^{-\frac{2\sqrt{2}}{3}k(1-\frac{2t}{k})^{3/2}}
\le e^{-\frac{2}{3}M^{3/2}},
\end{equation}
and for $0\le \arg z\le 3\pi/4$ or $5\pi/4\le\arg z< 2\pi$,
\begin{equation}\label{as12.34}
  |z^ke^{t(z-z^{-1})}| = e^{kF(\rho_{3\pi/4},\theta;\frac{2t}{k})}
\le e^{kF(\rho_{3\pi/4},\frac{3\pi}4;\frac{2t}{k})} \le
e^{-\frac{2\sqrt{2}}3k(1+\frac{2t}{k}\cos\frac{3\pi}{4})^{3/2}} \le
e^{-\frac1{24}k}.
\end{equation}
From \eqref{as12.12}, similar estimates hold for $z^{-k}e^{-t(z-z^{-1})}$,
$z\in\Sigma^{(3)}_{out}$. Thus we have
\begin{equation}\label{qz10.38}
  \| w^{(3)}\|_{L^\infty(\Sigma^{(3)})} \le
  Ce^{-\frac{2\sqrt2}{3}k(1-\frac{2t}{k})^{3/2}}.
\end{equation}
Also there exists $M_0$ such that for $M>M_0$,
$\|C_{w^{(3)}}\|_{L^2(\Sigma^{(3)})\to L^2(\Sigma^{(3)})} \le c_1< 1$, and
hence \eqref{as12.15} holds. This is precisely (5.18) of \cite{BDJ}. For this
derivation, we do not need the extension \eqref{as12.29} of (5.13) in
\cite{BDJ}. But for the improved $L^1$ norm estimate of $w^{(3)}$, which we do
now, we need \eqref{as12.29}.

Note that $|dz|\le C|d\theta|$ on $\Sigma^{(3)}$. Using the estimates in
\eqref{as12.33} and \eqref{as12.34},
\begin{equation}\label{as12.35}
   \int_{\Sigma^{(3)}_{in}} \bigl|z^ke^{t(z-z^{-1})}\bigr| |dz|
\le C\int_{\frac{3\pi}4}^{\frac{5\pi}4}
e^{-\frac{2\sqrt{2}}{3}k(1+\frac{2t}{k}\cos\theta)^{3/2}}
d\theta
+ C\int_{[0,2\pi)\setminus [3\pi/4,5\pi/4]}
e^{-\frac1{24}k} d\theta.
\end{equation}
The second integral is clearly less than $Ce^{-\frac1{24}k}$. For the first
integral, recall the inequality $(x+y)^a\ge x^a+y^a$, $x,y>0$, $a\ge 1$. Then
using the inequality $1+\cos\theta\ge \frac1{2\sqrt{2}}(\theta-\pi)^2$ for
$\theta\in[\frac{3\pi}4,\frac{5\pi}4]$, together with the condition $ak\le 2t$,
the first integral is less than or equal to
\begin{equation}
   \int_{\frac{3\pi}4}^{\frac{5\pi}4}
e^{-\frac{2\sqrt{2}}{3}k\bigl[\bigl(1-\frac{2t}{k}\bigr)^{3/2} +\bigl(
\frac{2t}{k}(1+\cos\theta)\bigr)^{3/2}\bigr]} d\theta \le
e^{-\frac{2\sqrt{2}}{3}k\bigl(1-\frac{2t}{k}\bigr)^{3/2}}
\int_{\frac{3\pi}4}^{\frac{5\pi}4}
e^{-\frac{a^{3/2}}{3\cdot2^{3/4}}k|\theta-\pi|^3} d\theta,
\end{equation}
where last inequality is less than or equal to $Ck^{-1/3}$ for some constant
$C>0$. Therefore adjusting constants, we obtain
\begin{equation}
   \int_{\Sigma^{(3)}_{in}} \bigl|z^ke^{t(z-z^{-1})}\bigr| |dz|
\le \frac{C}{k^{1/3}}
e^{-\frac{2\sqrt{2}}{3}k\bigl(1-\frac{2t}{k}\bigr)^{3/2}}.
\end{equation}
We have similar estimates on $\Sigma^{(3)}_{out}$. Therefore,
\begin{equation}\label{qz10.41}
  \|w^{(3)}\|_{L^1(\Sigma^{(3)})}
\le \frac{C}{k^{1/3}}
e^{-\frac{2\sqrt{2}}{3}k\bigl(1-\frac{2t}{k}\bigr)^{3/2}}.
\end{equation}
This is a refinement of (5.23) in \cite{BDJ}.

Now from \eqref{as12.15}, we have
\begin{equation}\label{as12.39}
      m^{(3)}(z)=I+
\frac1{2\pi i}\int_{\Sigma^{(3)}} \frac{w^{(3)}(s)}{s-z} ds
+ \frac1{2\pi i} \int_{\Sigma^{(3)}}
\frac{[(I-C_{w^{(3)}})^{-1}C_{w^{(3)}}I](s)w^{(3)}(s)}{s-z}ds,
\quad z\notin\Sigma^{(3)}.
\end{equation}
In \cite{BDJ}, only the term $m^{(3)}_{22}(z)$ (at $z=0$) was of interest. But
then $w^{(3)}_{22}=0$, and the first integral in \eqref{as12.39} was zero. As
computed in (5.20) of \cite{BDJ}, the second integral is bounded by the product
of the $L^\infty$ and $L^1$ norms of $w^{(3)}$, and then due to
\eqref{qz10.38}, the estimate $\|w^{(3)}\|_{L^1} \le Ck^{-1/3}$ in (5.23) of
\cite{BDJ} was enough to control $m^{(3)}_{22}$. But in the present paper, we
need estimates of $m^{(3)}_{12}$ for $\pi_k(z)$, and hence we need an estimate
of the first integral which is same as a bound on the $L^1$ norm of $w^{(3)}$.
The $L^1$ bound (5.23) obtained in \cite{BDJ} is not good enough for this
purpose, and we need an improved estimate on the $L^1$ norm of $w^{(3)}$. Now
by \eqref{qz10.41}, the first integral is less than or equal to
\begin{equation}\label{as12.40}
   \frac{C}{\dist(z,\Sigma^{(3)})k^{1/3}}
e^{-\frac{2\sqrt{2}}{3}k\bigl(1-\frac{2t}{k}\bigr)^{3/2}},
\end{equation}
while, as in (5.20) of \cite{BDJ}, the second integral is less than or equal
to, by using \eqref{qz10.38} and \eqref{qz10.41},
\begin{equation}\label{as12.41}
 \begin{split}
   &\frac1{2\pi \dist(z,\Sigma^{(3)})}
\|(I-C_{w^{(3)}})^{-1}C_{w^{(3)}}I\|_{L^2}\|w^{(3)}\|_{L^2} \\
   &\quad \le \frac1{2\pi \dist(z,\Sigma^{(3)})}
\|(I-C_{w^{(3)}})^{-1}\|_{L^2\rightarrow L^2}
\|C_{w^{(3)}}I\|_{L^2}\|w^{(3)}\|_{L^2} \\
   &\quad \le \frac{C}{\dist(z,\Sigma^{(3)})}
\|w^{(3)}\|_{L^2}^2 \\
   &\quad \le \frac{C}{\dist(z,\Sigma^{(3)})}
\|w^{(3)}\|_{L^\infty} \|w^{(3)}\|_{L^1}\\
   &\quad \le \frac{C}{\dist(z,\Sigma^{(3)})k^{1/3}}
e^{-\frac{4\sqrt{2}}{3}k\bigl(1-\frac{2t}{k}\bigr)^{3/2}}.
 \end{split}
\end{equation}

\medskip
When $z=0$, $\dist(z,\Sigma^{(3)}) \ge c_1>0$,
hence using \eqref{as12.31}, \eqref{as12.8}, \eqref{as12.9}
and \eqref{as7.7}, \eqref{as7.8},
we obtain Proposition \ref{thm1} (ii).

\medskip
As in \eqref{as12.17}--\eqref{as12.18},
from \eqref{as12.13}, \eqref{as12.8}, \eqref{as12.10} and \eqref{as7.8},
we have
\begin{eqnarray}
\label{as12.42}
   \pi_k(z;t) &=& -(-1)^ke^{-tz}m^{(3)}_{12}(z;k;t),
\quad\quad\qquad\qquad\qquad\qquad z\in\Omega^{(3)}_1, \\
\label{as12.43}
   \pi_k(z;t) &=& z^{k}e^{-tz^{-1}}m^{(3)}_{11}(z;k;t),
\quad\qquad\qquad\qquad\quad\qquad\qquad z\in\Omega^{(3)}_4,\\
\label{as12.44}
   \pi_k(z;t) &=& z^ke^{-tz^{-1}}m^{(3)}_{11}(z;k;t)
-(-1)^ke^{-tz}m^{(3)}_{12}(z;k;t),
\quad z\in\Omega^{(3)}_2\cup\Omega^{(3)}_3.
\end{eqnarray}
Define $x$ by
\begin{equation}
   \frac{2t}{k}= 1-\frac{x}{2^{1/3}k^{2/3}}
\end{equation}
as in Proposition \ref{thm4}.
Let $0<b<1$ be a fixed number.
Given $b$, from the beginning, we could have chosen $0<a<1$
such that
\begin{equation}\label{as12.50}
  \rho_{\theta_b} = \frac{1-\sqrt{1-a^2\cos^2\theta_b}}{-a\cos\theta_b}
\end{equation}
is strictly greater than $b$ for some $\pi/2\le \theta_b <\pi$. Note that in
\eqref{as12.30} defining $\Sigma^{(3)}$, the choice of $3\pi/4$ and $5\pi/4$
was arbitrary. Instead of $3\pi/4$ and $5\pi/4$, this time we use
$\theta_b$ and $2\pi-\theta_b$, and carry this forward through the later
calculations.  Thus we obtain the same estimates of \eqref{as12.40} and
\eqref{as12.41} with different constants $C$. Now $|z|\le b$ lies in
$\Omega^{(3)}_1$ and $|z|\ge b^{-1}$ lies in $\Omega^{(3)}_4$. Since the
distance $\dist(z,\Sigma^{(3)})\ge c_2>0$, using \eqref{as12.39},
\eqref{as12.40}, \eqref{as12.41},
\eqref{as12.42} and \eqref{as12.43}, we obtain \eqref{as7.42} and
\eqref{as7.43} in Proposition \ref{thm4} for $M\le x\le (1-a)2^{1/3}k^{2/3}$.
Since in \eqref{as12.21} and \eqref{as12.21.5} in Subsubsection \ref{subsub1},
the choice of $0<a<1$ was arbitrary, we obtain \eqref{as7.42} and
\eqref{as7.43} in Proposition \ref{thm4} for all $x\ge M$.

\medskip
On the other hand, let $0<L<2^{-3/2}\sqrt{M}$ be a fixed number.
Set $\alpha=1-2^{4/3}k^{-1/3}w$ with $-L\le w\le L$ as
in Proposition \ref{thm4}.
From the inequality $\frac{1-\sqrt{1-\gamma^2}}{\gamma}
\le 1-\sqrt{1-\gamma}$ for all $0\le \gamma\le 1$,
we have
\begin{equation}
  \rho_\pi=\frac{1-\sqrt{1-\bigl(\frac{2t}{k}\bigr)^2}}{\frac{2t}{k}}
\le 1-\sqrt{1-\frac{2t}{k}}
\le 1-\frac{\sqrt{M}}{2^{1/6}k^{1/3}}.
\end{equation}
But
\begin{equation}
  \alpha = 1-\frac{2^{4/3}w}{k^{1/3}} \ge 1-\frac{2^{4/3}L}{k^{1/3}}.
\end{equation}
Hence $\dist(-\alpha, \Sigma^{(3)}) \ge Ck^{-1/3}$.
Thus from \eqref{as12.39}, \eqref{as12.40} and \eqref{as12.41}, we have
\begin{equation}
   |m^{(3)}(-\alpha;k;t)-I| \le Ce^{-\frac{2}{3}x^{3/2}},
\end{equation}
which together with \eqref{as12.44}
(note that $-\alpha\in \Omega^{(3)}_2\cup\Omega^{(3)}_3$),
implies that
\begin{eqnarray}
\label{as12.53}
\bigl| e^{-t\alpha}\pi_k(-\alpha;t) \bigr| &\le&
\bigl| (-\alpha)^ke^{-t(\alpha-\alpha^{-1})} \bigr|
\bigl(1+Ce^{-c|x|^{3/2}}\bigr) + Ce^{-c|x|^{3/2}}, \\
\label{as12.54}
\bigl| e^{-t\alpha^{-1}}(-\alpha)^{-k} \pi_k(-\alpha;t) -1 \bigr|
&\le& Ce^{-c|x|^{3/2}} +
\bigl| (-\alpha)^{-k}e^{t(\alpha-\alpha^{-1})} \bigr|
Ce^{-c|x|^{3/2}}.
\end{eqnarray}
For large $k$, $\alpha\ge \frac12$ for all $-L\le w\le L$,
hence using \eqref{as12.24}, we obtain as in \eqref{as12.25}
and \eqref{as12.26},
\begin{equation}
  \bigl|(-\alpha)^ke^{-t(\alpha-\alpha^{-1})}\bigr|
= e^{kF(-\alpha;\frac{2t}{k})}
\le e^{-2wx + \frac{32}2|w|^3}
\le Ce^{c|x|},
\end{equation}
and
\begin{equation}
  \bigl|(-\alpha)^{-k}e^{t(\alpha-\alpha^{-1})}\bigr|
= e^{kF(-\alpha^{-1};\frac{2t}{k})}
\le Ce^{c|x|}.
\end{equation}
Thus from \eqref{as12.53} and \eqref{as12.54}, we obtain \eqref{as7.44},
\eqref{as7.45} in Proposition \ref{thm4}.

\subsubsection{The case $k-M2^{-1/3}k^{1/3}\le 2t\le k$ for some $M>0$.}
\label{subsub3}

In this case, as $k\to\infty$ the point $z=-\rho_\pi$ on the deformed
contour $\Sigma^{(3)}$ defined in \eqref{as12.30} approaches $z=-1$
rapidly, and so we need to pay special attention to the neighborhood of
$z=-1$. More precisely, we need to introduce the so-called parametrix for
the RHP around $z=-1$, which is an approximate local solution.

\begin{figure}[ht]
 \centerline{\epsfig{file=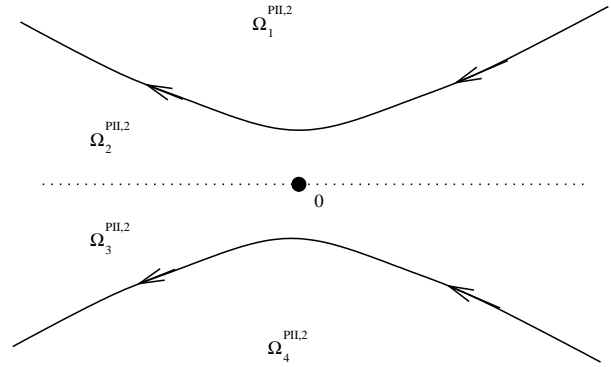, width=8cm}}
 \caption{$\Sigma^{PII,2}$ and $\Omega^{PII,2}_j$}\label{fig-PII2}
\end{figure}

Recall the RHP \eqref{as20} for the Painlev\'e II equation. Let
$\Sigma^{PII,2}=\Sigma^{PII,2}_1\cup\Sigma^{PII,2}_2$ be a contour of the
general shape indicated in Figure \ref{fig-PII2}. Asymptotically for large $z$,
the curves are straight lines of angle less than $\pi/3$. See the paragraph
after (2.18) in \cite{BDJ} for more precise discussions on the curve. We will
define the exact shape of $\Sigma^{PII,2}$ below. Define $m^{PII,2}(z;x)$ by
\begin{equation}\label{as12.59}
  \begin{cases}
    m^{PII,2}(z,x) = m(z;x)
\begin{pmatrix} 1&e^{-2i(\frac43z^3+xz)}\\0&1 \end{pmatrix}
&\text{in $\Omega^{PII,2}_2$},\\
    m^{PII,2}(z,x) = m(z;x)
\begin{pmatrix} 1&0\\e^{2i(\frac43z^3+xz)}&1 \end{pmatrix}
&\text{in $\Omega^{PII,2}_3$},\\
    m^{PII,2}(z,x) = m(z;x)
&\text{in $\Omega^{PII,2}_1,\Omega^{PII,2}_4$},
  \end{cases}
\end{equation}
where $m(z;x)$ is the solution of the RHP for the PII equation given in
\eqref{as20}.  Then $m^{PII,2}$ solves a new RHP (see (2.19) of \cite{BDJ})
\begin{equation}
  \begin{cases}
    m^{PII,2} \qquad\quad\text{is analytic in $\C\setminus\Sigma^{PII,2}$},\\
    m^{PII,2}_+ = m^{PII,2}_-
\begin{pmatrix} 1&-e^{-2i(\frac43z^3+xz)}\\0&1 \end{pmatrix}
\quad\text{in $\Sigma^{PII,2}_1$},\\
    m^{PII,2}_+ = m^{PII,2}_-
\begin{pmatrix} 1&0\\e^{2i(\frac43z^3+xz)}&1 \end{pmatrix}
\quad\text{in $\Sigma^{PII,2}_2$},\\
    m^{PII,2} = I+O\bigl(\frac1{z}\bigr) \qquad\text{as $z\to\infty$}.
  \end{cases}
\end{equation}
Also $m^{PII,2}_1(x)$ defined by
$m^{PII,2}(z;x)=I+\frac{m^{PII,2}_1(x)}{z}+O(z^{-2})$ satisfies
\begin{equation}
  m^{PII,2}_1(x)=m_1(x),
\end{equation}
where $m_1(x)$ is defined in a similar manner in \eqref{as7.27}.

Set $x$ by
\begin{equation}
   \frac{2t}{k}= 1-\frac{x}{2^{1/3}k^{2/3}}.
\end{equation}
We define $\Sigma^{(3)}$ and $m^{(3)}$ as in \eqref{as12.30} and
\eqref{as12.31}. Now we proceed as in (5.25)--(5.35) of \cite{BDJ}. Let
$\mathcal{O}$ be the ball of radius $\epsilon$ around $z=-1$, where
$\epsilon>0$ is a small fixed number. Define the map (see the equation
displayed between (5.26) and (5.27) of \cite{BDJ})
\begin{eqnarray}\label{as12.59.5}
   \lambda(z) := -i2^{-4/3}k^{1/3}\frac12(z-z^{-1})
\end{eqnarray}
in $\mathcal{O}$. Define $\Sigma^{PII,2}$ by
$\Sigma^{PII,2}\cap\lambda(\mathcal{O}) :=\lambda(\Sigma^{(3)}\cap\mathcal{O})$
and extend it smoothly outside $\lambda(\mathcal{O})$ as indicated in (5.29) of
\cite{BDJ}. And define $m^{PII,2}$ as above using this contour. Now we define
the parametrix by (see the equation displayed between (5.30) and (5.31) of
\cite{BDJ})
\begin{equation}\label{as12.60}
  \begin{cases}
    m_p(z;k;t)=m^{PII,2}(\lambda(z),x) \quad
&\text{in $\mathcal{O}\setminus\Sigma^{(3)}$,} \\
    m_p(z;k;t)=I \quad
&\text{in $\mathcal{\bar{O}}^c\setminus\Sigma^{(3)}$.}
  \end{cases}
\end{equation}
It is proved in (5.25)--(5.34) in \cite{BDJ} that
if we take $\epsilon$ small enough but fix it,
then the ratio
\begin{equation}\label{as12.61}
   R(z;k;t):= m^{(3)}m_p^{-1}
\end{equation}
solves a new RHP
\begin{equation}
  \begin{cases}
R(z;k;t) \qquad\quad\text{is analytic in
$\C\setminus\Sigma_R$,}\\
R_+(z;k;t)=R_-(z;k;t)v_R(z;k;t), \quad
\text{on $\Sigma_R$,}\\
R(z;k;t) = I+O(\frac1{z}), \qquad\text{as $z\to\infty$,}
  \end{cases}
\end{equation}
where $\Sigma_R:=\partial\mathcal{O}\cup\Sigma^{(3)}$, and the jump matrix
satisfies (see (5.34) of \cite{BDJ})
\begin{equation}\label{as12.62}
  \begin{cases}
    \|v_R-I\|_{L^\infty}\le\frac{C}{k^{2/3}}
&\text{on $\mathcal{O} \cap \Sigma^{(3)}$},\\
    \|v_R-I\|_{L^\infty}\le Ce^{-ck}
&\text{on $\mathcal{O}^c \cap \Sigma^{(3)}$},\\
    \|v_R-I+\frac{m^{PII,2}_1(x)}{\lambda(z)}\|_{L^\infty}
\le \frac{C}{k^{2/3}}
&\text{on $\partial \mathcal{O}$, as $k\to\infty$,}
  \end{cases}
\end{equation}
with some positive constants $C$ and $c$ which may depend on $M$. Set
$w_R:=v_R-I$. Using \eqref{as12.15} which holds generally, we have (see (5.35)
and the preceding calculations in \cite{BDJ})
\begin{equation}
 \begin{split}
   R(z;k;t)&=I+\frac1{2\pi i} \int_{\Sigma_R}
\frac{((I-C_{w_R})^{-1}I)(s)(w_R(s))}{s-z} ds \\
&=I+ \frac1{2\pi i} \int_{\Sigma_R} \frac{v_R(s)-I}{s-z} ds + \frac1{2\pi i}
\int_{\Sigma_R} \frac{[(I-C_{w_R})^{-1}C_{w_R}I](s)w_R(s)}{s-z} ds.
 \end{split}
\end{equation}
Now the absolute value of the second integral is less than or equal to (recall
that $|\lambda(z)|=O(k^{-1/3})$ for $z\in\partial\mathcal{O}$)
\begin{equation}
 \begin{split}
    &\frac{C}{\dist(z,\Sigma_R)}
\|(I-C_{w_R})^{-1}\|_{L^2(\Sigma_R)\to L^2(\Sigma_R)}
\|C_{w_R}I\|_{L^2(\Sigma_R)} \|w_R\|_{L^2(\Sigma_R)} \\
    &\le \frac{C}{\dist(z,\Sigma_R)}
\|w_R\|_{L^2(\Sigma_R)}^2 \\
    &\le \frac{\epsilon C}{\dist(z,\Sigma_R)k^{2/3}},
 \end{split}
\end{equation}
 and similarly, the first integral satisfies
\begin{equation}
    \biggl| \frac1{2\pi i} \int_{\Sigma_R}
\frac{v_R(s)-I}{s-z} ds
+ \frac{m^{PII,2}_1(x)}{2\pi i}
\int_{\partial\mathcal{O}} \frac1{\lambda(s)(s-z)} ds \biggr|
\le \frac{\epsilon C}{\dist(z,\Sigma_R)k^{2/3}}.
\end{equation}
Hence
\begin{equation}\label{as12.63}
   \biggl| m^{(3)}(z;k;t)\bigl(m_p(z;k;t)\bigr)^{-1}-I
+ \frac{m^{PII,2}_1(x)}{2\pi i}
\int_{\partial\mathcal{O}} \frac1{\lambda(s)(s-z)} ds \biggr|
\le \frac{\epsilon C}{\dist(z,\Sigma_R)k^{2/3}}.
\end{equation}
This is an extension of (5.35) to the case when $z\neq 0$.

For $z=0$, from \eqref{as12.60} and \eqref{as12.61},
$R(0)=m^{(3)}(0)$, and $\dist(0,\Sigma_R)\ge c_1>0$.
Note that $\lambda(s)$ is analytic in $\mathcal{O}$
except at $s=-1$, and
\begin{equation}
   \lambda(s) = -i2^{-4/3}k^{1/3}[(s+1)+\frac12(s+1)^2+\cdots],
\quad s\sim -1.
\end{equation}
By a residue calculation for \eqref{as12.63}, we have, as in (5.35) of
\cite{BDJ})
\begin{equation}\label{as12-71}
    m^{(3)}(0;k;t)
= I+ \frac{i2^{4/3}m_1^{PII,2}(x)}{k^{1/3}} + O\bigl( \frac1{k^{2/3}}\bigr).
\end{equation}
Thus using \eqref{as22} and \eqref{as23}, from \eqref{as7.7}, \eqref{as7.8},
\eqref{as12.8}, \eqref{as12.9} and \eqref{as12.31}, we obtain Proposition
\ref{thm1} (iii) of the case when $0\le x\le M$.

We now prove Proposition \ref{thm2} when $x\ge 0$. Since the choice of $M$ was
arbitrary in our calculations, for fixed $x$ we choose $M>0$ large enough so
that $x <M$. Let $z\in\C\setminus\Sigma$ be fixed. We first assume $|z|<1$. By
modifying the contour $\Sigma^{(3)}$, if necessary, as in \eqref{as12.50} and
the following paragraph, we have $z\in\Omega^{(3)}_1$ and
$\dist(z,\Sigma^{(3)})\ge c_1>0$. Thus from \eqref{as12.63}, $|R(z)-I|\le
Ck^{-1/3}$ with some constant $C$ which depends on $x$. Thus from
\eqref{as7.8}, \eqref{as12.8}, \eqref{as12.31} and \eqref{as12.61}, we obtain
the first limit of \eqref{as27}. Similar calculation applies to the case
$|z|>1$, and we obtain the first limit of \eqref{as28}. The second limits of
\eqref{as27} and \eqref{as28} follow from the first limits of \eqref{as28} and
\eqref{as27} respectively by replacing $z\to 1/z$. Hence this extends the
calculation (5.35) of \cite{BDJ} to the case when $z$ is bounded away from the
contour.

Finally, we prove Proposition \ref{thm3} when $x> 0$.
Set
\begin{equation}\label{as13.72}
  \alpha= 1-\frac{2^{4/3}w}{k^{1/3}},
\qquad \text{$w$ fixed,}
\end{equation}
and
\begin{equation}
  \frac{2t}{k} = 1-\frac{x}{2^{1/3}k^{2/3}},
\qquad \text{$x>0$ fixed.}
\end{equation}
In this case, $-\alpha\in\mathcal{O}$.
By a residue calculation again, for $w$ not equal to $0$,
\begin{equation}
 \begin{split}
   \frac1{2\pi i}
\int_{\partial\mathcal{O}} \frac1{\lambda(s)(s+\alpha)} ds
&= \frac{i2^{4/3}}{(-1+\alpha)k^{1/3}}
+ \frac{1}{\lambda(-\alpha)} \\
&= -\frac{i}{w} + \frac{i}{w+\frac{2^{1/3}}{k^{1/3}}w^2+\cdots}
= O\bigl(\frac1{k^{1/3}}\bigr).
 \end{split}
\end{equation}
When $w=0$, we have the same order $O(k^{-1/3})$ by a similar calculation. On
the other hand, since
\begin{equation}
  \rho_\pi=
\frac{1-\sqrt{1-\bigl(\frac{2t}{k}\bigr)^2}}{\frac{2t}{k}}
=1-\frac{2^{1/3}\sqrt{x}}{k^{1/3}}
+ \frac{x}{2^{1/3}k^{2/3}}+ O\bigl(\frac1{k}\bigr),
\end{equation}
using \eqref{as13.72} we have $\dist(-\alpha,\Sigma_R)\ge Ck^{-2/3}$. Thus we
obtain from \eqref{as12.63}
\begin{equation}\label{as12-76}
   |R(-\alpha;k;t) - I| \le \epsilon C.
\end{equation}
Using $\lambda(-\alpha)\sim -iw$,
from \eqref{as12.60} and \eqref{as12.61}, we have
\begin{equation}
  \lim_{k\to\infty} m^{(3)}(-\alpha;k;t) = m^{PII,2}(-iw,x)
\end{equation}
since $\epsilon$ is arbitrarily small. On the other hand, from the conditions
on $t$ and $\alpha$, we have
\begin{equation}\label{as12.79}
   \lim_{k\to\infty} \alpha^ke^{-t(\alpha-\alpha^{-1})}
= e^{\frac83w^3-2xw}.
\end{equation}
Thus using \eqref{as12.31},  we obtain
\begin{eqnarray}
  \lim_{k\to\infty} m^{(2)}(-\alpha;k;t)&=&
m^{PII,2}(-iw,x),
\quad\qquad -\alpha\in\Omega^{(3)}_1,\Omega^{(3)}_4,\\
  \lim_{k\to\infty} m^{(2)}(-\alpha;k;t)&=&
m^{PII,2}(-iw,x) \begin{pmatrix} 1&-e^{\frac83w^3-2xw}\\0&1
\end{pmatrix}, \quad -\alpha\in\Omega^{(3)}_2,\\
  \lim_{k\to\infty} m^{(2)}(-\alpha;k;t)&=&
m^{PII,2}(-iw,x) \begin{pmatrix} 1&0\\-e^{-\frac83w^3+2xw}&1
\end{pmatrix}, \quad -\alpha\in\Omega^{(3)}_3.
\end{eqnarray}
Now finally using \eqref{as12.59}, for each fixed $w$ and $x$, we have
\begin{equation}\label{as13.82}
  \lim_{k\to\infty} m^{(2)}(-\alpha;k;t)
= m(-iw;x).
\end{equation}
From \eqref{as12.8}--\eqref{as12.11} and \eqref{as12.79}, this implies
Proposition \ref{thm3} of the case when $x>0$. This was a new computation we
had to do in the present paper in order to include the case when $\alpha \to
1$.

%
%

\subsection{When $(2t)/k>1$.}\label{sub2}

Throughout this subsection, we set
\begin{equation}
   \gamma:=\frac{2t}{k}>1.
\end{equation}
We need some definitions from \cite{BDJ}. Set $0<\theta_c<\pi$ by
$\sin^2\frac{\theta_c}2=\frac1\gamma$. Define a probability measure on an arc
(see (4.13) of \cite{BDJ}),
\begin{equation}\label{as12-84}
  d\mu(\theta) := \frac{\gamma}{\pi} \cos\bigl(\frac{\theta}2\bigr)
\sqrt{\frac1{\gamma}-\sin^2\bigl(\frac{\theta}2\bigr)}
d\theta,
\qquad -\theta_c\le\theta\le\theta_c,
\end{equation}
and define a constant (see (4.14) of \cite{BDJ})
\begin{equation}
  l:= -\gamma+\log\gamma+1.
\end{equation}
Now we introduce the so-called $g$-function (see (4.8) of \cite{BDJ}),
\begin{equation}\label{as13.86}
  g(z;k;t) := \int_{-\theta_c}^{\theta_c}
\log(z-e^{i\theta}) d\mu(\theta),
\qquad z\in\C\setminus\Sigma\cup (-\infty,-1].
\end{equation}
The measure $d\mu(\theta)$ is the equilibrium measure of a certain variational
problem and the constant $l$ is a related constant (see Section 4 of
\cite{BDJ}). For each $|\theta|\le \theta_c$, the branch is chosen such that
$\log(z-e^{i\theta})$ is analytic in
$\C\setminus(-\infty,-1]\cup\{e^{i\phi} : -\pi\le\phi\le\theta\}$ and
behaves like $\log z$ as $z\in\R\to+\infty$. The basic properties of $g(z)$
are summarized in Lemma 4.2 of \cite{BDJ}. In general, the role of the
$g$-function in RHP analysis, first introduced in \cite{DZ2} and then
generalized in \cite{DVZ}, is to replace exponentially growing terms in the
jump matrix by oscillating or exponentially decaying terms. In
\cite{DKMVZ}, the authors introduced a $g$-function of a form similar to
\eqref{as13.86} to analyze an RHP associated to orthogonal polynomials on the
real line. The above $g$-function \eqref{as13.86} introduced in \cite{BDJ}
is an adaptation of their work to the circle case. When $0\le \gamma\le 1$,
the related equilibrium measure is (see (4.12) of \cite{BDJ})
\begin{equation}
   d\mu(\theta)=\frac1{2\pi} (1+\gamma\cos\theta)d\theta,
\qquad -\pi\le\theta<\pi,
\end{equation}
with the related constant $l=0$, and hence (see (4.15) of \cite{BDJ})
\begin{equation}
   g(z) = \begin{cases} \log z-\frac\gamma{2z},
\quad &|z|>1,z\notin (-\infty,-1)\\
-\frac\gamma2 z+\pi i, \quad&|z|<1.\end{cases}
\end{equation}
Since $g(z)$ is explicit in this case, we did not introduce it of
the form \eqref{as13.86} in the previous subsection.

Then up to \eqref{qz10.104}, we follows the procedure in Section 6 of
\cite{BDJ}. Write $\Sigma=C_1\cup\overline{C_2}$ where $C_2:=\{ e^{i\theta} :
-\theta_c<\theta<\theta_c \}$ and $C_1:= \Sigma\setminus\overline{C_2}$. Define
$m^{(1)}(z;k;t)$ by
\begin{equation}
   m^{(1)}(z;k;t) :=
e^{\frac{kl}2\sigma_3}Y(z;k;t)
e^{-kg(z;k;t)\sigma_3}e^{-\frac{kl}2\sigma_3},
\end{equation}
where $\sigma_3=\bigl( \begin{smallmatrix} 1&0\\0&-1
\end{smallmatrix}\bigr)$.
Then $m^{(1)}$ solves
(see (6.1) in \cite{BDJ}) a new RHP
\begin{equation}
\begin{cases}
    m^{(1)}(z;k;t) \qquad\text{is analytic in}
\quad \C\setminus\Sigma,\\
    m^{(1)}_+(z;k;t) = m^{(1)}_-(z;k;t)
\begin{pmatrix}e^{-2k\wt{\alpha}(z;k;t)}& (-1)^k\\ 0& e^{2k\wt{\alpha}(z;k;t)}
\end{pmatrix} \ \ \text{on} \ \ C_2,\\
    m^{(1)}_+(z;k;t) = m^{(1)}_-(z;k;t)
\begin{pmatrix} 1& (-1)^ke^{-2k\wt{\alpha}-(z;k;t)}\\ 0&1 \end{pmatrix}
 \ \ \text{on} \ \ C_1,\\
    m^{(1)}(z;k;t)=I+O\bigl(\frac1z\bigr)
\quad\text{as}\quad z\rightarrow \infty,
\end{cases}
\end{equation}
where $\wt{\alpha}(z;k;t)$ is defined by (see Lemma 6.1 in \cite{BDJ})
\begin{equation}
   \wt{\alpha}(z;k;t):= -\frac{\gamma}4 \int^z_{e^{i\theta_c}}
\frac{s+1}{s^2}\sqrt{(s-e^{i\theta_c})(s -e^{-i\theta_c})}ds,
\qquad \xi:=e^{i\theta_c}.
\end{equation}
(Notation: We use $\wt{\alpha}$ here instead of $\alpha$ in \cite{BDJ}
to avoid confusion with $\alpha$ in \eqref{as12-109} below.)
The branch is chosen such that $\sqrt{(s-e^{i\theta_c})(s -e^{-i\theta_c})}$
is analytic in $\C\setminus\overline{C_1}$
and behaves like $s$ as $s\in\R\to+\infty$.

Define $m^{(2)}(z;k;t)$ as in \eqref{as12.4}. Then $m^{(2)}$ solves a new RHP,
normalized as $z\to\infty$, with the jump matrices ((6.2) of \cite{BDJ})
\begin{equation}
\begin{cases}
v^{(2)}(z;k;t)= \begin{pmatrix} 1& -e^{-2k\wt{\alpha}(z;k;t)}\\
e^{2k\wt{\alpha}(z;k;t)}&0 \end{pmatrix} \ \ &\text{on} \ \ C_2,\\
v^{(2)}(z;k;t)= \begin{pmatrix} e^{-2k\wt{\alpha}-(z;k;t)}& -1\\ 1&0
\end{pmatrix} \ \ &\text{on} \ \ C_1.
\end{cases}
\end{equation}
Through the changes $Y\to m^{(1)}\to m^{(2)}$,
we have
\begin{eqnarray}
\label{as12-91}   Y_{11}(z;k;t) &=& -e^{kg(z;k;t)}m^{(2)}_{12}(z;k;t),
\qquad\qquad\qquad |z|<1,\\
\label{as12-92}   Y_{21}(z;k;t) &=& -(-1)^ke^{kg(z;k;t)+kl}m^{(2)}_{22}(z;k;t),
\qquad |z|<1,\\
\label{as12-93}   Y_{11}(z;k;t) &=& e^{kg(z;k;t)}m^{(2)}_{11}(z;k;t),
\qquad\quad\qquad \qquad |z|>1,\\
\label{as12-94}   Y_{21}(z;k;t) &=& (-1)^ke^{kg(z;k;t)+kl}m^{(2)}_{21}(z;k;t),
\qquad |z|>1.
\end{eqnarray}

\begin{figure}[ht]
 \centerline{\epsfig{file=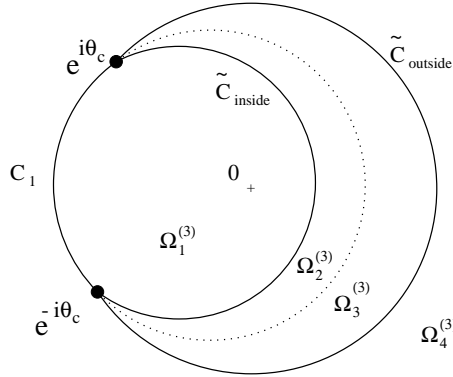, width=6cm}}
 \caption{$\Sigma^{(3)}$ and $\Omega^{(3)}$when $\gamma>1$}
\label{fig-gammage1}
\end{figure}

Set $\Sigma^{(3)}:=\overline{C_1}\cup \wt{C}_{inside} \cup \wt{C}_{outside}$ as
in the Figure \ref{fig-gammage1}, which divides $\C$ into four regions,
$\Omega^{(3)}_j$, $j=1,\cdots,4$. Again, there is a certain freedom choosing
the shape of $\wt{C}_{inside}$ and $\wt{C}_{outside}$. For example,
$\wt{C}_{inside}$ (resp., $\wt{C}_{outside}$) can be any smooth curve lying in
$\Omega^{(3)}_2$ (resp. $\Omega^{(3)}_3$) connecting $e^{i\theta_c}$ and
$e^{-i\theta_c}$; the precise requirement is given in \cite{BDJ} (see also
\eqref{as13.99}--\eqref{as13.101} below). Define $m^{(3)}(z;k;t)$ by (see p.1151
of \cite{BDJ})
\begin{equation}\label{as12-95}
  \begin{cases}
     m^{(3)}=m^{(2)} \begin{pmatrix} 1&-e^{-2k\wt{\alpha}(z;k;t)}\\0&1 \end{pmatrix}^{-1}
\quad&\text{in}\quad \Omega^{(3)}_2,\\
     m^{(3)}=m^{(2)} \begin{pmatrix} 1&0\\e^{2k\wt{\alpha}(z;k;t)}&1 \end{pmatrix}
\qquad&\text{in}\quad \Omega^{(3)}_3,\\
     m^{(3)}=m^{(2)} \qquad\qquad&\text{in}\quad \Omega^{(3)}_1,\Omega^{(3)}_4.
  \end{cases}
\end{equation}
Then $m^{(3)}$ solves a RHP, normalized as $z\to\infty$, with the jump matrix
given by
\begin{equation}
  v^{(3)}(z;k;t)= \begin{cases}
     \begin{pmatrix} 1&-e^{-2k\wt{\alpha}(z;k;t)}\\0&1 \end{pmatrix}
\quad&\text{on}\quad \wt{C}_{inside},\\
     \begin{pmatrix} 1&0\\e^{2k\wt{\alpha}(z;k;t)}&1 \end{pmatrix}
\qquad&\text{on}\quad \wt{C}_{outside},\\
     \begin{pmatrix} e^{-2k\wt{\alpha}-(z;k;t)}&-1\\1&0 \end{pmatrix}
\quad&\text{on}\quad C_1.
  \end{cases}
\end{equation}
From the properties of $g(z)$, it is proved in \cite{BDJ} (between (6.4) and
(6.4)) that
\begin{eqnarray}
\label{as13.99}
   &e^{-k\wt{\alpha}-(z;k;t)}\to 0 \quad \text{as} \quad
k\to\infty, \qquad z\in C_1, \\
\label{as13.100}
   &e^{-k\wt{\alpha}(z;k;t)}\to 0 \quad \text{as} \quad
k\to\infty, \qquad z\in \wt{C}_{inside}, \\
\label{as13.101}
   &e^{k\wt{\alpha}(z;k;t)}\to 0 \quad \text{as} \quad
k\to\infty, \qquad z\in \wt{C}_{outside}.
\end{eqnarray}
The choice of $\wt{C}_{inside}$ and $\wt{C}_{outside}$ is precisely for these
properties. Here the convergence is uniform for any compact part of the each
contour \emph{away} from the end points $e^{i\theta_c}$ and $e^{-i\theta_c}$,
but is not uniform on the whole contour. This gives rise the technical
difficulty which will be overcome below using the idea of parametrix. Formally
$v^{(3)}\to v^{\infty}$ as $k\to\infty$ where
\begin{equation}
\begin{cases}
v^\infty(z) =
 \begin{pmatrix} 1&0\\0&1 \end{pmatrix} \ \ &\text{on} \ \
\wt{C}_{inside}\cup\wt{C}_{outside},\\
v^\infty(z) =
 \begin{pmatrix} 0&-1\\1&0 \end{pmatrix} \ \ &\text{on} \ \
C_1.
\end{cases}
\end{equation}
Thus we expect that $m^{(3)}$ converges to $m^\infty$, the solution of the
RHP $m^\infty_+=m^\infty_-v^\infty$ with $m^\infty\to I$ as $z\to\infty$.
The solution $m^\infty$ is easily given by (see Lemma 6.2 in \cite{BDJ})
\begin{equation}\label{qz10.104}
   m^\infty(z) = \begin{pmatrix} \frac12(\beta+\beta^{-1})&
\frac1{2i}(\beta-\beta^{-1})\\ -\frac1{2i}(\beta-\beta^{-1})&
\frac12(\beta+\beta^{-1}) \end{pmatrix},
\end{equation}
where $\beta(z):=\bigl(\frac{z-e^{i\theta_c}}{z-e^{-i\theta_c}}\bigr)^{1/4}$,
which is analytic $\C\setminus\overline{C}_1$ and $\beta\sim +1$
as $z\in\R\rightarrow +\infty$

\subsubsection{The case $2t\ge ak$ for some $a>1$.}
\label{subsub4}

The parametrix $m_p(z)$ is introduced in (6.25)--(6.31) in \cite{BDJ}, and
has the following properties. In the neighborhood $\mathcal{O}$ of size
$\epsilon$ around the points $e^{i\theta_c}$ and $e^{-i\theta_c}$, $m_p(z)$
is constructed using the Airy function in such a way that
$(m_p(z))_+=(m_p(z))_-v^{(3)}(z)$ for $z\in\Sigma^{(3)}\cap\mathcal{O}$,
and $\|m_p(m^\infty)^{-1} -I\|_{L^\infty(\partial\mathcal{O})}=
O(k^{-1})$. In $\C\setminus\mathcal{O}$, we set $m_p(z):=m^\infty(z)$. Then
the ratio $R(z):=m^{(3)}(z)m_p^{-1}(z)$ has no jump on
$\Sigma^{(3)}\cap\mathcal{O}$, and has a jump $v_R:=w_R+I$ converging to
$I$ uniformly of order $O(k^{-1})$ on $\partial\mathcal{O}$, and of order
$O(e^{-ck})$ on $\Sigma^{(3)}\cap\mathcal{O}^c$ as $k\to\infty$. This
implies that $R(z)=I+O(k^{-1})$ for any $z\in\C\setminus\Sigma_p$,
$\Sigma_p:=(\Sigma^{(3)}\cap\mathcal{O}^c)\cup\partial\mathcal{O}$. Moreover
following the arguments in Section 8 of \cite{DKMVZ2}, the error is uniform
up to the boundary in each open regions in $\C\setminus\Sigma_p$. In
particular, for
$z\in\overline{\Omega^{(3)}_1}\cup\overline{\Omega^{(3)}_4}$, (see
(6.34)-6.40) in \cite{BDJ})
\begin{equation}\label{e-temp}
   m^{(3)}(z)=\bigl( I+O(k^{-1}) \bigr) m^\infty(z).
\end{equation}
Here the error is uniform for $ak\le 2t\le bk$ for some $0<a<b$. For the
case $(2t)/k\to\infty$, by shrinking the size of $\mathcal{O}$ properly, we
again obtain uniform error (see \cite{B99}). Therefore for any $a>0$, we
obtain uniformity in \eqref{e-temp} for $ak\le 2t$.

\medskip
When $z=0$, $\beta(0)=-ie^{i\theta_c/2}$ and $g(0)=\pi i$ (Lemma 4.2 (vi) of
\cite{BDJ}), Also \eqref{as12-95} says $m^{(3)}(0)=m^{(2)}(0)$. Thus
Proposition \ref{thm1} (v) follows from \eqref{as12-91} and \eqref{as12-92}, as
in (6.4) of \cite{BDJ}.

\medskip
Now we consider Proposition \ref{thm4} of the case where $x\le
-2^{1/3}(a-1)k^{2/3}$. For $z=-\alpha$ real, $|\beta(-\alpha)|=1$, so
$m^\infty(-\alpha)$ is bounded. Hence from \eqref{as12-95} and
\eqref{as12-91},\eqref{as12-93}, we have for $\alpha\ge 1$,
\begin{equation}
   |\pi_k(-\alpha;k;t)| \le C|e^{kg(-\alpha;k;t)}|.
\end{equation}
Then we proceed as in \eqref{qz10.109}--\eqref{qz10.121} of the following
subsubsection to obtain the proper estimate.

\subsubsection{The case $k+M2^{-1/3}k^{1/3}\le 2t\le ak$
for some $a>1$ and $M>M_0$.}
\label{subsub5}

In this case, the points $e^{i\theta_c}$ and $e^{-i\theta_c}$
are allowed to approach $-1$, but the rate is restricted:
\begin{equation}\label{as13.107}
   |e^{i\theta_c}+1|=2\biggl(1-\frac{k}{2t}\biggr)^{1/2}
\ge \frac{\sqrt{M}}{k^{1/3}},
\end{equation}
for $k$ large. We now take the neighborhood $\mathcal{O}$ to be of size
$\epsilon\sqrt{\frac{2t}{k}-1}$ around $e^{i\theta_c}$ and
$e^{-i\theta_c}$.  From \eqref{as13.107}, $\mathcal{O}$ consists of two
disjoint disks and their boundaries do not touch the real axis. We
introduce the same parametrix $m_p$ as in the previous subsubsection. Then
we have a similar result: there is $M_0>0$ such that for $M>M_0$,
\begin{equation}
   m^{(3)}(z)= \biggl( I+O\biggl(\frac1{k(\frac{2t}{k}-1)}\biggr) \biggr)
m^\infty(z),
\end{equation}
for $z\in\overline{\Omega^{(3)}_1}\cup\overline{\Omega^{(3)}_4}$; this is
proved in (6.34)--(6.40) of \cite{BDJ}.

\medskip
When $z=0$, as in the previous subsubsection, we obtain Proposition \ref{thm1}
(iv), as in \cite{BDJ}.

\medskip
Now we prove Proposition \ref{thm4} when $x\le -M$.
As in the previous subsubsection, we have
\begin{equation}\label{qz10.109}
   |\pi_k(-\alpha;k;t)| \le
C|e^{kg(-\alpha;k;t)}|.
\end{equation}
Now we need an estimate of $g(-\alpha;k;t)$; as we mentioned at the
beginning of Section \ref{RHP}, this is the second way in which we must
extend \cite{BDJ}. Note that
\begin{equation}\label{qz10.110}
 \begin{split}
   \Re \bigl(g(-\alpha)\bigr)
&=\frac12 \int_{-\theta_c}^{\theta_c}
\log(1+\alpha^2+2\alpha\cos\theta)d\mu(\theta)\\
&= \log 2 +\frac12 \log\frac\alpha\gamma +I(s),
 \end{split}
\end{equation}
where
\begin{equation}\label{as12-110}
   I(s) = \frac1\pi \int_{-1}^1 \log(s^2-x^2)\sqrt{1-x^2}dx,
\qquad s:=\frac{\sqrt\gamma(1+\alpha)}{2\sqrt\alpha}> 1.
\end{equation}
The inequality $s>1$ follows from the arithmetic-geometric mean inequality
and the assumption $\gamma>1$.
A residue calculation gives us
\begin{equation}
  I'(y)= \frac1\pi\int_{-1}^1 \frac{2y}{y^2-x^2}\sqrt{1-x^2}dx
= 2y-2\sqrt{y^2-1}, \qquad y>1.
\end{equation}
Integrating from $1$ to $s>1$, we have
\begin{equation}
  I(s) =s^2-1-2\int_1^s \sqrt{y^2-1}dy +I(1).
\end{equation}
The constant $I(1)$ can be evaluated (cf. Lemma 4.3 (ii)-(a) in
\cite{BDJ}):
\begin{equation}
   I(1)= \frac1\pi \int_{-\frac\pi2}^{\frac\pi2}
\log(\sin^2\theta) \sin^2\theta d\theta
= \frac12 -\log 2.
\end{equation}
Thus we have
\begin{equation}\label{as13.113}
  \Re \bigl(g(-\alpha)\bigr) = -\frac12 +\frac12 \log\frac\alpha\gamma
+s^2 -s\sqrt{s^2-1}+\log(s+\sqrt{s^2-1}).
\end{equation}

Assume $0<\alpha\le 1$.
We change the variables $\gamma$, $\alpha$ into
$s$, $\xi$ where $s$ is defined in \eqref{as12-110} and
\begin{equation}
   \xi:= \biggl(\frac\gamma\alpha\biggr)^{1/2} > 1.
\end{equation}
Then
\begin{equation}
   F(\xi):= g(-\alpha)-\frac\gamma2\alpha =
-\frac12-\log\xi -\frac12\xi^2+2s\xi-s^2
-s\sqrt{s^2-1}+\log(s+\sqrt{s^2-1}).
\end{equation}
Differentiating with respect to $\xi$, we find
\begin{equation}
  F'(\xi)= -\frac1\xi-\xi+2s.
\end{equation}
Thus the maximum of $F$ occurs at $\xi=s+\sqrt{s^2-1}$.
But $F(s+\sqrt{s^2-1})=0$, hence $F(\xi)\le 0$.
Thus we obtain
\begin{equation}\label{as13.117}
   |e^{-t\alpha}\pi_k(-\alpha;k)| \le
Ce^{k\Re\bigl(g(-\alpha;k;t)-\frac\gamma2 \alpha\bigr)} \le C,
\qquad 0<\alpha\le 1.
\end{equation}
For $\alpha\ge 1$, note that
\begin{equation}
  \Re(g(-\alpha))=\log \alpha+\Re(g(-\alpha^{-1})).
\end{equation}
Thus using \eqref{as13.117},  we have
\begin{equation}\label{qz10.121}
   |e^{-t\alpha^{-1}}(-\alpha)^{-k}\pi_k(-\alpha;k)| \le
Ce^{k\Re\bigl(g(-\alpha;k;t)-\frac\gamma2\alpha^{-1}-\log\alpha\bigr)}
= Ce^{k\Re\bigl(g(-\alpha^{-1};k;t)-\frac\gamma2\alpha^{-1}\bigr)}
\le C.
\end{equation}

\subsubsection{Case $k<2t\le k+M2^{-1/3}k^{1/3}$
for some $M>0$.}
\label{subsub6}

\begin{figure}[ht]
 \centerline{\epsfig{file=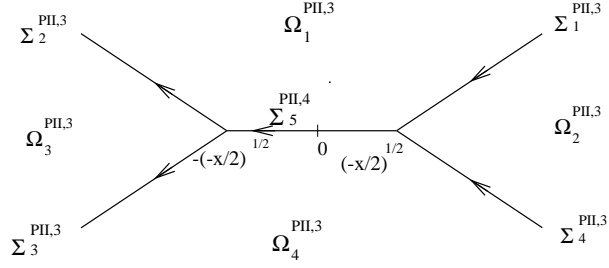, width=8cm}}
 \caption{$\Sigma^{PII,3}$ and $\Omega^{PII,3}_j$}\label{fig-PII3}
\end{figure}

First we introduce $m^{PII,3}$ as in (2.22)--(2.28) of \cite{BDJ}. Set
\begin{equation}
  g^{PII}(z):= \frac43\bigl(z^2+\frac{x}2\bigr)^{3/2},
\end{equation}
which is analytic in $\C\setminus\bigl[-\sqrt{\frac{-x}2},
\sqrt{\frac{-x}2}\bigr]$
and behaves like $\frac43z^3+xz + \frac{x^2}{8z}+O(z^{-3})
=:\theta_{PII}(z)+O(z^{-1})$
as $z\to+\infty$.
Let $\Sigma^{PII,3}:=\cup_{j=1}^5 \Sigma^{PII,3}_j$ as shown
in Figure \ref{fig-PII3}.
The angles of the rays with the real line are between $0$ and $\pi/3$.
Recall that $m(z;x)$ solves \eqref{as20} the RHP for PII equation.
Define $m^{PII,3}(z;x)$ by
\begin{equation}\label{as12-103}
 \begin{cases}
   m^{PII,3}= m(z;x) e^{i(g^{PII}-\theta_{PII})\sigma_3},
&z\in\Omega^{PII,3}_1, \Omega^{PII,3}_4, \\
   m^{PII,3}= m(z;x)\begin{pmatrix} 1&e^{-2i\theta_{PII}}\\0&1
\end{pmatrix}
e^{i(g^{PII}-\theta_{PII})\sigma_3},
&z\in (\Omega^{PII,3}_2\cup\Omega^{PII,3}_3)\cap\C_-, \\
   m^{PII,3}= m(z;x)\begin{pmatrix} 1&0\\e^{2i\theta_{PII}}&1
 \end{pmatrix}
e^{i(g^{PII}-\theta_{PII})\sigma_3},
&z\in (\Omega^{PII,3}_2\cup\Omega^{PII,3}_3)\cap\C_+.
 \end{cases}
\end{equation}
Then $m^{(3)}$ solves the RHP (see (2.25) in \cite{BDJ})
normalized at $\infty$ with the jump matrix
\begin{equation}
 v^{(3)}(z;k;t)=
  \begin{cases}
     \begin{pmatrix}1&0\\e^{2ig^{PII}}&1 \end{pmatrix}
\qquad&\text{on}\quad \Sigma_1^{PII,3},\Sigma_2^{PII,3}\\
     \begin{pmatrix}1&-e^{-2ig^{PII}}\\0&1 \end{pmatrix}
\quad&\text{on}\quad \Sigma_3^{PII,3},\Sigma_4^{PII,3}\\
     \begin{pmatrix}e^{-2ig_{-}^{PII}}&-1\\1&0 \end{pmatrix}
\quad&\text{on}\quad \Sigma_5^{PII,3}.\\
  \end{cases}
\end{equation}
Also we have
\begin{equation}\label{as12-105}
   m_1(x)=m_1^{PII,3}(x)-\bigl( \frac{ix^2}{8} \bigr)\sigma_3.
\end{equation}
where $m^{PII,3}(z;x)= I+\frac{m^{PII,3}_1(x)}{z}+O(z^{-2})$
as $z\to\infty$.

Now as before, set $x$ by
\begin{equation}
   \frac{2t}{k}=1-\frac{x}{2^{1/3}k^{2/3}}.
\end{equation}
Hence we have $-M\le x<0$ in this subsubsection. We now proceed as in case
(iii) of Section 6 in \cite{BDJ}.  Define the parametrix
\begin{equation}
  \begin{cases}
    m_p(z;k;t)=m^{PII,3}\bigl(\lambda(z),x\bigr)&
\quad\text{in} \quad\mathcal{O}\setminus\Sigma^{(3)},\\
    m_p(z;k;t)= I& \quad\text{in}\quad\bar{\mathcal{O}}^c\setminus\Sigma^{(3)},
  \end{cases}
\end{equation}
where $\lambda(z)$ is defined in \eqref{as12.59.5} and
$\mathcal{O}$ is a small neighborhood of size $\epsilon>0$ around $z=-1$
(see \cite{BDJ} case (iii) of Section 6 for details).
As in Subsubsection \ref{subsub3}, the ratio
$R(z;k;t):= m^{(3)}m_p^{-1}$ satisfies a new RHP, normalized at $\infty$,
with jump matrix $v_R$ satsifying the estimate
\eqref{as12.62} where $m^{PII,2}_1(x)$ is replaced by $m^{PII,3}_1(x)$.
Hence we have
\begin{equation}
   \biggl| m^{(3)}(z;k;t)\bigl(m_p(z;k;t)\bigr)^{-1}-I
+ \frac{m^{PII,3}_1(x)}{2\pi i} \int_{\partial\mathcal{O}}
\frac1{\lambda(s)(s-z)} ds \biggr| \le \frac{\epsilon
C}{\dist(z,\Sigma_R)k^{2/3}},
\end{equation}
which is hidden in the derivation of (6.19) of \cite{BDJ}.

\medskip
Then as in \eqref{as12-71}, we have
\begin{equation}\label{as12-108}
    m^{(3)}(0;k;t)
= I+ \frac{i2^{4/3}m_1^{PII,3}(x)}{k^{1/3}} + O\bigl( \frac1{k^{2/3}}\bigr),
\end{equation}
which is a (direct) extension of (6.19) in \cite{BDJ}. Now from \eqref{as12-95}
and \eqref{as12-105} (see (6.19) of \cite{BDJ}),
\begin{equation}
   m^{(2)}(0;k;t)= I+\frac{i2^{4/3}m_1(x)-2^{-5/3}x^2\sigma_3}{k^{1/3}}
+ O\bigl(\frac1{k^{2/3}}\bigr).
\end{equation}
Hence using $e^{kl}=1-\frac{x^2}{2^{5/3}k^{1/3}}+O(k^{-2/3})$ and $g(0)=\pi i$,
\eqref{as12-91} and \eqref{as12-92} yield Proposition \ref{thm1} (iii) in the
case $-M\le x<0$, as in \cite{BDJ}.

\medskip
For the proof of Proposition \ref{thm2}, note that, as before, for each fixed
$z\in\C\setminus\Sigma$, we can use the freedom of the shape of $\Sigma^{(3)}$
(and $\Sigma_R$) so that $z\in\Omega^{(3)}_1\cup\Omega^{(3)}_4$, and
$\dist(z,\Sigma_R)\ge c_1>0$. Thus we obtain
\begin{equation}\label{as13.128}
  \lim_{k\to\infty} m^{(2)}(z;k;t)= I,
\qquad z\in\C\setminus\Sigma \ \ fixed.
\end{equation}
From \eqref{as12-91} and \eqref{as12-93}, we have
\begin{eqnarray}
  \lim_{k\to\infty} e^{-kg(z;k;t)}\pi_k(z)=0,
\qquad &|z|<1,\\
  \lim_{k\to\infty} e^{-kg(z;k;t)}\pi_k(z)=1,
\qquad &|z|>1.
\end{eqnarray}
This is an extension of the calculation of \cite{BDJ} where $z=0$ is given. Now
in order to prove Proposition \ref{thm2}, we need further analysis which is an
extension in the second category as mentioned as the beginning of Section
\ref{RHP}. Since $\gamma=(2t)/k$, for the proof of Proposition \ref{thm2}, it
is enough to show that
\begin{eqnarray}
\label{as12-115}
\lim_{k\to\infty} (-1)^ke^{k[g(z;k;t)+\frac{\gamma}2z]}=1,
\qquad &|z|<1,\\
\label{as12-116}
\lim_{k\to\infty} e^{k[g(z;k;t)+\frac{\gamma}2z^{-1}-\log z]}=1,
\qquad &|z|>1.
\end{eqnarray}
But the proof of Lemma 4.3 (ii) of \cite{BDJ} says that for $|z|>1$,
$z\notin(-\infty,-1)$,
\begin{equation}
  g(z)= \frac12 \log z -\frac{\gamma}4(z+z^{-1})+\frac\gamma2
+\frac\gamma4 \int_{1_{+0}}^z \frac{s+1}{s^2}
\sqrt{(s-e^{i\theta_c})(s-e^{-i\theta_c})} ds +g_-(1),
\end{equation}
where the integral is taken over a curve from $1_{+0}$ to $z$
lying in $\{z\in\C : |z|>1, z\notin(-\infty,-1)\}$.
Here $\sqrt{(s-e^{i\theta_c})(s-e^{-i\theta_c})}$
is analytic in $\C\setminus C_1$ and behaves like $s$ as $s\in\R\to+\infty$,
and $\log z$ is analytic in $\C\setminus(-\infty,0]$ and is real
for $z\R_+$.
Calculations in the same proof, together with Lemma 4.2 (viii) of \cite{BDJ},
give us $g_-(1)=-\frac12-\frac12\log\gamma$.
Also using $\sin^2\frac{\theta_c}2=\frac1\gamma$,
for $|s|>1$, $s\notin(-\infty,-1)$,
\begin{equation}
   \sqrt{(s-e^{i\theta_c})(s-e^{-i\theta_c})}
=(s+1)-\frac{2s}{s+1}(\gamma-1)+O((\gamma-1)^2).
\end{equation}
Thus expanding in $\gamma-1$, we have
\begin{equation}
   g(z)+\frac{\gamma}2z^{-1}-\log z
=O((\gamma-1)^2) = O(k^{-4/3}),
\end{equation}
which implies \eqref{as12-116}.
Similar calculations using
for $|z|<1$, $z\notin(-1,0]$,
\begin{equation}\label{as13.136}
   g(z)= \frac12 \log z -\frac{\gamma}4(z+z^{-1})+\frac\gamma2
+\frac\gamma4 \int_{1_{+0}}^z \frac{s+1}{s^2}
\sqrt{(s-e^{i\theta_c})(s-e^{-i\theta_c})} ds +g_+(1),
\end{equation}
and $g_+(1)=-\frac12-\frac12\log z+\pi i$ yield \eqref{as12-115}.

\medskip
For the proof of Proposition \ref{thm3} when $x<0$, set
\begin{equation}\label{as12-109}
  \alpha =1-\frac{2^{4/3}w}{k^{1/3}}.
\end{equation}
In this case, we need a new argument as $\alpha\to 1$. When $w$ and $x$ are
fixed, again as in \eqref{as12-76}, we have
$\lim_{k\to\infty}R(-\alpha;k;t)=I$, which implies that
\begin{equation}
   \lim_{k\to\infty} m^{(3)}(-\alpha;k;t) = m^{PII,3}(-iw,x).
\end{equation}
From (6.8) of \cite{BDJ}, we have
\begin{equation}
   \lim_{k\to\infty} k\wt{\alpha}(-\alpha;k;t)
= i\frac43\bigl((-iw)^2+\frac{x}2\bigr)^{3/2}
=ig^{PII}(-iw),
\end{equation}
which from \eqref{as12-95} and \eqref{as12-103} implies
\begin{equation}\label{as12.115}
  \lim_{k\to\infty} m^{(2)}(-\alpha;k;t)
= m(-iw,x)e^{i(g^{PII}(-iw)-\theta_{PII}(-iw))\sigma_3}.
\end{equation}
Now we compute the large $k$ limit of
$kg(-\alpha;k;t)-t\alpha$ when $w>0$,
and of $kg(-\alpha;k;t)-t\alpha^{-1}-\log\alpha$ when $w<0$.
For $-\pi<\theta<\pi$,
$\lim_{\epsilon\downarrow 0} \argum(-\alpha+i\epsilon-e^{i\theta})
=\pi+\tan^{-1}(\frac{\sin\theta}{\alpha+\cos\theta})$
where $-\pi<\tan^{-1}\phi<\pi$.
Since $\tan^{-1}(\frac{\sin\theta}{\alpha+\cos\theta})$
is odd in $\theta$, we have from \eqref{as13.113},
\begin{equation}\label{as12-131}
 \begin{split}
   \lim_{\epsilon\downarrow 0} g(-\alpha+i\epsilon)
&= \lim_{\epsilon\downarrow 0} \Re g(-\alpha+i\epsilon) +\pi i\\
&= -\frac12 +\frac12\log\frac\alpha\gamma
+ s^2-s\sqrt{s^2-1} + \log(s+\sqrt{s^2-1})
+\pi i
 \end{split}
\end{equation}
where $s=\frac{\sqrt\gamma(1+\alpha)}{2\sqrt\alpha}>1$.
Under the stated conditions on $\gamma$ and $\alpha$, as $k\to\infty$,
\begin{eqnarray}
   \frac{\sqrt\gamma(1+\alpha)}{2\sqrt\alpha}
&=&1+\biggl(\frac{w^2}{2^{1/3}}-\frac{x}{2^{4/3}} \biggr)
\cdot\frac1{k^{2/3}}
+\frac{2w^3}{k} + O\bigl(\frac1{k^{-4/3}}\bigr),\\
\frac\alpha\gamma&=&1-\frac{2^{4/3}w}{k^{1/3}}
+\frac{x}{2^{1/3}k^{2/3}}-\frac{2xw}{k}
+O\bigl(\frac1{k^{-4/3}}\bigr).
\end{eqnarray}
Note that $\lim_{\epsilon\downarrow 0}g(-\alpha+i\epsilon)
-\lim_{\epsilon\downarrow 0}g(-\alpha-i\epsilon)$
is $2\pi i$ for $\alpha>1$, and is $0$ for $0<\alpha<1$.
Therefore we obtain
\begin{eqnarray}
   \lim_{k\to\infty} (-1)^ke^{kg(-\alpha)-t\alpha}&=&
e^{\frac43w^3-xw-\frac43(w^2-\frac{x}2)^{3/2}},
\qquad w>0,\\
   \lim_{k\to\infty} (-1)^ke^{kg(-\alpha)-t\alpha^{-1}-k\log\alpha}&=&
e^{-\frac43w^3+xw-\frac43(w^2-\frac{x}2)^{3/2}},
\qquad w<0.
\end{eqnarray}
Also being careful of branch, we have
\begin{eqnarray}
\label{as13.152}
   i(g^{PII}(-iw)-\theta_{PII}(-iw))
&=& \frac43w^3-xw-\frac43\biggl(w^2-\frac{x}2\biggr)^{3/2},
\qquad w>0,\\
\label{as13.153}
   i(g^{PII}(-iw)-\theta_{PII}(-iw))
&=& \frac43w^3-xw+\frac43\biggl(w^2-\frac{x}2\biggr)^{3/2},
\qquad w<0.
\end{eqnarray}
Since $\lim_{k\to\infty}e^{kl}=1$, using \eqref{as12-91}--\eqref{as12-94}
and \eqref{as12.115}, this implies \eqref{as7.32}--\eqref{as7.35}
when $x<0$.

%
%

\subsection{Proof of Proposition \ref{prop7.7}}\label{sub3}

The analysis in this section is new, and is needed for the proof of Proposition
\ref{prop7.7} . Let $\alpha>1$ be fixed. Set
\begin{equation}\label{as12.88}
   \frac{t}{k} = \frac{\alpha}{\alpha^2+1}
-\frac{\alpha(\alpha^2-1)^{1/2}}{(\alpha^2+1)^{3/2}}
\cdot\frac{x}{\sqrt{k}}, \qquad \text{$x\in\R\setminus\{0\}$ is fixed.}
\end{equation}
We are interested in the asymptotics of
$e^{-\alpha t}(-\alpha)^l\pi_k(-\alpha^{-1};t)$.
Since $\frac{2\alpha}{\alpha^2+1}<1$
and $x$ is fixed, we are in the case of
Subsubsection \ref{subsub1} and/or \ref{subsub2}.
We define $m^{(1)}$ and $m^{(2)}$ as in
\eqref{as12.1} and \eqref{as12.4}.
Recall that we have certain freedom in the choice of
$\Sigma^{(3)}$.
We choose a contour passing through the saddle points
of
\begin{equation}
  f(z):=\frac{t}{k}(z-z^{-1})+\log (-z),
\end{equation}
the exponent of the $(12)$ entry of $v^{(2)}$
divided by $k$.
The saddle points (see \eqref{as12.28}
and the following discussion) are $-\rho_\pi$
and $-\rho_\pi^{-1}$, where
\begin{equation}\label{as12.88.5}
   \rho_\pi=
\frac{1-\sqrt{1-(2t/k)^2}}{2t/k}
= \frac1{\alpha} -
\frac{(\alpha^2+1)^{1/2}}{\alpha(\alpha^2-1)^{1/2}}
\cdot\frac{x}{\sqrt{k}} + O\bigl(\frac1{k}\bigr)
=: \rho_c+ O\bigl(\frac1{k}\bigr).
\end{equation}
Take $\delta>0$ and $\epsilon>0$ small such that
$\Sigma_c:=\{-\rho_c+is :
-k^{\delta-1/2}\le s\le\epsilon k^{\delta-1/2}\}$
lies inside the open unit disc for all $k\ge 1$.
\begin{figure}[ht]
 \centerline{\epsfig{file=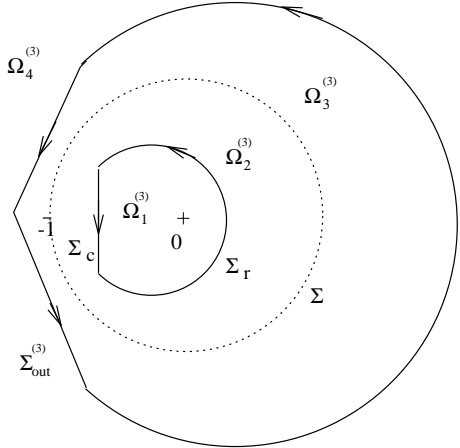, width=6cm}}
 \caption{$\Sigma^{(3)}$ and $\Omega^{(3)}$}
\label{fig-alpha}
\end{figure}
Define (see Figure \ref{fig-alpha})
$\Sigma^{(3)}:=\Sigma_{in}^{(3)}\cup\Sigma^{(3)}_{out}$
by $\Sigma_{in}^{(3)} := \Sigma_c\cup
\Sigma_r$
and $\Sigma_{out}^{(3)} := \{ r^{-1}e^{i\phi} :
re^{i\phi}\in\Sigma^{(3)}_{in} \}$,
where
$\Sigma_r:= \{ |-\rho_c+i\epsilon k^{\delta-1/2}|e^{i\theta}
: |\theta|<\theta_0\}$,
$-\rho_c+i\epsilon k^{\delta-1/2}=
|-\rho_c+i\epsilon k^{\delta-1/2}|e^{i\theta_0}$.
Let $\Omega^{(3)}_j$ be as in Figure \ref{fig-alpha}
and define $m^{(3)}$ as in \eqref{as12.31}.
As in \eqref{as12.42}-\eqref{as12.44},
the quantities we are interested in are
\begin{eqnarray}
\label{as12.89}
   e^{-\alpha t}(-\alpha)^k\pi_k(-\alpha^{-1};t)
&=& -\alpha^ke^{-t(\alpha-\alpha^{-1})}
m^{(3)}_{12}(-\alpha^{-1};k;t),
\quad\qquad\qquad\qquad\qquad
x<0,\\
\label{as12.90}
   e^{-\alpha t}(-\alpha)^k\pi_k(-\alpha^{-1};t)
&=& -\alpha^ke^{-t(\alpha-\alpha^{-1})}
m^{(3)}_{12}(-\alpha^{-1};k;t)
+ m^{(3)}_{11}(-\alpha^{-1};k;t),
\quad x>0.
\end{eqnarray}

For the estimates of $w^{(3)}:=v^{(3)}-I$,
note that for any fixed $0<\rho<1$,
$\Re f(\rho e^{i\theta})=F(\rho e^{i\theta};\frac{2t}{k})$
(recall \eqref{as13.11})
is increasing
in $0<\theta<\pi$, and is decreasing
in $\pi<\theta<2\pi$, hence
$\|e^{kf(z)}\|_{L^\infty(\Sigma_r)}
= e^{k\Re f(-\rho_c+\epsilon k^{\delta-1/2})}$.
But we have
\begin{equation}\label{as12.93}
   f(-\rho_c+ia)
= \frac{t}{k}(\alpha-\alpha^{-1})-\log\alpha
-\frac12\biggl(
\frac{\alpha^2(\alpha^2-1)}{\alpha^2+1}a^2
+ \frac{x^2}{k}\biggr)
+ O(k^{-\frac32+2\delta}),
\qquad |a|\le \epsilon k^{\delta-1/2}.
\end{equation}
Thus,
\begin{equation}\label{as12.94}
\|\alpha^ke^{-t(\alpha-\alpha^{-1})}e^{kf(z)}
\|_{L^\infty(\Sigma_r)}
\le Ce^{-ck^{2\delta}},
\end{equation}
and also using
\begin{equation}\label{as12.95}
   \alpha^{-k}e^{t(\alpha-\alpha^{-1})}
=e^{k\bigl[\log\alpha-\frac{\alpha^2-1}{\alpha^2+1}
+ O\bigl(\frac1{\sqrt{k}}\bigr)\bigr]},
\qquad \log\alpha-\frac{\alpha^2-1}{\alpha^2+1} >0,
\ \ \text{for $\alpha>1$,}
\end{equation}
we have
\begin{equation}
   \|e^{kf(z)}\|_{L^\infty(\Sigma_r)}
\le Ce^{-ck}.
\end{equation}
On the other hand, one can directly check that
$\Re f(-\rho_c+a)$ has its maximum at $a=0$ for
$-\epsilon k^{\delta-1/2}\le a\le \epsilon k^{\delta-1/2}$, hence
$\|e^{kf(z)}\|_{L^\infty(\Sigma_c)}
= e^{k \Re f(-\rho_c)}$.
Again \eqref{as12.93} and \eqref{as12.95} yield
\begin{equation}\label{as12.95.5}
   \|e^{kf(z)}\|_{L^\infty(\Sigma_c)} \le e^{-ck}.
\end{equation}
Similarly, we have
$\|e^{-kf(z)}\|_{L^\infty(\Sigma^{(3)}_{out})} \le e^{-ck}$.
Now calculations as in Subsubsection \ref{subsub2}
give us the result \eqref{as12.39}.
Hence using
\eqref{as12.95} and \eqref{as12.95.5}
and noting $\dist(-\alpha^{-1},\Sigma^{(3)})
=\frac{(\alpha^2+1)^{1/2}}{\alpha(\alpha^2-1)^{1/2}}
\cdot\frac{x}{\sqrt{k}}$, we have
\begin{eqnarray}
\label{as12.91}
   m_{11}^{(3)}(-\alpha^{-1};k;t) &=&
1+ O(\sqrt{k}e^{-2(1-\epsilon_1)c_1k}),\\
\label{as12.92}
   \alpha^ke^{-t(\alpha-\alpha^{-1})}
m_{12}^{(3)}(-\alpha^{-1};k;t) &=&
\alpha^ke^{-t(\alpha-\alpha^{-1})}
\int_{\Sigma_c}
\frac{(-s)^ke^{t(s-s^{-1})}}{s+\alpha^{-1}}
\frac{ds}{2\pi i}
+ O(\sqrt{k}e^{-ck^{2\delta}}).
\end{eqnarray}

To evaluate the integral asymptotically, first we change the variable by
$s=-\rho_c- \frac{i(\alpha^2+1)^{1/2}}{\alpha(\alpha^2-1)^{1/2}}
\cdot\frac{y}{\sqrt{k}}$. Then from \eqref{as12.93}, the numerator of the
integrand becomes
\begin{equation}
   \alpha^{-k}e^{t(\alpha-\alpha^{-1})}
e^{-\frac12(y^2+x^2)+O(k^{-\frac12+2\delta})}.
\end{equation}
Hence setting
$A:=\frac{\epsilon\alpha(\alpha^2-1)^{1/2}}
{(\alpha^2+1)^{1/2}}$,
\begin{equation}
  \alpha^ke^{-t(\alpha-\alpha^{-1})}
\frac{1}{2\pi i}\int_{\Sigma_c} \frac{(-s)^ke^{t(s-s^{-1})}}{s+\alpha^{-1}} ds
= \frac1{2\pi i} \int_{-A k^\delta}^{A k^\delta}
\frac{e^{-\frac12(y^2+x^2)}}{y+ix} dy \bigl(1+O(k^{-\frac12+2\delta})\bigr) +
O(e^{-ck^{2\delta}}).
\end{equation}
Thus from \eqref{as12.89} and \eqref{as12.90},
we obtain
\begin{eqnarray}
\label{as12.101}
   \lim_{k\to\infty}
e^{-\alpha t}(-\alpha)^k\pi_k(-\alpha^{-1};t)
&=&
\frac1{2\pi i} \int_{-\infty}^{\infty}
\frac{e^{-\frac12(y^2+x^2)}}{y+ix} dy,
\qquad\qquad
x<0,\\
\label{as12.102}
   \lim_{k\to\infty}
e^{-\alpha t}(-\alpha)^k\pi_k(-\alpha^{-1};t)
&=&
\frac1{2\pi i} \int_{-\infty}^{\infty}
\frac{e^{-\frac12(y^2+x^2)}}{y+ix} dy
+1,
\qquad x>0.
\end{eqnarray}
The function
$h(x):= \frac1{2\pi i} \int_{-\infty}^{\infty}
\frac{e^{-\frac12(y^2+x^2)}}{y+ix} dy$ is
smooth in $x>0$ and $x<0$.
The derivative is $h'(x)=\frac1{\sqrt{2\pi}}e^{-\frac12x^2}$.
As $x\to\pm\infty$, $h(x)\to 0$.
Therefore we see that
\begin{equation}
   \frac1{2\pi i} \int_{-\infty}^{\infty}
\frac{e^{-\frac12(y^2+x^2)}}{y+ix} dy
=
\begin{cases}
  \frac1{\sqrt{2\pi}}\int_{-\infty}^x e^{-\frac12y^2}dy,
\qquad &x<0,\\
  \frac1{\sqrt{2\pi}}\int_\infty^x e^{-\frac12y^2}dy,
\qquad &x>0.
\end{cases}
\end{equation}
Thus we proved Proposition \ref{prop7.7}.

\bibliographystyle{plain}
\bibliography{part2}
\end{document}